\documentclass[12pt,a4paper]{amsart}

\usepackage{lscape}
\usepackage{graphicx}
\usepackage[all]{xy}

\usepackage{DLdef}

\newcommand{\uN}{\underline{N}}

\newcommand {\fbr}{{\mathfrak{br}}}
\newcommand {\fbrj}{{\mathfrak{brj}}}
\newcommand {\fbgl}{{\mathfrak{bgl}}}

\newcommand {\fel}{{\mathfrak{el}}}

\newcommand {\fwk}{{\mathfrak{wk}}}

\newcommand {\foo}{{\mathfrak{oo}}}

\newcommand {\size}{{\text{size}}}

\newcommand{\fhov}{\overline{\fh}}
\newcommand{\nc}{_2(\fn;\Kee)}
\newcommand{\ww}{\wedge}

\rmnameii{Charr}{char}

\newcounter{notes}
\newenvironment{arab}{%
 \begin{list}{\arabic{notes})}{\usecounter{notes}%
  \settowidth{\labelwidth}{0.0em}
  \setlength{\labelsep}{0.0em}
  \setlength{\leftmargin}{0em}
  \setlength{\itemindent}{0.0em}
  \setlength{\itemsep}{1mm}
  \setlength{\topsep}{2mm}
  \setlength{\parsep}{3mm}
  \setlength{\partopsep}{0mm}
  }}
 {\end{list}}



\begin{document}

\graphicspath{{./figs/}}

\title[Divided power (co)homology. Presentations]
{Divided power (co)homology. Presentations of simple finite
dimensional modular Lie superalgebras with Cartan matrix}

\author{Sofiane Bouarroudj}
\email{Bouarroudj.sofiane@uaeu.ac.ae}
\address{Department of Mathematics, United Arab Emirates University, Al
Ain, PO. Box: 17551}

\author{Pavel Grozman}
\email{pavel.grozman@bredband.net}
\address{Equa Simulation AB, Stockholm, Sweden}
\author{Alexei Lebedev}
 \email{yorool@mail.ru}
 \address{Nizhegorodskij Univ., RU-603950, N. Novgorod, pr. Gagarina, 23,
 Russia}
\author{Dimitry Leites}
\email{mleites@math.su.se}

\address{Department of Mathematics,
Stockholm University, Roslagsv. 101, Kr\"aft\-riket hus 6, SE-106 91
Stockholm, Sweden}




\begin{abstract}  For modular Lie superalgebras, new notions are
introduced: Divided power homology and divided power cohomology. For
illustration, we give presentations (in terms of analogs of
Chevalley generators) of finite dimensional  Lie (super)algebras
with indecomposable Cartan matrix in characteristic 2 (and in other
characteristics for completeness of the picture).

We correct the currently available in the literature  notions of
Chevalley generators and Cartan matrix in the modular and super
cases, and an auxiliary notion of the Dynkin diagram.

In characteristic 2, the defining relations of simple classical Lie
algebras of the A, D, E types are not only Serre ones; these
non-Serre relations are same for Lie superalgebras with the same
Cartan matrix and any distribution of parities of the generators.
Presentations of simple orthogonal Lie algebras having no Cartan
matrix are also given.

\end{abstract}

\maketitle

\section{Introduction}
In what follows $\Kee$ is a field of characteristic $p>0$,
algebraically closed unless otherwise stated. The Lie
(super)algebras considered are of finite dimension.

\ssec{Divided power (co)homology} Over $\Kee$, the notion of Lie
{\bf super}algebra (co)homology obtains one more dimension --- the
shearing parameter $\underline{N}$. Indeed, since the (co)chains,
with trivial coefficients and differential forgotten, form a
supercommutative superalgebra
--- an analog of the polynomial superalgebra (with values in a
module for non-trivial coefficients), and the polynomial algebra has
divided power analogs in the modular case, so does Lie {\bf
super}algebra (co)homology. For Lie algebras, this phenomenon does
not exists since the supercommutative superalgebra of polynomials is
generated by purely odd elements only.

This being the main idea, the only thing to do is to define the
differential. The appropriate definitions are given in the text and
even implemented in the package \texttt{SuperLie}, see \cite{Gr}.

For an illustration, we consider defining relations (here),
deformations of (finite dimensional) Lie superalgebras with
indecomposable Cartan matrix and of queer series (in \cite{BGL2});
in the sequels (in preparation) we consider deformations of their
representations. In {\bf these} problems, the effect of divided
power (co)homology is only visible for $p=2$; for completeness,
however, we consider presentations for all new (previously not
covered in the literature) cases for $p>2$.

\ssec{Presentations of simple Lie (super)algebras: Overview}{}~{}

$\bullet$ {\bf Over $\Cee$}, the most studied type of simple Lie
algebras are finite dimensional ones and the $\Zee$-graded of
polynomial growth. The latter type splits into (twisted) loop
algebras, vectorial Lie algebras (with polynomial coefficients) and
Witt algebra (the vectorial Lie algebra with Laurent polynomials as
coefficients).

The finite dimensional and (twisted) loop algebras can be defined by
means of their Cartan matrix and Chevalley generators (we recall
these notions in what follows). The explicit presentation was first
published by Serre and certain relations, whose sufficiency was most
difficult to prove, are referred to as {\it Serre relations}.

For simple vectorial Lie algebras, it was not even clear (until
implicitly by V.~Uf\-na\-rov\-sky in late 1970s for some cases) if
they were finitely presented; for the explicit presentations
eventually obtained, and references, see \cite{GLP}, where the Lie
superalgebras are also considered. The relations are passable for a
computer, but rather ugly for humans; the only message we can deduce
from their description at the moment is that, in addition to the
relations in the linear part of the vectorial Lie algebra, there are
only finitely many relations (less than 10 for any type of the
algebras).

In the super case, in addition to the Lie superalgebras of vectorial
type and those with Cartan matrix, there are also the queer series
whose presentation is clear in principle, but whose explicit form is
even less appealing than that of vectorial Lie (super)algebras, see
\cite{LSe}.

$\bullet$ {\bf Over $\Kee$}, the vectorial Lie (super)algebras
acquire one more parameter (a shearing vector $\underline{N}$, see
\eqref{u;N}) and even the description of generators becomes too
complicated, to say nothing of relations. For the {\it restricted
case} and {\it sufficiently large characteristic} and dimension of
the space on which the vectorial Lie algebra is realized, the answer
is identical to that obtained in \cite{GLP}.

The Lie (super)algebra with more roots of one sign (say, positive
than negative) is said to be {\it skew-symmetric} and {\it
symmetric} otherwise.\footnote{For $p=2$, there are skew-symmetric
Lie (super)algebras distinct from vectorial Lie (super)algebras.}
For vectorial Lie (super)algebras, as well as for queer Lie
superalgebras, the general picture of their presentations is clear
and as long as the explicit answer is not really needed (as in
\cite{LSg}, where somewhat awful relations found in \cite{GL1} are
used), we see no point in deriving it. In every particular case, it
is easy by means of \texttt{SuperLie} \cite{Gr} to anybody capable
to use {\it Mathematica}.

For $p=2$, several more types of simple Lie (super)algebras appear:
Symmetric but without Cartan matrix (such as $\fo_I(n)$ and
$\fq(\fo_I(n))$, see \cite{LeD}), various deforms of the
above-listed types. So, we arrive at the last cases left:
\begin{enumerate}
  \item Lie (super)algebras with Cartan
matrix;
  \item Lie (super)algebras without Cartan
matrix but not of vectorial type.
\end{enumerate}
In this paper, we consider case 1 (and a series of examples of case
2: The Lie algebras\footnote{The derived of $\fg(A)$ (or any other
algebra with a commentary in parentheses like $(A)$ after a \lq\lq
family name" $\fg$) should be denoted $\fg(A)^{(i)}$ but it is
usually more convenient to denote it $\fg^{(i)}(A)$ (and similarly
treat other commentaries).} $\fo^{(1)}_I(2n)$). The first thing to
do is to define the basic notions sufficiently clear: Unlike humans,
computers can not work otherwise whereas we can not write this text
without computer's assistance.

\ssec{Lie superalgebras with Cartan matrix} The classification of
{\it finite dimensional modular Lie algebras with indecomposable
Cartan matrix} over  $\Kee$ was obtained in \cite{WK} with a gap
corrected in \cite{Br3, Sk1} (not even mentioned in \cite{KWK}).
Although in \cite{WK} some notions used in the description of the
classification were left undefined, the strategy was impeccable. In
\cite{BGL}, we clarified the notions left somewhat vague during the
time elapsed since publication of \cite{WK} (Cartan matrix,
restrictedness, Dynkin diagram) and superized them, as well as the
key notion --- that of Lie superalgebra --- for the case where
$p=2$. Following ideas of Weisfeiler and Kac \cite{WK}, and with the
help of \texttt{SuperLie} package \cite{Gr}, we classified {\it
finite dimensional modular Lie {\bf super}algebras with
indecomposable Cartan matrix}, see \cite{BGL}.

If a given indecomposable Cartan matrix $A$ is invertible, the Lie
(super)algebra $\fg(A)$ is simple, otherwise $\fg^{(i)}(A)/\fc$
--- the quotient of its first (for $i=1$) or second (for $i=2$)
derived algebra modulo the center $\fc$ --- is simple if $\size(A)
>1$ (we say that the {\it size} of an $n\times n$ matrix is equal to $n$).

The simple Lie algebra $\fg^{(i)}(A)/\fc$ --- in what follows in
such situation $i$ is equal to 1 or 2  (meaning that the derived
series of algebras stabilizes) --- does not possess any Cartan
matrix although the conventional sloppy practice is to refer to the
simple Lie (super)algebra $\fg^{(i)}(A)/\fc$ as \lq\lq possessing a
Cartan matrix''.

Elduque interpreted about a dozen of exceptional (when their
exceptional nature was only conjectured; now it is proved) simple
Lie superalgebras in characteristic 3 \cite{CE2} in terms of super
analogs of division algebras and collected them into a Supermagic
Square (an analog of Freudenthal's Magic Square); the rest of the
exceptional examples for $p=3$ and $p=5$, not entering the
Elduque\footnote{Although the first, as far as we know, superization
of Freudenthal's Square was performed by Martinez \cite{Mz} (for
$p=0$), Elduque went much further. It is instructive, however, to
compare the two squares.} Supermagic Square (the ones described for
the first time in \texttt{arXiv:math/0611391, math/0611392} and
\cite{BGL}) are, nevertheless, somehow affiliated to the Elduque
Supermagic Square \cite{El3}.

Very interesting, we think, is the situation in characteristic 2.
{\it A posteriori} we see that the list of Lie {\bf super}algebras
in characteristic 2 of the form $\fg(A)$ with an indecomposable
matrix $A$ is as follows:

{\sl In characteristic 2, take any finite dimensional simple Lie
algebra of the form $\fg(A)$  with indecomposable Cartan matrix $A$
(\cite{WK}) and declare some of the Chevalley generators of $\fg(A)$
odd (the corresponding diagonal elements of $A$ should be changed
accordingly $\ev$ to $0$ and $\od$ to $1$, see subsect.
\ref{normA}). Do this for each of the inequivalent Cartan matrices
of $\fg(A)$ and for any distribution of parities $I$ of the
Chevalley generators. Construct the Lie {\bf super}algebra
$s(\fg)(A, I)$ from these generators by the rules (\ref{myst1})
explicitly described in this paper. For the Lie superalgebra
$s(\fg)(A, I)$, list all its inequivalent Cartan matrices.}

Such superization may turn a given orthogonal Lie algebra into
ortho-orthogonal or periplectic Lie superalgebra; the three
exceptional Lie algebras of $\fe$ type turn into seven
non-isomorphic Lie superalgebras of $\fe$ type, whereas the $\fwk$
type Lie algebras turn into $\fbgl$ type Lie superalgebras.

The Lie superalgebra $s(\fg)(A, I)$ is simple if $A$ is invertible,
otherwise pass to $s(\fg)(A, I)^{(i)}/\fc$, where $i$ can be equal
to $1$ or $2$. We normalize the Cartan matrix so as to make the
parameter $I$ redundant and do not mention it in what follows.

In \cite{BGL}, we also listed {\it all} inequivalent Cartan matrices
$A$ for each given Lie (super)algebra $\fg(A)$. Although the number
of inequivalent Cartan matrices grows with the size of $A$, it is
easy to describe all possibilities for serial Lie (super)algebras.
Certain exceptional Lie superalgebras have dozens of inequivalent
Cartan matrices; nevertheless, there are several reasons to list all
of them: To classify all $\Zee$-gradings of a given $\fg(A)$ (in
particular, inequivalent Cartan matrices) is a very natural problem.
Besides, sometimes the knowledge of the best, for the occasion,
$\Zee$-grading is important, cf. \cite{RU} (all simple roots
non-isotropic), \cite{LSS} (all simple roots odd); for computations
\lq\lq by hand'' the cases where only one simple root is odd are
useful. In particular, the defining relations between the natural
(Chevalley) generators of $\fg(A)$ are of completely different form
for inequivalent $\Zee$-gradings and this is used in \cite{RU}.

\ssec{More on motivations} Being on the war path in the quest for
simple finite dimensional Lie superalgebras, we are interested in
complete description of their deformations. If $p=2$, these
deformations can only be described by means of divided power
homology introduced here for any $p>0$, together with divided power
cohomology. These notions are, clearly, of independent interest.

Since our results on deformation are not final, we decided to
illustrate our definitions with a final result of independent
interest --- presentations of simple modular Lie superalgebras with
Cartan matrix as well as presentations of simple modular Lie
algebras with Cartan matrix in the cases neglected so far: for $p=3$
and 2.

Recently we observe a rise of interest in presentations (by means of
generators and defining relations) of simple (and close to simple)
Lie (super)algebras occasioned by various applications of this
technical result, see \cite{GL1, Sa, Di, iPR} and references
therein, where presentations in terms of various other  types of
generators ({\it Jacobson, Silvester-t'Hooft, extremal}, etc.) are
given. Sometimes these other types of generators can be used as an
alternative to Chevalley generators; it is desirable, however, to
know the situations in which some of them are better (use less time
to construct the basis of the algebra they generate) than the others
or unavoidable as seems to be the case for Lie algebras of \lq\lq
matrices of complex size" (\cite{GL1}). Kornyak compared time needed
to present a given simple finite dimensional Lie algebra (over
$\Cee$) in terms of Chevalley generators and Serre relations  with
same in terms of Jacobson generators and Grozman-Leites relations,
see \cite{GL1}; the usefulness (in the above sense) of extremal
generators \cite{Di, iPR} is not yet compared with other
presentations, which is a pity: presentation in terms of them is
rather cumbersome.

For $p=2$, non-Serre relations appear even between the Chevalley
generators of simple Lie algebras. This is a new result.

Representations of quantum groups --- the deforms $U_q(\fg)$ of the
enveloping algebras --- at $q$ equal to a root of unity resemble,
even over $\Cee$, representations of Lie algebras in positive
characteristic and this is one more application that brings the
modular Lie (super)algebras and an {\bf explicit} form of their
presentations to the limelight.

\footnotesize

\ssec{Disclaimer}\label{discl}   Although presentation
--- description in terms of generators and relations --- is one of
the accepted ways to represent a given algebra, it seems that an
{\it explicit} form of the presentation is worth the trouble to
obtain only if this presentation is often in need, or (which is
usually the same) is sufficiently neat. The Chevalley generators of
simple finite dimensional Lie algebras over $\Cee$ satisfy simple
and neat relations (\lq\lq Serre relations'') and are often needed
for various calculations and theoretical discussions. Relations
between their analogs in the super case, although not so neat
(certain \lq\lq non-Serre relations'' appear), are still tolerable,
at least, for most Cartan matrices.

The defining relations expressed in terms of other generators,
different from Chevalley ones, are a bit too complicated to be used
by humans and were of academic interest until recently Grozman's
package \texttt{SuperLie} (\cite{Gr}) made the task of finding the
explicit expression of the defining relations for many types of Lie
algebras and superalgebras a routine exercise for anybody capable to
use \textit{Mathematica}.

What we usually need to know about defining relations is that there
are finitely many of them; hence the fact that some simple loop
superalgebras with Cartan matrix are {\bf not} finitely presentable
in terms of Chevalley generators was unexpected (although obvious as
an afterthought). The explicit form of defining relations for the
dozens or hundreds of systems of simple roots for the Lie
superalgebras of $e$ type (for $p=2$) can be easily obtained using
\texttt{SuperLie}, whereas for the exceptional simple Lie
superalgebras for $p>2$, it seems natural to list the relations
explicitly.

\normalsize \ssec{Main results} The definitions of new and
clarification of classical\footnote{The reader might be interested
in related problems, especially those posed by Deligne, see
\cite{LL}.} notions, especially, the definition of divided power
(co)homology.

We also define {\it Chevalley generators} and describe presentations
of finite dimensional modular Lie (super)al\-gebras of the form
$\fg(A)$ and $\fg(A)^{(i)}/\fc$ with indecomposable Cartan matrix
$A$ in terms of these generators.

If $p=2$, the non-Serre defining relations  for each Lie
superalgebra with indecomposable Cartan matrix are the same as for
the Lie algebra with the same (assuming $0=\ev$ and $1=\od$ on the
main diagonal) Cartan matrix. (This is proved for the exceptional
cases and $\fsl$ series; for the other series this is a conjecture
backed up by numerous examples.)\\

{\bf Acknowledgement.} We thank I.~Shchepochkina for help,
A.~Elduque for comments, O.~Shirokova for a \TeX pert help, and
A.~Protopopov for his help with the graphics, see \cite{Pro}. AL and
DL are thankful to MPIMiS, Leipzig, where part of the text was
written, for financial support and most creative environment. We are
very thankful to the referee for clarifying comments and help.

\section{What Lie superalgebra in characteristic
$2$ is}\label{Ssalgin2}

Let us give a naive definition of a Lie superalgebra for $p=2$. (For
a scientific one, as a Lie algebra in the category of
supervarieties, needed, for example, for a rigorous study and
interpretation of odd parameters of deformations, see \cite{LSh}.)
We define a Lie superalgebra as a superspace
$\fg=\fg_\ev\oplus\fg_\od$ such that the even part $\fg_\ev$ is a
Lie algebra, the odd part $\fg_\od$ is a $\fg_\ev$-module (made into
the two-sided one by symmetry; more exactly, by {\it anti}-symmetry,
but if $p=2$, it is the same) and on $\fg_\od$ a {\it squaring}
(roughly speaking, the halved bracket) is defined as a map
\begin{equation}\label{squaring}
\begin{array}{c}
x\mapsto x^2\quad \text{such that $(ax)^2=a^2x^2$ for any $x\in
\fg_\od$ and $a\in \Kee$, and}\\
{}(x+y)^2-x^2-y^2\text{~is a bilinear form on $\fg_\od$ with values
in $\fg_\ev$.}
\end{array}
\end{equation}
(We use a minus sign, so the definition also works for $p\neq 2$.)
The origin of this operation is as follows: If $\Char \Kee\neq 2$,
then for any Lie superalgebra $\fg$ and any odd element
$x\in\fg_\od$, the Lie superalgebra $\fg$ contains the element $x^2$
which is equal to the even element $\frac12 [x,x]\in\fg_\ev$. It is
desirable to keep this operation for the case of $p=2$, but, since
it can not be defined in the same way, we define it separately, and
then define the bracket of odd elements to be (this equation is
valid for $p\neq 2$ as well):
\begin{equation}\label{bracket}
{}[x, y]:=(x+y)^2-x^2-y^2.
\end{equation}
We also assume, as usual, that
\begin{itemize}
  \item if $x,y\in\fg_\ev$, then $[x,y]$ is the bracket on the Lie algebra;
  \item if $x\in\fg_\ev$ and $y\in\fg_\od$, then
$[x,y]:=l_x(y)=-[y,x]=-r_x(y)$, where $l$ and $r$ are the left and
right $\fg_\ev$-actions on $\fg_\od$, respectively.
\end{itemize}

The Jacobi identity involving odd elements now takes the following
form:
\begin{equation}\label{JI}
~[x^2,y]=[x,[x,y]]\text{~for any~} x\in\fg_\od, y\in\fg.
\end{equation}
If $\Kee\neq \Zee/2\Zee$, we can replace the condition (\ref{JI}) on
two odd elements by a simpler one:
\begin{equation}\label{JIbest}
[x,x^2]=0\;\text{ for any $x\in\fg_\od$.}
\end{equation}

Because of the squaring, the definition of derived algebras should
be modified. For any Lie superalgebra $\fg$, set $\fg^{(0)}:=\fg$
and
\begin{equation}\label{deralg}
\fg^{(1)}:=[\fg,\fg]+\Span\{g^2\mid g\in\fg_\od\} ,\quad
\fg^{(i+1)}:=[\fg^{(i)},\fg^{(i)}]+\Span\{g^2\mid
g\in\fg^{(i)}_\od\}.
\end{equation}

An even linear map $r\colon  \fg\tto\fgl(V)$ is said to be a {\it
representation of the Lie superalgebra}\index{representation of the
Lie superalgebra} $\fg$ (and the superspace $V$ is said to be a {\it
$\fg$-module})\index{$\fg$-module} if
\begin{equation}\label{repres}
\begin{array}{l}
r([x, y])=[r(x), r(y)]\quad \text{ for any $x, y\in
\fg$;}\\
r(x^2)=(r(x))^2\text{~for any $x\in\fg_\od$.}
\end{array}
\end{equation}

\ssec{Examples: Lie superalgebras preserving non-degenerate
(anti-)sym\-met\-ric forms} We say that two bilinear forms $B$ and
$B'$ on a superspace $V$ are {\it equivalent} if there is an even
invertible linear map $M\colon V\tto V$ such that
\begin{equation}\label{eqform}
B'(x,y)=B(Mx,My) \text{~~for any~~}x,y\in V.
\end{equation}
{\bf We fix some basis in $V$ and identify a given bilinear form
with its Gram matrix in this basis; we also identify any linear
operator on $V$ with its supermatrix in a fixed basis}.

Then two bilinear forms (rather supermatrices) are equivalent if and
only if there is an even invertible matrix $M$ such that
\begin{equation}\label{eqformM}
B'=MBM^T,\text{~~where $T$ is for transposition}.
\end{equation}

A bilinear form $B$ on $V$ is said to be {\it
symmetric}\index{bilinear form! symmetric} if $B(v, w)=B(w, v)$ for
any $v, w\in V$; a bilinear form is said to be {\it anti-symmetric}
if $B(v, v)=0$ for any $v\in V$.\index{bilinear form!
anti-symmetric}

A homogeneous\footnote{Hereafter, as always in Linear Algebra in
superspaces, all formulas of linear algebra defined on homogeneous
elements only are supposed to be extended to arbitrary ones by
linearity.} linear map $F$ is said to preserve a bilinear form $B$,
if\footnote{Hereafter, $p$ denotes both parity defining a
superstructure and the characteristic of the ground field; the
context is, however, always clear.}
\[
 B(F x, y)+(-1)^{p(x)p(F)}B(x, Fy)=0\quad\text{for any~}x,y\in V.
\]
All linear maps preserving a given bilinear form constitute a Lie
sub(super)algebra $\faut_B(V)$ of $\fgl(V)$ denoted
$\faut_B(n)\subset\fgl(n)$ in matrix realization and consisting of
the supermatrices $X$ such that
\[
BX+(-1)^{p(X)}X^{st}B=0, \]
where the {\it supertransposition} $st$
acts as follows (in the standard format):
\[
st\colon  \begin{pmatrix}A&B\\C&D\end{pmatrix}\tto
\begin{pmatrix}A^t&-C^t\\B^t&D^t\end{pmatrix}.
\]

Consider the case of purely even space $V$ of dimension $n$ over a
field of characteristic $p\neq 2$. Every non-zero form $B$ can be
uniquely represented as the sum of a symmetric and an anti-symmetric
form and it is possible to consider automorphisms and equivalence
classes of each summand separately.

If the ground field $\Kee$ of characteristic $p>2$
satisfies\footnote{Aside: We thought that one should require
perfectness of $\Kee$, i.e., $\Kee^p=\Kee$ but the referee suggested
a simple counterexample for $\Kee=\Zee/3$ with 2 non-equivalent
types of non-degenerate symmetric forms. In this paper $\Kee$ is
algebraically closed; over fields  algebraically non-closed, there
are more types of symmetric forms.} $\Kee^2=\Kee$, then there is
just one equivalence class of non-degenerate symmetric even forms,
and the corresponding Lie algebra $\faut_B(V)$ is called {\it
orthogonal} and denoted $\fo_B(n)$ (or just $\fo(n)$).
Non-degenerate anti-symmetric forms over $V$ exist only if $n$ is
even; in this case, there is also just one equivalence class of
non-degenerate antisymmetric even forms; the corresponding Lie
algebra $\faut_B(n)$ is called  {\it symplectic} and denoted
$\fsp_B(2k)$ (or just $\fsp(2k)$). Both algebras $\fo(n)$ and
$\fsp(2k)$ are simple.

{\bf If $p=2$, the space of anti-symmetric bilinear forms is a
subspace of symmetric bilinear forms}. Also, instead of a unique
representation of a given form as a sum of an anti-symmetric and
symmetric form, we have a subspace of symmetric forms and the
quotient space of non-symmetric forms; it is not immediately clear
what to take for a representative of a given non-symmetric form. For
an answer and classification, see Lebedev's thesis \cite{LeD} and
\cite{Le1}. There are no new simple Lie superalgebras associated
with non-symmetric forms, so we confine ourselves to symmetric ones.

Instead of orthogonal and symplectic Lie algebras we have two
different types of orthogonal Lie algebras (see Theorem
\ref{SymForm}). Either the derived algebras of these algebras or
their quotient modulo center are simple if $n$ is large enough, so
the canonical expressions of the forms $B$ are needed as a step
towards classification of simple Lie algebras in characteristic $2$
which is an open problem, and as a step towards a version of this
problem for Lie superalgebras, even less investigated.

In \cite{Le1}, Lebedev showed that, with respect to the above
natural equivalence of forms (\ref{eqformM}), every {\bf even}
symmetric non-degenerate form on a superspace of dimension
$n_{\ev}|n_{\od}$ over a perfect  (i.e., such that every element of
$\Kee$ has a square root\footnote{Since $a^2-b^2=(a-b)^2$ if $p=2$,
it follows that no element can have two distinct square roots.})
field of characteristic $2$ is equivalent to a form of the shape
(here: $i=\bar 0$ or $\bar 1$ and each $n_i$ may equal to 0),
\[
B=\mat{ B_{\ev}&0\\0&B_{\od}}, \quad \text{where
$B_i=\begin{cases}1_{n_i}&\text{if $n_i$ is odd,}\\
\text{either $1_{n_i}$ or $\Pi_{n_i}$}&\text{if $n_i$ is
even,}\end{cases}$}
\] and where
\[\Pi_n=\begin{cases} \mat{0&1_k\\
1_k&0}&\text{if $n=2k$},\\[3mm]
\mat{0&0&1_k\\
0&1&0\\
1_k&0&0}&\text{if $n=2k+1$}.\end{cases}
\]
(In other words, the bilinear forms with matrices $1_{n}$ and
$\Pi_{n}$ are equivalent if $n$ is odd and non-equivalent if $n$ is
even.) The Lie superalgebra preserving $B$
--- by analogy with the orthosymplectic Lie superalgebras $\fosp$ in
characteristic $0$ we call it {\it ortho-orthogonal} and denote
$\foo_B(n_\ev|n_\od)$
--- is spanned by the supermatrices which in the standard format are
of the form
\[
\mat{ A_{\ev}&B_{\ev}C^TB_{\od}^{-1}\\C&A_{\od}}, \;
\begin{matrix}\text{where
$A_{\ev}\in\fo_{B_{\ev}}(n_\ev)$, $A_{\od}\in\fo_{B_{\od}}(n_\od)$, and}\\
\text{$C$ is arbitrary $n_{\od}\times n_{\ev}$ matrix.}\end{matrix}
\]
Since, as is easy to see, \[\foo_{\Pi I}(n_\ev|n_\od)\simeq
\foo_{I\Pi}(n_\od|n_\ev),\] we do not have to consider the Lie
superalgebra $\foo_{\Pi I}(n_\ev|n_\od)$ separately unless we study
Cartan prolongations where the difference between these two
incarnations of one algebra is vital: For the one, the prolong is
finite dimensional (the automorphism algebra of the $p=2$ analog of
the Riemann geometry), for the other one it is infinite dimensional
(an analog of the Lie superalgebra of Hamiltonian vector fields).

For an {\bf odd} symmetric form $B$ on a superspace of dimension
$(n_{\ev}|n_{\od})$ over a field of characteristic $2$ to be
non-degenerate, we need $n_{\ev}=n_{\od}$, and every such form $B$
is equivalent to $\Pi_{k|k}$, where $k=n_{\ev}=n_{\od}$, and which
is same as $\Pi_{2k}$ if the superstructure is forgotten. This form
is preserved by linear transformations with supermatrices in the
standard format of the shape
\begin{equation}
\label{pe} \mat{ A&C\\D&A^T}, \quad \text{where $A\in\fgl(k)$, $C$
and $D$ are symmetric
$k\times k$ matrices}. 
\end{equation}
As over $\Cee$ or $\Ree$, the Lie superalgebra of linear maps
preserving $B$ will be referred to as {\it periplectic}, as A.~Weil
suggested, and denoted $\fpe_B(k)$ or just $\fpe(k)$. Note that even
the superdimensions of the characteristic $2$ versions of the Lie
(super)algebras $\faut_B(k)$ differ from their analogs in other
characteristics for both even and odd forms $B$.

Now observe that
\begin{equation}\label{nt}
\begin{tabular}{c}
{\bf The fact that two bilinear forms are inequivalent does}\\
{\bf not, generally, imply that the Lie (super)algebras that}\\
{\bf preserve them are not isomorphic}.
\end{tabular}
\end{equation} In \cite{Le1}, Lebedev proved that for the {\it non-degenerate
symmetric} forms, the implication spoken about in (\ref{nt}) is,
however, true (bar a few exceptions), and therefore we have several
types of non-isomorphic Lie (super) algebras (except for occasional
isomorphisms intermixing the types, e.g.,
$\foo_{I\Pi}\simeq\foo_{\Pi I}$ and
$\foo_{\Pi\Pi}^{(1)}(6|2)\simeq\fpe^{(1)}(4)$).

The problem of describing preserved bilinear forms has two levels:
we can consider {\it linear transformations} (Linear Algebra) and
{\it arbitrary coordinate changes} (Differential Geometry). 
In the literature, both levels are completely investigated, except
for the case where $p=2$. More precisely, the fact that the
non-split and split forms of the Lie algebras that preserve the
symmetric bilinear forms are not always isomorphic was never
mentioned. (Although known for the Chevalley groups preserving these
forms, cf. \cite{St}, these facts do not follow from each other
since there is no analog of Lie theorem on the correspondence
between Lie groups and Lie algebras.) Here we consider the Linear
Algebra aspect, for the Differential Geometry related to the objects
considered here, see \cite{Le2}.


\sssec{Known facts: The case $p= 2$}\label{s2.1.5} The following
facts are given for clarity: lecturing on these results during the
past several years we have encountered incredulity of the listeners
based on several false premises intermixed with correct statements:
\lq\lq the question sounds classical and so had been solved by
classics without doubt (the solution just has to be dug out from
paper diluvium)", \lq\lq this is known for quadratic forms (don't
you know about Arf invariant?!)", \lq\lq there are two
non-isomorphic types of simple finite orthogonal
groups\footnote{Observe that these groups are defined as preserving
quadratic forms, not bilinear ones.} acting on $2n$-dimensional
space, so what's new?", and so on.

1) With any symmetric bilinear form $B$ the quadratic form
$Q(x):=B(x, x)$ is associated. Arf has discovered {\it the Arf
invariant} --- an important invariant of non-degenerate quadratic
forms in characteristic $2$; for an exposition, see \cite{Dye}. Two
such forms are equivalent if and only if their Arf invariants are
equal.

The other way round, given a quadratic form $Q$, we define a
symmetric bilinear form, called {\it the polar form}\index{bilinear
form! polar} of $Q$, by setting
\[
B_Q(x, y)=Q(x+y)-Q(x)-Q(y).
\]
The Arf invariant can not, however, be used for classification of
symmetric bilinear forms because one symmetric bilinear form can
serve as the polar form for two non-equivalent (and having different
Arf invariants) quadratic forms. Moreover, {\bf not every symmetric
bilinear form can be represented as a polar form. If $p= 2$, the
correspondence $Q\longleftrightarrow B_Q$ is not one-to-one}.

2) Recall that the space of anti-symmetric forms (their matrices are
zero-diagonal ones) is a subspace in the space of symmetric forms.
Albert \cite{A} classified symmetric bilinear forms over a field of
characteristic $2$ and proved that (we have in mind symmetric forms
only)
\begin{equation*}\label{123}
\begin{minipage}[l]{12cm}
\begin{description}
                 \item[(1)] two anti-symmetric forms of equal ranks are equivalent;
                 \item[(2)] every non-anti-symmetric form has a matrix which is equivalent
to a diagonal matrix;
                 \item[(3)] if $\Kee$ is perfect, then every two non-anti-symmetric forms of
equal ranks are equivalent.
               \end{description}

\end{minipage}
\end{equation*}

\parbegin{Remarks} 1) Over a field $\Kee$ of
characteristic $2$, Albert also obtained certain results on the
classification of quadratic forms  (considered as elements of the
quotient space of all bilinear forms modulo the space of
anti-symmetric forms). In particular, he showed that if $\Kee$ is
algebraically closed, then every quadratic form is equivalent to
exactly one of the forms
\begin{equation}
\label{aforms} x_1x_{r+1}+\dots+x_rx_{2r}\;\text{ or }\;
x_1x_{r+1}+\dots+x_rx_{2r}+x_{2r+1}^2,
\end{equation}
where $2r$ is the rank of the form. Lebedev \cite{Le1} used this
result in the study of Lie algebras preserving the contact
structure.

2) Lebedev \cite{Le1} also suggested canonical forms (or rather of
their classes modulo the subspace of symmetric forms) of
non-symmetric bilinear forms and classified them. This result is
also related to a result of Albert and --- rather unexpectedly
--- with contact structures on superspaces.
\end{Remarks}

\parbegin{Theorem}[\cite{Le1}]\label{SymForm}\label{s2.2.1} Let $\Kee$ be a
perfect field of characteristic $2$. Let $V$ be an
$n$-di\-men\-sio\-nal space over $\Kee$.

\textup{1)} For $n$ odd, there is only one equivalence class of
non-degenerate symmetric bilinear forms on $V$.

\textup{2)} For $n$ even, there are two equivalence classes of
non-degenerate symmetric bilinear forms, one --- with at least one
non-zero element on the main diagonal --- contains $1_n$ and the
other one --- all its Gram matrices are {\it zero-diagonal}
--- contains $S_n:=\antidiag(1, \dots, 1)$ and $\Pi_n$.
\end{Theorem}

In view of (\ref{nt}) the statement of the next Lemma (proved in
\cite{Le1, BGL}) is non-trivial.

\parbegin{Lemma}\label{noniso} \textup{1)} The Lie algebras
$\fo_I(2k)$ and $\fo_\Pi(2k)$ are not isomorphic (though are of the
same dimension); the same applies to their derived algebras:

\textup{2)} $\fo_I^{(1)}(2k)\not\simeq \fo_\Pi^{(1)}(2k)$, though
$\dim\fo_I^{(1)}(2k)=\dim\fo_\Pi^{(1)}(2k)$;

\textup{3)} $\fo_I^{(2)}(2k)\not\simeq \fo_\Pi^{(2)}(2k)$ unless
$k=1$.
\end{Lemma}

Based on these results, Lebedev described all the (four) possible
analogs of the Poisson bracket, and (there exists just one) contact
bracket. Similar results for the odd bilinear form yield a
description of the anti-bracket (a.k.a. Buttin bracket), and the
(peri)contact bracket, compare \cite{Le2} with \cite{LSh}. The
quotients of the Poisson and Buttin Lie (super)algeb\-ras modulo
center --- analogs of Lie algebras of Hamiltonian vector fields, and
their divergence-free subalgebras
--- are also described in \cite{Le2}.

\section{Analogs of functions and vector fields for $p>0$}
\ssec{Divided powers} Let us consider the supercommutative
superalgebra $\Cee[x]$ of polynomials in $a$ indeterminates $x =
(x_1,...,x_a)$, for convenience ordered in a \lq\lq standard
format'', i.e., so that the first $m$ indeterminates are even and
the rest $n$ ones are odd ($m+n=a$). Among the integer bases of
$\Cee[x]$ (i.e., the bases, in which the structure constants are
integers), there are two canonical ones,
--- the usual, monomial, one and the basis of {\it divided powers},
\index{divided power} which is constructed in the following way.

For any multi-index $\underline{r}=(r_1, \ldots , r_a)$, where
$r_1,\dots,r_m$ are non-negative integers, and $r_{m+1},\dots,r_n$
are $0$ or $1$, we set
\[
u_i^{(r_{i})} :=
\frac{x_i^{r_{i}}}{r_i!}\quad \text{and}\quad u^{(\underline{r})} :=
\prod\limits_{i=1}^a u_i^{(r_{i})}.
\]
These $u^{(\underline{r})}$ form an integer basis of $\Cee[x]$.
Clearly, their multiplication relations are
\begin{equation}
\label{divp}
\renewcommand{\arraystretch}{1.4}
\begin{array}{l}  u^{(\underline{r})} \cdot u^{(\underline{s})} =
\prod\limits_{i=m+1}^n
\min(1,2-r_i-s_i)\cdot(-1)^{\sum\limits_{m<i<j\leq a} r_js_i}\cdot
\binom {\underline{r} + \underline{s}} {\underline{r}}
u^{(\underline{r} + \underline{s})},  \\
\text{where}\quad \binom {\underline{r} + \underline{s}}
{\underline{r}}:=\prod\limits_{i=1}^m\binom {r_{i} + s_{i}} {r_{i}}.
\end{array}
\end{equation}
In what follows, for clarity, we will write exponents of divided
powers in parentheses, as above, especially if the usual exponents
might be encountered as well.

Now, for an arbitrary field $\Kee$ of characteristic $p>0$, we may
consider the supercommutative superalgebra $\Kee[u]$ spanned by
elements $u^{(\underline{r})}$ with multiplication relations
(\ref{divp}). For any $m$-tuple $\underline{N} = (N_1,..., N_m)$,
where $N_i$ are either positive integers or infinity, denote (we
assume that $p^\infty=\infty$)
\begin{equation}
\label{u;N} \cO(m; \underline{N}):=\Kee[u;
\underline{N}]:=\Span_{\Kee}\left(u^{(\underline{r})}\mid r_i
\begin{cases}< p^{N_{i}}&\text{for $i\leq m$}\\
=0\text{ or 1}&\text{for $i>m$}\end{cases}\right).
\end{equation}
 As is clear from (\ref{divp}), $\Kee[u; \underline{N}]$ is a
subalgebra of $\Kee[u]$. The algebra $\Kee[u]$ and its subalgebras
$\Kee[u; \underline{N}]$ are called the {\it algebras of divided
powers;} they can be considered as analogs of the polynomial
algebra.

Only one of these numerous algebras of divided powers
$\cO(n;\underline{N})$ are indeed generated by the indeterminates
declared: If $N_i=1$ for all $i$. Otherwise, in addition to the
$u_i$, we have to add $u_i^{(p^{k_i})}$ for all $i\leq m$ and all
$k_i$ such that $1<k_i<N_i$ to the list of generators. Since any
derivation $D$ of a given algebra is determined by the values of $D$
on the generators, we see that $\fder(\cO[m; \underline{N}])$ has
more than $m$ functional parameters (coefficients of the analogs of
partial derivatives) if $N_i\neq 1$ for at least one $i$. Define
{\it distinguished partial derivatives}\index{derivative! partial,
distinguished}
 by setting
\[
\partial_i(u_j^{(k)})=\delta_{ij}u_j^{(k-1)}\;\text{ for any $k<p^{N_j}$}.
\]

The simple vectorial Lie algebras over $\Cee$ have only one
parameter: the number of indeterminates. If $\Charr~ \Kee =p>0$, the
vectorial Lie algebras acquire one more parameter: $\underline{N}$.
For Lie superalgebras, $\underline{N}$ only concerns the even
indeterminates.

The Lie (super)algebra of all derivations $\fder(\cO[m;
\underline{N}])$ turns out to be not so interesting as its {\it Lie
subsuperalgebra of distinguished derivations}: Let
\begin{equation}\label{vect-super}
\begin{array}{c} \fvect (m; \underline{N}|n) \;\text{ a.k.a }W(m;
\underline{N}|n)\; \text{ a.k.a }\\
\fder_{dist} \Kee[u;
\underline{N}]=\Span_{\Kee}\left(u^{(\underline{r})}\partial_k\mid
r_i\begin{cases}< p^{N_{i}}&\text{for $i\leq m$},\\
=0\text{ or 1}&\text{for $i>m$};\end{cases}
\quad  1\leq k\leq
n\right)\end{array}
\end{equation} be the {\it general vectorial Lie algebra of
distinguished derivations}. The next notions are analogs of the
polynomial algebra of the dual space.

\ssec{Symmetric differential forms and exterior differential forms}
In what follows, as is customary in modern geometry, we use
the\index{$\wedge$ product of differential forms} antisymmetric
$\wedge$ product for the analogs of the {\it exterior differential}
forms, and the symmetric $\circ$\index{$\circ$ product of
differential forms} product for the {\it symmetric differential}
forms, e.g., analogs of the metrics. We can also consider the
divided power versions of the exterior and symmetric forms because
both types of forms generate (in the divided sense) supercommutative
superalgebras depending not only on the $u_i$, as above, but also on
$du_i$, such that $p(du_i)=p(u_i)$ in the symmetric case and
$p(du_i)=p(u_i)+\od$ in the exterior case. Usually we suppress the
$\wedge$ or $\circ$ signs, since all is clear from the context,
unless both multiplications are needed simultaneously. We have,
however, to distinguish the non-divided $\wedge$ or $\circ$ from
their divided counterparts! This is important since both non-divided
and divided products are often needed simultaneously. Fortunately,
in this paper, we only need divided products, so we do not use more
appropriate notation $\stackrel{d}{\wedge}$, although use just it,
not $\wedge$.

Considering {\it exterior} differential forms, we use divided powers
$dx_i^{(\wedge k)}$\index{divided powers of differential forms} with
multiplication relations (\ref{divp}), where the indeterminates are
now the $dx_i$ of parity $p(x_i)+\od$, and the Lie derivative along
the vector field $X$ is given by the formula
\[
L_X(dx_i^{(\wedge k)})=(L_Xdx_i)\wedge dx_i^{(\wedge k-1)}.
\]
Note that if we consider divided power differential forms in
characteristic $2$, then, for $x_i$ odd, we have $dx_i\wedge
dx_i=2(dx_i^{(\wedge 2)})=0$. (If $x_i$ is even, then $dx_i\wedge
dx_i=0$, anyway.)

Considering  divided powers of chains and cochains of Lie
superalgebras affects the formula for the (co)chain differentials.
For cochains of a given Lie superalgebra $\fg$, this only means that
a divided power of an odd element must be differentiated as a whole:
\begin{equation}\label{dcoch}
d(\varphi^{(\wedge k)})=d\varphi\wedge \varphi^{(\wedge (k-1))}
\text{~for any~}\varphi\in(\fg^*)_\od.
\end{equation}
For chains, the modification is a little more involved: Let $g_1,
\dots,  g_n$ be a basis of $\fg$. Then for chains of $\fg$ with
coefficients in a right module $A$, and $a\in A$, we have
\begin{equation}\label{dchane}
\footnotesize
\renewcommand{\arraystretch}{1.10}
\begin{array}{ll}
d\left(a\otimes\bigwedge\limits_{i=1}^n g_i^{(\wedge
r_i)}\right)=&\sum\limits_{p(g_k)=\od,~r_k\geq
2}a\otimes\bigwedge\limits_{i<k} g_i^{(\wedge r_i)}\wedge
g_k^2\wedge g_k^{(\wedge(r_k-2))}\wedge \bigwedge\limits_{i>k}
g_i^{(\wedge r_i)}+\\
&\sum\limits_{1\leq k<l\leq n,~~r_k,r_l\geq 1}
(-1)^{\sum\limits_{k<m<l} r_mp(g_m)}a\otimes\bigwedge\limits_{i<k}
g_i^{(\wedge r_i)}\wedge \\
& [g_k,g_l]\wedge
g_k^{(\wedge(r_k-1))}\wedge\bigwedge\limits_{k<i<l} g_i^{(\wedge
r_i)}\wedge g_l^{(\wedge(r_l-1))}\wedge \bigwedge\limits_{i>l}
g_i^{(\wedge
r_i)}+\\
&\sum\limits_{r_k\geq 1} (-1)^{p(g_k)\sum\limits_{m<k}r_mp(g_m)}
(ag_k)\otimes \bigwedge\limits_{i<k} g_i^{(\wedge r_i)}\wedge
g_k^{(\wedge(r_k-1))}\wedge \bigwedge\limits_{i>k} g_i^{(\wedge
r_i)}.\end{array}
\end{equation}
Denote the divided power cohomology by $DPH^{i, \uN}(\fg; M)$ and
divided power homology by $DPH_{i, \uN}(\fg; M)$. Note that if $\fg$
is a Lie {\bf super}algebra and $p=2$, we can not interpret its
generating relations in terms of the 2nd homology $H_{2}(\fg)$, as
we do for $p\neq 2$: Instead, we {\bf must} use divided powers
homology $DPH_{2, \uN}(\fg):=DPH_{2, \uN}(\fg;\Kee)$ (with $\uN$
such that $N_i\geq 2$ for all $i$) since otherwise we won't be able
to take into account the relations of the form $x^2=0$ for $x$ odd.
For the same reason, to interpret deformations in the same
situation, we need $DPH^{2, \uN}(\fg;\fg)$, not $H^{2}(\fg;\fg)$.

\sssbegin{Problem} To define the divided power (co)homology as the
derived functor, we have to completely modify the representation
theory and, in particular, the notion of the universal enveloping
algebra. {\em We do not know a precise definition of the \lq \lq
divided power universal enveloping algebra" but conjecture that is
can be found along the way hinted at in \cite{LL}.}
\end{Problem}


\sssec{A useful Lemma} We computed cohomology using Grozman's
\texttt{Mathemati\-ca}-based package \texttt{SuperLie}. The formula
of the following lemma was helpful in the computations. For any
finite dimensional Lie (super)algebra $\mathfrak{g}$, all cochains
with non-trivial coefficients in a $\fg$-module $M$ can be expressed
as sums of tensor products of the form $m\otimes\omega$, where $m\in
M$ and $\omega\in\bigwedge^{\bcdot}(\mathfrak{g}^*)$. We are working
with a fixed basis of $M$ and the dual basis of $\mathfrak{g}^*$.
For simplicity, the following Lemma is formulated for Lie algebras,
its superization is routine, by means of the Sign Rule.

\begin{Lemma} \label{lem1}
For any $c=m\otimes\omega$, where $m\in M$ and
$\omega\in\bigwedge^{r}(\mathfrak{g}^*)$, let $dc$ denote the
coboundary of $c$ in the complex with coefficients in $M$, while
$d\omega$ denotes the coboundary in the complex with trivial
coefficients and $dm$ denotes the coboundary of $m\in M$ considered
as a $0$-cochain in the complex with coefficients in $M$. If
$c=m\otimes\omega$, then $dc=m\otimes d\omega + dm\wedge\omega$.
\end{Lemma}

\begin{proof}
For any $x_1,\ldots,x_{r+1}\in \fg$, we have:
\[
\begin{split}
 dc(x_1,\ldots,x_{r+1})&=\sum_{1\leq i<j\leq r+1}(-1)^{i+j-1}a\otimes
 \omega([x_i,x_j],x_1,\ldots,\hat{x_i},\ldots,\hat{x_j},\ldots,x_{r+1})+{}\\
 &+\sum_{1\leq i\leq r+1}(-1)^{i}x_i(m)\otimes
 \omega(x_1,\ldots,\hat{x_i},\ldots,x_{r+1})=\\
 &=(a\otimes d\omega)(x_1,\ldots,x_{r+1})+(dm\wedge\omega)
 (x_1,\ldots,x_{r+1}).\qed
\end{split}
\]
\noqed\end{proof}

\paragraph{In characteristic $2$}. Let now $p=2$. In this subsection $-1=1$, of course;
the signs are kept to make expressions look like in characteristics
0. The following definition of Lie algebra cohomology in $\Char=2$
is implemented in \texttt{SuperLie}:

For $1$-cochains with trivial coefficients, the codifferential  is
defined as an operation dual to the Lie bracket:
\[
 d\colon \mathfrak{g}^*\rightarrow \mathfrak{g}^* \wedge \mathfrak{g}^*.
\]
For $q$-cochains with trivial coefficients, $d$ is defined via the
Leibniz rule. For cochains with coefficients in a module $M$, we set
\[
\begin{split}
 &d(m):=-\mathop{\sum}\limits_{1\leq i\leq \dim M} g_i(m) \otimes g_i^*,\\
 &d(m \otimes \omega):=d(m)\wedge\omega +  m \otimes d(\omega)
\end{split}
\]
for any $m\in M$, any $r$-cochain $\omega$, where $r>0$, and any
basis $g_i$ of $M$, cf. Lemma \ref{lem1}.

\section{What $\fg(A)$ is}\label{Sg(A)}

\ssec{Warning: $\fpsl$ has no Cartan matrix. The relatives of $\fsl$
and $\fpsl$ that have Cartan matrices}\label{warn} For the most
reasonable definition of Lie algebra with Cartan matrix over $\Cee$,
see \cite{K}. The same definition applies, practically literally, to
Lie superalgebras and to modular Lie algebras and to modular Lie
superalgebras. However, the usual sloppy practice is to attribute
Cartan matrices to (usually simple) Lie (super)algebras none of
which, strictly speaking, has a Cartan matrix!

Although it may look strange for those with non-super experience
over $\Cee$, neither the simple modular Lie algebra $\fpsl(pk)$, nor
the simple modular Lie superalgebra $\fpsl(a|pk+a)$, nor --- in
characteristic $0$ --- the simple Lie superalgebra $\fpsl(a|a)$
possesses a Cartan matrix. Their central extensions ($\fsl(pk)$, the
modular Lie superalgebra $\fsl(a|pk+a)$, and --- in characteristic
$0$ --- the Lie superalgebra $\fsl(a|a)$) do not have Cartan matrix,
either.

Their relatives possessing a Cartan matrix are, respectively,
$\fgl(pk)$, $\fgl(a|pk+a)$, and $\fgl(a|a)$, and for the grading
operator we take the matrix unit $E_{1,1}$.

Since all the Lie (super)algebras involved (the simple one, its
central extension, the derivation algebras thereof) are often needed
simultaneously (and only representatives of one of these types of
Lie (super)algebras are of the form $\fg(A)$), it is important to
have (preferably short and easy to remember) notation for each of
them. For example, in addition to $\fpsl$, $\fsl$, $\fp\fgl$ and
$\fgl$, we have:

{\bf for $p=3$}: $\fe(6)$ is of dimension 79, then
$\dim\fe(6)^{(1)}=78$, whereas the \lq\lq simple core'' is
$\fe(6)^{(1)}/\fc$ of dimension 77;

$\fg(2)$ is not simple, its \lq\lq simple core'' is isomorphic to
$\fpsl(3)$;

{\bf for $p=2$}: $\fe(7)$ is of dimension 134, then
$\dim\fe(7)^{(1)}=133$, whereas the \lq\lq simple core'' is
$\fe(7)^{(1)}/\fc$ of dimension 132;

$\fg(2)$ is not simple, its \lq\lq simple core'' is isomorphic to
$\fpsl(4)$;

{\bf the orthogonal Lie algebras and their super analogs} are
considered in detail later.

In our main examples, $\sdim \fg(A)^{(i)}/\fc=d|B$  whereas the
notation $D/d|B$ means that $\sdim \fg(A)=D|B$. The general formula
is
\begin{equation}\label{dims}
d=D-2(\size(A)-\rk(A))\text{~~and~~}i=\size(A)-\rk(A).
\end{equation}

\ssec{Generalities} Let $A=(A_{ij})$ be an $n\times n$-matrix with
elements in $\Kee$ with $\rk A=n-l$. Complete $A$ to an $(n+l)\times
n$-matrix $\begin{pmatrix} A\\ B \end{pmatrix}$ of rank $n$. (Thus,
$B$ is an $l\times n$-matrix.)

Let the elements $e_i^\pm, h_i$, where $i=1,\dots,n$, and $d_k$,
where $k=1,\dots,l$, generate a Lie superalgebra denoted $\tilde
\fg(A, I)$,\index{$\tilde\fg(A, I)$} where $I=(p_1, \dots
p_n)\in(\Zee/2)^n$ is a collection of parities ($p(e_i^\pm)=p_i$,
the parities of the $d_k$'s being $\ev$), free except for the
relations
\begin{equation}\label{gArel_0}
\renewcommand{\arraystretch}{1.4}
\begin{array}{l}
{}[e_{i}^+, e_{j}^-] = \delta_{ij}h_i; \quad [h_i, e_{j}^\pm]=\pm
A_{ij}e_{j}^\pm;\quad [d_k, e_{j}^\pm]=\pm
B_{kj}e_{j}^\pm;\\
{}[h_i, h_j]=[h_i, d_k]=[d_k, d_m]=0\text{~~for any $i, j, k,
m$}.\end{array}
\end{equation}
The Lie superalgebra  $\tilde \fg(A, I)$ is $\Zee^n$-graded with
\begin{equation}\label{deg}
\renewcommand{\arraystretch}{1.4}
\begin{array}{l}
\deg e_{i}^\pm =(0, \ldots, 0, \pm 1, 0, \ldots, 0)\\
\deg h_i=\deg d_k=(0, \ldots, 0)\text{~~for any $i, k$}.\end{array}
\end{equation}

Let $\fh$ denote the linear span of the $h_i$'s and $d_k$'s. Let
$\tilde \fg(A, I)^\pm$ denote the Lie subsuperalgebras in $\tilde
\fg(A, I)$ generated by  $e_{1}^\pm, \ldots, e_{n}^\pm$. Then
\[\tilde \fg(A, I)=\tilde \fg(A, I)^-\oplus \fh\oplus\tilde \fg(A, I)^+,\]
where the homogeneous component of degree $(0, \ldots, 0)$ is just
$\fh$.

The Lie subsuperalgebras $\tilde \fg(A, I)^\pm$ are homogeneous in
this $\Zee^n$-grading, and there is a
\begin{equation}\label{r}
\text{maximal homogeneous (in this $\Zee^n$-grading) ideal $\fr$
such that $\fr\cap\fh=0$.}
\end{equation} The ideal $\fr$ is just the sum of homogeneous ideals
whose homogeneous components of degree $(0, \ldots, 0)$ is trivial.

As $\rk A =n-l$, there exists an $l\times n$-matrix $T=(T_{ij})$ of
rank $l$ such that
\begin{equation}\label{TA=0}
TA=0.
\end{equation}
Let
\begin{equation}\label{central}
c_i=\sum_{j=1}^n T_{ij}h_j, \text{~~where~~} i=1,\dots,l.
\end{equation}
Then, from the properties of the matrix $T$, we deduce that
\begin{equation}\label{central1}
\begin{tabular}{l}
a) the elements $c_i$ are linearly independent; let $\fc$ be the space they span;\\
b) the elements $c_i$ are central, because\\
$[c_i,e_j^\pm]=\pm\left(\sum\limits_{k=1}^n T_{ik}A_{kj}\right)
e_j^\pm=\pm (TA)_{ij} e_j^\pm \stackrel{\eqref{TA=0}}{=}0$.
\end{tabular}
\end{equation}

The Lie (super)algebra $\fg(A, I)$\index{$\fg(A, I)$} is defined as
the quotient $\tilde \fg(A, I)/\fr$ and is called the {\it Lie
(super)algebra with Cartan matrix $A$ (and parities $I$)}. Note that
this coincides with the definition in \cite{CE} of the {\it
contragredient} Lie superalgebras, although written in a slightly
different way. Condition \eqref{r} modified as
\begin{equation}\label{rc}
\text{maximal homogeneous (in this $\Zee^n$-grading) ideal $\fs$
such that $\fs\cap\fh=\fc$}
\end{equation}
leads to what in \cite{CE} is called the {\it centerless
contragredient} Lie superalgebra, cf. \cite{Bi}.

By abuse of notation we denote by $e_i^\pm, h_i, d_k$ and $\fc$
their images in $\fg(A,I)$ and $\fg(A,I)^{(1)}$.

The Lie superalgebra $\fg(A,I)$ inherits, clearly, the
$\Zee^n$-grading of $\tilde \fg(A, I)$. The non-zero elements
$\alpha\in\Zee^n\subset \Ree^n$ such that the homogeneous component
$\fg(A,I)_\alpha$ is non-zero are called {\it roots}.
\index{Root}The set $R$ of all roots is called {\it the root
system}\index{System!root} of $\fg$. Clearly, the subspaces
$\fg_\alpha$ are purely even or purely odd, and the corresponding
roots are said to be \textit{even} or \textit{odd}.

The additional to \eqref{gArel_0} relations that turn $\tilde \fg(A,
I)^\pm$ into $\fg(A, I)^\pm$ are of the form $R_i=0$ whose left
sides are implicitly described as follows:
\begin{equation}
\label{myst1}
\begin{split}
 &\text{the $R_i$ that generate the maximal ideal $\fr$.}
 \end{split}
\end{equation}
{\bf The explicit description of these additional relations} forms
the main bulk of this paper.

\ssec{Roots and weights}\label{roots} In this subsection, $\fg$
denotes one of the algebras $\fg(A,I)$ or $\tilde{\fg}(A,I)$.

The elements of $\fh^*$ are called {\it weights}.\index{weight} For
a given weight $\alpha$, the {\it weight subspace} of a given
$\fg$-module $V$ is defined as
\[
V_\alpha=\{x\in V\mid \text{an integer $N>0$ exists such that
$(\alpha(h)-\ad_h)^N x=0$ for any $h\in\fh$}\}.
\]

Any non-zero element $x\in V$ is said to be {\it of weight
$\alpha$}. For the roots, which are particular cases of weights if
$p=0$, the above definition is inconvenient: In the modular analog
of the following useful statement summation should be over roots
defined in the previous subsection.

\sssbegin{Statement}[\cite{K}] Over $\Cee$, the space $\fg$ can be
represented as a direct sum of subspaces
\[
\fg=\mathop{\bigoplus}\limits_{\alpha\in \fh^*} \fg_\alpha.
\]
\end{Statement}

Note that $\fh\subsetneq\fg_0$ over $\Kee$, e.g., all weights of the
form $p\alpha$ over $\Cee$ become 0.



\ssec{Systems of simple and positive roots} In this subsection,
$\fg=\fg(A,I)$, and $R$ is the root system of $\fg$.

For any subset $B=\{\sigma_{1}, \dots, \sigma_{m}\} \subset R$, we
set  (we denote by $\Zee_{+}$ the set of non-negative integers):
\[
R_{B}^{\pm} =\{ \alpha \in R \mid \alpha = \pm \sum n_{i}
\sigma_{i},\;\;n_{i} \in \Zee_{+} \}.
\]

The set $B$ is called a {\it system of simple roots~}\index{Root!
simple system of}\; of $R$ (or $\fg$) if $ \sigma_{1}, \dots ,
\sigma_{m}$ are linearly independent and $R=R_B^+\cup R_B^-$. Note
that $R$ contains basis coordinate vectors, and therefore spans
$\Ree^n$; thus, any system of simple roots contains exactly $n$
elements.

Let $(\cdot,\cdot)$ be the standard Euclidean inner product in
$\Ree^n$. A subset $R^+\subset R$ is called a {\it system of
positive roots~}\index{Root! positive system of}\; of $R$ (or $\fg$)
if there exists $x\in\Ree^n$ such that
\begin{equation}\label{x}
\begin{split}
 &(\alpha,x)\in\Ree\backslash \{0\}\text{ for any $\alpha\in R$},\\
 &R^+=\{\alpha\in R\mid (\alpha,x)>0\}.
\end{split}
\end{equation} Since $R$ is a finite (or, at least, countable
if $\dim \fg(A)=\infty$) set, so the set
\[\{y\in\Ree^n\mid\text{there exists $\alpha\in R$ such that }
(\alpha,y)=0\}
\]
is a finite/countable union of $(n-1)$-dimensional subspaces in
$\Ree^n$, so it has zero measure. So for almost every $x$, condition
(\ref{x}) holds.

By construction, any system $B$ of simple roots is contained in
exactly one system of positive roots, which is precisely $R_B^+$.

\sssbegin{Statement} Any finite system $R^+$ of positive roots of
$\fg$ contains exactly one system of simple roots. This system
consists of all the positive roots (i.e., elements of $R^+$) that
can not be represented as a sum of two positive
roots.\end{Statement}

We can not give an {\it a priori} proof of the fact that each set of
all positive roots each of which is not a sum of two other positive
roots consists of linearly independent elements. This is, however,
true for finite dimensional Lie algebras and superalgebras $\fg(A,
I)$ if $p\neq 2$.

\ssec{Normalization convention}\label{normA} Clearly,
\begin{equation}
\label{rescale} \text{the rescaling
$e_i^\pm\mapsto\sqrt{\lambda_i}e_i^\pm$, sends $A$ to $A':=
\diag(\lambda_1, \dots , \lambda_n)\cdot A$.} 
\end{equation}
Two pairs $(A, I)$ and $(A', I')$ are said to be {\it equivalent}
(and we write $(A, I)\sim(A', I')$) if $(A', I')$ is obtained from
$(A, I)$ by a composition of a permutation of parities and a
rescaling $A' = \diag (\lambda_{1}, \dots, \lambda_{n})\cdot A$,
where $\lambda_{1}\dots \lambda_{n}\neq 0$. Clearly, equivalent
pairs determine isomorphic Lie superalgebras.

The rescaling affects only the matrix $A_B$, not the set of parities
$I_B$. The Cartan matrix $A$ is said to be {\it
normalized}\index{Cartan matrix, normalized} if
\begin{equation}
\label{norm} A_{jj}=0\text{~~ or 1, or 2,}
\end{equation}
where we let $A_{jj}=2$ only if $p_j=\ev$; in order to distinguish
between the cases where $p_j=\ev$ and $p_j=\od$, we write
$A_{jj}=\ev$ or $\od$, instead of 0 or 1, if $p_j=\ev$. {\bf We will
only consider normalized Cartan matrices; for them, we do not have
to describe $I$.}

The row with a $0$ or $\ev$ on the main diagonal can be multiplied by
any nonzero factor; usually (not only in this paper) we multiply the
rows so as to make $A_{B}$ symmetric, if possible.

{\it A posteriori}, for each {\bf finite dimensional} Lie
(super)algebra of the form $\fg(A)$ with indecomposable Cartan
matrix $A$, the matrix $A$ is symmetrizable (i.e., it can be made
symmetric by operation (\ref{rescale})) for any $p$. For affine and
almost affine Lie (super)algebra of the form $\fg(A)$ this is not
so, cf. \cite{CCLL}

\ssec{Equivalent systems of simple roots} \label{EqSSR} Let
$B=\{\alpha_1,\dots,\alpha_n\}$ be a system of simple roots. Choose
non-zero elements $e_i^\pm$ in the 1-dimensional (by definition)
superspaces $\fg_{\pm\alpha_i}$; set $h_{i}=[e_{i}^{+}, e_{i}^-]$,
let $A_{B} =(A_{ij})$, where the entries $A_{ij}$ are recovered from
relations \eqref{gArel_0}, and let $I_{B}=\{p(e_{1}), \cdots,
p(e_{n})\}$. Lemma \ref{serg} claims that all the pairs $(A_B,I_B)$
are equivalent to each other.

Two systems of simple roots $B_{1}$ and $B_{2}$ are said to be {\it
equivalent} if the pairs $(A_{B_{1}}, I_{B_{1}})\sim(A_{B_{2}},
I_{B_{2}})$.

For the role of the \lq\lq best'' (first among equals) order of
indices we propose the one that minimizes the value
\begin{equation}\label{minCM}
\max\limits_{i,j\in\{1,\dots,n\}\text{~such that~}(A_B)_{ij}\neq
0}|i-j|
\end{equation}
(i.e., gather the non-zero entries of $A$ as close to the main
diagonal as possible).

\sssec{Chevalley generators and Chevalley bases}\label{SsChev} We
often denote the set of generators corresponding to a normalized
matrix by $X_{1}^{\pm},\dots , X_{n}^{\pm}$ instead of
$e_{1}^{\pm},\dots , e_{n}^{\pm}$; and call them, together with the
elements $H_i:=[X_{i}^{+}, X_{i}^{-}]$, and the derivatives $d_j$
added for convenience for all $i$ and $j$, the {\it Chevalley
generators}.\index{Chevalley generator}

For $p=0$ and normalized Cartan matrices of simple finite
dimensional Lie algebras, there exists only one (up to signs) basis
containing $X_i^\pm$ and $H_i$  in which $A_{ii}=2$ for all $i$ and
all structure constants are integer, cf. \cite{St}. Such a basis is
called the {\it Chevalley}\index{Basis! Chevalley}  basis.

Observe that, having normalized the Cartan matrix of $\fo(2n+1)$ so
that $A_{ii}=2$ for all $i\neq n$ but $A_{nn}=1$, we get {\bf
another} basis with integer structure constants. We think that this
basis also qualifies to be called {\it Chevalley  basis}; for the
Lie superalgebras, the basis normalized as in (\ref{norm}) is even
more appropriate.

\begin{Conjecture} If $p>2$, then for finite dimensional Lie
(super)algebras with indecomposable Cartan matrices normalized as in
$(\ref{norm})$, there also exists only one (up to signs) analog of
the Chevalley basis. \end{Conjecture}

We had no idea how to describe analogs of Chevalley bases for $p=2$
until appearance of the recent paper \cite{CR}; clearly, its methods
should solve the problem.

\section{Ortho-orthogonal and periplectic Lie superalgebras}\label{Soo}

In this section, $p=2$ and $\Kee$ is perfect. We also assume that
$n_\ev,n_\od>0$.

\ssec{Non-degenerate even supersymmetric bilinear forms and
ortho\Defis or\-tho\-gonal Lie superalgebras} For $p=2$, there are,
in general, four equivalence classes of inequivalent non-degenerate
even supersymmetric bilinear forms on a given superspace. Any such
form $B$ on a superspace $V$ of superdimension $n_\ev|n_\od$ can be
decomposed as follows:
\[
B=B_\ev\oplus B_\od,
\]
where $B_\ev$, $B_\od$ are symmetric non-degenerate forms on $V_\ev$
and $V_\od$, respectively. For $i=\ev,\od$, the form $B_i$ is
equivalent to $1_{n_i}$ if $n_i$ is odd, and equivalent to $1_{n_i}$
or $\Pi_{n_i}$ if $n_i$ is even. So every non-degenerate even
symmetric bilinear form is equivalent to one of the following forms
(some of them are defined not for all dimensions):
\[
\begin{array}{ll}
B_{II}=1_{n_\ev}\oplus 1_{n_\od}; &B_{I\Pi}=1_{n_\ev}\oplus
\Pi_{n_\od}\text{~if~}n_\od \text{~is even;}\\
B_{\Pi I}=\Pi_{n_\ev}\oplus 1_{n_\od}\text{~if~}n_\ev \text{~is
even}; &B_{\Pi\Pi}=\Pi_{n_\ev}\oplus \Pi_{n_\od}\text{~if~}n_\ev,
n_\od \text{~are even.}
\end{array}
\]
We denote the Lie superalgebras that preserve the respective forms
by $\fo\fo_{II}(n_\ev|n_\od)$, $\fo\fo_{I\Pi}(n_\ev|n_\od)$,
$\fo\fo_{\Pi I}(n_\ev|n_\od)$, $\fo\fo_{\Pi\Pi}(n_\ev|n_\od)$,
respectively. Now let us describe these algebras.

\sssec{$\fo\fo_{II}(n_\ev|n_\od)$} If $n\geq 3$, then the Lie
superalgebra $\fo\fo_{II}^{(1)}(n_\ev|n_\od)$ is simple. This Lie
superalgebra {\bf has no Cartan matrix}.

\sssec{$\fo\fo_{I\Pi}(n_\ev|n_\od)$~~($n_\od=2k_\od$)} The Lie
superalgebra $\fo\fo_{I\Pi}^{(1)}(n_\ev|n_\od)$ is simple, it has
Cartan matrix if and only if $n_\ev$ is odd; this matrix has the
following form (up to a format; all possible formats ---
corresponding to $\ast=0$ or $\ast=\ev$ --- are described in Table
\ref{tbl} below): \begin{equation}\label{oowith1}
\begin{pmatrix} \ddots&\ddots&\ddots&\vdots\\
\ddots&\ast&1&0\\
\ddots&1&\ast&1\\
\cdots&0&1&1\end{pmatrix}
\end{equation}

\sssec{$\fo\fo_{\Pi\Pi}(n_\ev|n_\od)~~(n_\ev=2k_\ev,
n_\od=2k_\od)$}\label{oo-oo_PP} If $n=n_\ev+n_\od\geq 6$, then
\begin{equation}
\label{oopipi}
\begin{split}
&\text{if $k_\ev+k_\od$ is odd, then the Lie superalgebra
$\fo\fo_{\Pi\Pi}^{(2)}(n_\ev|n_\od)$ is simple;} \cr &\text{if
$k_\ev+k_\od$ is even, then the Lie superalgebra
$\fo\fo_{\Pi\Pi}^{(2)}(n_\ev|n_\od)/\Kee 1_{n_\ev|n_\od}$ is
simple.}
\end{split}
\end{equation}

Each of these simple Lie superalgebras is also close to a Lie
superalgebra with Cartan matrix. To describe this Cartan matrix Lie
superalgebra in most simple terms, we will choose a slightly
different realization of $\fo\fo_{\Pi\Pi}(2k_\ev|2k_\od)$: Let us
consider it as the algebra of linear transformations that preserve
the bilinear form $\Pi(2k_\ev+2k_\od)$ in the format
$k_\ev|k_\od|k_\ev|k_\od$. Then the algebra
$\fo\fo_{\Pi\Pi}^{(i)}(2k_\ev|2k_\od)$ is spanned by supermatrices
of format $k_\ev|k_\od|k_\ev|k_\od$ and of the form
\begin{equation}\label{matform}
\begin{pmatrix}A&C\\D&A^T\end{pmatrix}\text{~where~}
\begin{array}{l}
A\in\begin{cases}\fgl(k_\ev|k_\od)&\text{if~}i\leq
1,\\\fsl(k_\ev|k_\od)&\text{if~}i\geq 2,\end{cases} \\ 
C,D\text{~are~}\begin{cases}\text{symmetric matrices}&\text{if~}
i=0;\\\text{symmetric zero-diagonal matrices}&\text{if~} i\geq
1.\end{cases}\end{array}
\end{equation}
If $i\geq 1$, these derived algebras have a non-trivial central
extension given by the following cocycle:
\begin{equation}\label{cocycle}
F\left(\begin{pmatrix}A&C\\D&A^T\end{pmatrix},
\begin{pmatrix}A'&C'\\D'&A'^T\end{pmatrix}\right)=\sum\limits_{1\leq
i<j\leq k_\ev+k_\od} (C_{ij}D'_{ij}+C'_{ij}D_{ij})
\end{equation}
(note that this expression resembles $\frac 12\tr(CD'+C'D)$). We
will denote this central extension of
$\fo\fo_{\Pi\Pi}^{(i)}(2k_\ev|2k_\od)$ by
$\fo\fo\fc(i,2k_\ev|2k_\od)$.

Let\index{$I_0:=\diag(1_{k_\ev\vert k_\od},0_{k_\ev\vert k_\od})$}
\begin{equation}
\label{I_0osp}I_0:=\diag(1_{k_\ev|k_\od},0_{k_\ev|k_\od}).
\end{equation}
Then the corresponding Cartan matrix Lie superalgebra is
\begin{equation}
\label{ooc}
\begin{split}
&\fo\fo\fc(2,2k_\ev|2k_\od)\subplus\Kee I_0\quad\text{ if
$k_\ev+k_\od$ is odd;} \cr &\fo\fo\fc(1,2k_\ev|2k_\od)\subplus\Kee
I_0\quad\text{ if $k_\ev+k_\od$ is even.}
\end{split}
\end{equation}

The corresponding Cartan matrix has the following form (up to a
format; all possible formats --- corresponding to $\ast=0$ or
$\ast=\ev$ --- are described in Table \ref{tbl} below):
\begin{equation}\label{ooPPCM}
\begin{pmatrix}
\ddots&\ddots&\ddots&\vdots&\vdots\\
\ddots&\ast&1&0&0\\
\ddots&1&\ast&1&1\\
\cdots&0&1&\ev&0\\
\cdots&0&1&0&\ev\end{pmatrix}
\end{equation}

\ssec{The non-degenerate odd supersymmetric bilinear forms.
Periplectic Lie superalgebras} \label{peLS} In this subsect., $m\geq
3$.

\begin{equation}
\label{pe0}
\begin{split}
&\text{If $m$ is odd, then the Lie superalgebra $\fpe_B^{(2)}(m)$ is
simple;} \cr &\text{If $m$ is even, then the Lie superalgebra
$\fpe_B^{(2)}(m)/\Kee 1_{m|m}$ is simple.}
\end{split}
\end{equation}

If we choose the form $B$ to be $\Pi_{m|m}$, then the algebras
$\fpe_B^{(i)}(m)$ consist of matrices of the form (\ref{matform});
the only difference from $\fo\fo_{\Pi\Pi}^{(i)}$ is the format which
in this case is $m|m$.

Each of these simple Lie superalgebras has a $2$-structure. Note
that if $p\neq 2$, then the Lie superalgebra $\fpe_B(m)$ and its
derived algebras are not close to Cartan matrix Lie superalgebras
(because, for example, their root system is not symmetric). If $p=2$
and $m\geq 3$, then they {\bf are} close to Cartan matrix Lie
superalgebras; here we describe them.

The algebras $\fpe_B^{(i)}(m)$, where $i>0$, have non-trivial
central extensions with cocycles (\ref{cocycle}); we denote these
central extensions by $\fpe\fc(i,m)$. Let us introduce another
matrix \index{$I_0:=\diag(1_m,0_m)$}
\begin{equation}
\label{I_0pe}I_0:=\diag(1_m,0_m).
\end{equation} Then the Cartan matrix Lie
superalgebras are
\begin{equation}
\label{pec}
\begin{split}
&\fpe\fc(2,m)\subplus\Kee I_0\text{ if $m$ is odd;} \cr
&\fpe\fc(1,m)\subplus\Kee I_0\text{ if $m$ is even.}
\end{split}
\end{equation}

The corresponding Cartan matrix has the form (\ref{ooPPCM}); the
only condition on its format is that the last two simple roots must
have distinct parities. The corresponding Dynkin diagram is shown in
Table \ref{tbl}; all its nodes, except for the \lq\lq horns'', may
be both $\otimes$ or~$\odot$, see (\ref{cm1}).

\ssec{Superdimensions} The following expressions (with a $+$ sign)
are the superdimensions of the
 relatives of the ortho-orthogonal and
periplectic Lie superalgebras that possess Cartan matrices. To get
the superdimensions of the simple relatives, one should replace $+2$
and $+1$ by $-2$ and $-1$, respectively, in the two first lines and
the four last ones:
\begin{equation}\label{dimtable}
\begin{array}{lll}
\dim \fo\fc (1;2k)\subplus\Kee I_0&=2k^2-k\pm 2&\text{if $k$ is even;}\\
\dim \fo\fc (2;2k)\subplus\Kee I_0&=2k^2-k\pm 1&\text{if $k$ is odd;}\\
\dim\fo^{(1)}(2k+1)&=2k^2+k&\\
\sdim\fo\fo^{(1)}(2k_\ev+1|2k_\od)&=2k_\ev^2+k_\ev+
2k_\od^2+k_\od\mid
2k_\od(2k_\ev+1)&\\
\sdim \fo\fo\fc (1;2k_\ev|2k_\od)\subplus\Kee I_0&=2k_\ev^2-k_\ev+
2k_\od^2-k_\od\pm 2\mid 4k_\ev k_\od&\text{if $k_\ev+k_\od$ is even;}\\
\sdim \fo\fo\fc (2;2k_\ev|2k_\od)\subplus\Kee I_0&=2k_\ev^2-k_\ev+
2k_\od^2-k_\od\pm 1\mid 4k_\ev k_\od&\text{if $k_\ev+k_\od$ is odd;}\\
\sdim \fpe\fc (1;m)\subplus\Kee I_0&=m^2\pm 2\mid m^2-m&\text{if $m$ is even;}\\
\sdim \fpe\fc (2;m)\subplus\Kee I_0&=m^2\pm 1\mid m^2-m&\text{if $m$
is odd}
\end{array}
\end{equation}

\sssec{Summary: The types of Lie superalgebras preserving
non-degene\-ra\-te symmetric forms} In addition to the isomorphisms
$\foo_{\Pi I}(a|b)\simeq\foo_{I\Pi }(b|a)$, there is the only \lq\lq
occasional" isomorphism intermixing the types of Lie superalgebras
preserving non-degene\-ra\-te symmetric forms:
$\foo_{\Pi\Pi}^{(1)}(6|2)\simeq\fpe^{(1)}(4)$.

Let $\widehat{\fg}:=\fg\subplus\Kee I_0$.  We have the following
types of non-isomorphic Lie (super)algebras:
\begin{equation}\label{oandoo}\renewcommand{\arraystretch}{1.4}
\begin{tabular}{|l|l|}
\hline no relative has Cartan matrix&with Cartan matrix\\
\hline
$\foo_{II}(2n+1|2m+1),\;\;\foo_{II}(2n+1|2m)$&$\widehat{\fo\fc(i;2n)},\;\;\fo^{(1)}(2n+1);\;
\;\widehat{\fpe\fc(i;k)}$\\
$\foo_{II}(2n|2m),
\;\;\foo_{I\Pi}(2n|2m);\;\;\fo_I(2n);$&$\widehat{\foo\fc(i;2n|2m)},
\;\;\foo_{I\Pi}^{(1)}(2n+1|2m)$\\
\hline\end{tabular}
\end{equation}

\section{Dynkin diagrams}

A usual way to represent simple Lie algebras over $\Cee$ with
integer Cartan matrices is via graphs called, in the finite
dimensional case, {\it Dynkin diagrams}. The Cartan matrices of
certain interesting infinite dimensional simple Lie {\it
super}algebras $\fg$ (even over $\Cee$) can be non-symmetrizable or
(for any $p$ in the super case and for $p>0$ in the non-super case)
have entries belonging to the ground field $\Kee$. Still, it is
always possible to assign an analog of the Dynkin diagram to each
(modular) Lie (super)algebra with Cartan matrix, of course) provided
the edges and nodes of the graph (Dynkin diagram) are rigged with an
extra information. Although these analogs of the Dynkin graphs are
not uniquely recovered from the Cartan matrix (and the other way
round), they give a graphic presentation of the Cartan matrices and
help to observe some hidden symmetries.

Namely, the {\it Dynkin diagram}\index{Dynkin diagram} of a
normalized $n\times n$ Cartan matrix $A$ is a set of $n$ nodes
connected by multiple edges, perhaps endowed with an arrow,
according to the usual rules (\cite{K}) or their modification, most
naturally formulated by Serganova: compare \cite{Se, FLS} with
\cite{FSS}. In what follows, we recall these rules, and further
improve them to fit the modular case.

\ssec{Nodes} To every simple root there corresponds
\begin{equation}\label{cm1}
\begin{cases}
\text{a node}\; \mcirc\; &\text{if $p(\alpha_{i})= \ev$ and $
A_{ii}=2$},\\
\text{a node}\; \ast \;&\text{if $p(\alpha_{i}) =\ev$ and
$A_{ii}=\od$};\\
\text{a node}\; \mbullet \;&\text{if $p(\alpha_{i}) =\od$ and
$A_{ii}=1$};\\
\text{a node}\; \motimes \;& \text{if $p(\alpha_{i}) =\od$ and $
A_{ii}=0$},\\
\text{a node}\; \odot\; &\text{if $p(\alpha_{i})= \ev$ and $
A_{ii}=\ev$}.\\
\end{cases}
\end{equation}

The Lie algebras $\fsl(2)$ and $\fo(3)^{(1)}$ with Cartan matrices
$(2)$ and $(\od)$, respectively, and the Lie superalgebra
$\fosp(1|2)$ with Cartan matrix $(1)$ are simple.

The Lie algebra $\fgl(2)$ with Cartan matrix $(\ev)$ and the Lie
superalgebra $\fgl(2|2)$ with Cartan matrix $(0)$ are solvable of
$\dim 4$ and $\sdim 2|2$, respectively. Their derived algebras are
the {\it Heisenberg algebra} $\fhei(2):=\fhei(2|0)\simeq\fsl(2)$ and
the {\it Heisenberg superalgebra} $\fhei(0|2)\simeq\fsl(1|1)$ of
(super)dimension 3 and $1|2$, respectively.

\sssbegin{Remark} {\it A posteriori} (from the classification of
simple Lie superalgebras with Cartan matrix and of polynomial
growth) we find out that {\bf for $p=0$}, the simple root~$\odot$
can only occur if $\fg(A, I)$ grows faster than polynomially. Thanks
to classification again, if $\dim \fg<\infty$, the simple root
$\odot$ can not occur if $p>3$; whereas for $p=3$, the Brown Lie
algebras are examples of $\fg(A)$ with a simple root of type
$\odot$; for $p=2$, such roots are routine.
\end{Remark}

\ssec{Edges} If $p=2$ and $\dim \fg(A)<\infty$, the Cartan matrices
considered are symmetric. If $A_{ij}=a$, where $a\neq 0$ or 1, then
we rig the edge connecting the $i$th and $j$th nodes by a label $a$.

If $p>2$ and $\dim \fg(A)<\infty$, then $A$ is symmetrizable, so let
us symmetrize it, i.e., consider $DA$ for an invertible diagonal
matrix $D$. Then, if $(DA)_{ij}=a$, where $a\neq 0$ or $-1$, we rig
the edge connecting the $i$th and $j$th nodes by a label $a$.

If all off-diagonal entries of $A$ belong to $\Zee/p$ and their
representatives are selected to be non-positive integers, we can
draw the Dynkin diagram as for $p=0$, i.e., connect the $i$th node
with the $j$th one by $\max(|A_{ij}|, |A_{ji}|)$ edges rigged with
an arrow $>$ pointing from the $i$th node to the $j$th if
$|A_{ij}|>|A_{ji}|$ or in the opposite direction if
$|A_{ij}|<|A_{ji}|$.

\ssec{Reflections} Let $R^+$ be a system of positive roots of Lie
superalgebra $\fg$, and let $B=\{\sigma_1,\dots,\sigma_n\}$ be the
corresponding system of simple roots with some corresponding pair
$(A=A_B,I=I_B)$. Then the set
$(R^+\backslash\{\sigma_k\})\coprod\{-\sigma_k\}$ is a system of
positive roots for any $k\in \{1, \dots, n\}$. This operation is
called {\it the reflection in $\sigma_k$}; it changes the system of
simple roots by the formulas
\begin{equation}
\label{oddrefl}
r_{\sigma_k}(\sigma_{j})= \begin{cases}{-\sigma_j}&\text{if~}k=j,\\
\sigma_j+B_{kj}\sigma_k&\text{if~}k\neq j,\end{cases}\end{equation}
where
\begin{equation}
\label{Boddrefl}B_{kj}=\begin{cases}
-\displaystyle\frac{2A_{kj}}{A_{kk}}& \text{~if~}p_k=\ev,\;\;
A_{kk}\neq 0,\text{~and~}
-\displaystyle\frac{2A_{kj}}{A_{kk}}\in \Zee/p\Zee,\\
p-1&\text{~if~}p_k=\ev,\;\; A_{kk}\neq 0\text{~and~}
 -\displaystyle\frac{2A_{kj}}{A_{kk}}\not\in \Zee/p\Zee,\\
-\displaystyle\frac{A_{kj}}{A_{kk}}& \text{~if~}p_k=\od,\;\;
A_{kk}\neq 0,\text{~and~}
-\displaystyle\frac{A_{kj}}{A_{kk}}\in \Zee/p\Zee,\\
p-1&\text{~if~}p_k=\od,\;\; A_{kk}\neq 0,
\text{~and~} -\displaystyle\frac{A_{kj}}{A_{kk}}\not\in \Zee/p\Zee,\\
1&\text{~if~}p_k=\od,\;\; A_{kk}=0,\;\;A_{kj}\neq 0,\\
0&\text{~if~}p_k=\od,\;\; A_{kk}=A_{kj}=0,\\
p-1&\text{~if~}p_k=\ev,\;\; A_{kk}=\ev,\;\;A_{kj}\neq 0,\\
0&\text{~if~}p_k=\ev,\;\; A_{kk}=\ev,\;\;A_{kj}=0,\end{cases}
\end{equation}
where we consider $\Zee/p\Zee$ as a subfield of $\Kee$.

\sssbegin{Remark} The description of the numbers $B_{ik}$ is
empirical and based on classification \cite{BGL}: For
infinite-dimensional Lie (super)algebras these numbers might be
different. In principle, in the second, fourth and penultimate
cases, the matrix \eqref{Boddrefl} can be equal to $kp-1$ for any
$k\in\Nee$, and in the last case any element of $\Kee$ may occur.
For $\dim\fg<\infty$, this does do not happen (and it is of interest
to investigate at least the simplest infinite dimensional case ---
the modular analog of \cite{CCLL}).

The values $-\displaystyle\frac{2A_{kj}}{A_{kk}}$ and
$-\displaystyle\frac{A_{kj}}{A_{kk}}$ are elements of $\Kee$, while
the roots are elements of a vector space over $\Ree$. Therefore {\bf
These expressions in the first and third cases in} (\ref{Boddrefl})
{\bf should be understood as} \lq\lq {\bf the minimal non-negative
integer congruent to $-\displaystyle\frac{2A_{kj}}{A_{kk}}$ or
$-\displaystyle\frac{A_{kj}}{A_{kk}}$, respectively''. (If
$\dim\fg<\infty$, these expressions are always congruent to
integers.)}

{\bf There is known just one exception: If $p=2$ and
$A_{kk}=A_{jk}$, then} $-\displaystyle\frac{2A_{jk}}{A_{kk}}$ {\bf
should be understood as $2$, not 0.}
\end{Remark}

The name \lq\lq reflection'' is used because in the case of
(semi)simple finite-dimensional Lie algebras this action extended on
the whole $R$ by linearity is a map from $R$ to $R$, and it does not
depend on $R^+$, only on $\sigma_k$. This map is usually denoted by
$r_{\sigma_k}$ or just $r_{k}$. The map $r_{\sigma_i}$ extended to
the $\Ree$-span of $R$ is reflection in the hyperplane orthogonal to
$\sigma_i$ relative the bilinear form dual to the Killing form.

The reflections in the even (odd) roots are said to be {\it even}
({\it odd}).\index{Reflection! odd}\index{Reflection! even!
non-isotropic} \index{Reflection! even! isotropic} A simple root is
called {\it isotropic}, if the corresponding row of the Cartan
matrix has zero on the diagonal, and {\it non-isotropic} otherwise.
The reflections that correspond to isotropic or non-isotropic roots
will be referred to accordingly.

If there are isotropic simple roots, the reflections $r_\alpha$ do
not, as a rule, generate a version of the {\it Weyl group} because
the product of two reflections in nodes not connected by one
(perhaps, multiple) edge is not defined. These reflections just
connect pair of \lq\lq neighboring'' systems of simple roots and
there is no reason to expect that we can multiply two distinct such
reflections. In the general case (of Lie superalgebras and $p>0$),
the action of a given isotropic reflections (\ref{oddrefl}) can not,
generally, be extended to a linear map $R\tto R$. For Lie
superalgebras over $\Cee$, one can extend the action of reflections
by linearity to the root lattice but this extension preserves the
root system only for $\fsl(m|n)$ and $\fosp(2m+1|2n)$, cf.
\cite{Se1}.

If $\sigma_i$ is an odd isotropic root, then the corresponding
reflection sends one set of Chevalley generators into a new one:
\begin{equation}
\label{oddrefx} \tilde X_{i}^{\pm}=X_{i}^{\mp};\;\; \tilde
X_{j}^{\pm}=\begin{cases}[X_{i}^{\pm},
X_{j}^{\pm}]&\text{if $A_{ij}\neq 0, \ev$},\\
X_{j}^{\pm}&\text{otherwise}.\end{cases}
\end{equation}

\sssec{On neighboring root systems} Serganova \cite{Se}
proved (for $p=0$) that there is always a chain of reflections
connecting $B_1$ with some system of simple roots $B'_2$ equivalent
to $B_2$ in the sense of definition \ref{EqSSR}. Here is the modular
version of Serganova's Lemma. Observe that Serganova's statement is
not weaker: Serganova used only odd reflections.

\begin{Lemma}[\cite{LCh}]\label{serg} For any two systems of simple roots
$B_1$ and $B_2$ of any simple finite dimensional Lie superalgebra
with Cartan matrix, there is always a chain of reflections
connecting $B_1$ with $B_2$.\end{Lemma}

\section{Presentations of $\fg(A)$}\label{Spresent}

\ssec{Serre relations, see \cite{GL2}}

Let $A$ be an $n\times n$ matrix. We find the defining relations by
induction on $n$ with the help of the Hochschild--Serre spectral
sequence (for its description for Lie superalgebras, which has
certain subtleties, see \cite{Po}). For the basis of induction
consider the following cases of Dynkin diagrams with one vertex and
no edges:
\begin{equation}
\label{3.4.1}
\begin{array}{lll}
 \mcirc \text{ or } \mbullet &\text{no relations, i.e., $\fg^{\pm}$
are free Lie superalgebras}&\text{if $p\neq 3$;}\\
\mbullet &\ad_{X^{\pm}}^2(X^{\pm})=0&\text{if $p=3$;}\\
{\motimes}& [X^{\pm}, X^{\pm}]=0.&
\end{array}
\end{equation}
Set $\deg X_{i}^{\pm} = 0$ for $ 1\leq i\leq n-1$ and $\deg
X_{n}^{\pm} = \pm 1$. Let $\fg^{\pm} = \oplus \fg_{ i}^{\pm}$ and
$\fg= \oplus \fg_{i} $ be the corresponding $\Zee $-gradings. Set
$\fg_{\pm} =\fg^{\pm}/\fg_{ 0}^{\pm}$. From the Hochschild--Serre
spectral sequence for the pair $\fg_{ 0}^{\pm} \subset \fg ^{\pm}$
we get (for more detail, see \cite{LCh}):
\begin{equation}
\label{3.4.2} H_{2}(\fg_{\pm})\subset H_{2}(\fg_{0}^{\pm})\oplus
H_{1}(\fg_{ 0}^{\pm}; H_{ 1}(\fg_{\pm}))\oplus H_{0}(\fg_ {
0}^{\pm}; H_{ 2}(\fg_{\pm})).
\end{equation}
It is clear that
\begin{equation}
\label{3.4.3} H_{1}(\fg_{\pm})= \fg_{1}^{\pm} , \;\;\; H_{
2}(\fg_{\pm}) = \wedge^{2}(\fg_ {1}^{\pm})/\fg_{ 2}^{\pm}.
\end{equation}
So, the second summand in~(\ref{3.4.2}) provides us with relations
of the form:
\begin{equation}
\label{3.4.4} \begin{array}{ll} (\ad_{X_{n}^{\pm}})^{k_{ni}}
(X_{i}^{\pm})=0&\text{if the $n$-th simple root is not}\;\;
{\motimes}\\
{} [X_{n}, X_{n}]=0&\text{if the $n$-th simple root is}
\;\;\motimes.
\end{array}
\end{equation}
while the third summand in (\ref{3.4.2}) is spanned by the
$\fg_{0}^{\pm}$-lowest vectors in
\begin{equation}
\label{3.4.6} \wedge^{2}(\fg_{1}^{\pm})/(\fg_{2}^{\pm} + \fg^{\pm}
\wedge^{2}(\fg_{1}^{\pm})).
\end{equation}

Let the matrix $B=(B_{ij})$ be as in formula (\ref{Boddrefl}). The
following proposition, whose proof is straightforward, illustrates
the usefulness of our normalization of Cartan matrices as compared
with other options:

\sssbegin{Proposition} The numbers $k_{in}$ and $k_{ni}$ in
$(\ref{3.4.4})$ are expressed in terms of $(B_{ij})$ as follows:
\begin{equation}
\label{srpm}
\begin{matrix} (\ad_{X_{i}^{\pm}})^{1+B_{ij}}(X_{j}
^{\pm})=0 &
\text{ for $ i \neq j$} \\ & \\
{}[X_{i}^{\pm}, X_{i}^{\pm}]=0& \text{if $A_{ii} =0$ and $p\neq
2$}\\
\left(X_{i}^{\pm}\right)^2=0& \text{if $A_{ii} =0$ and $p=
2$}.\end{matrix}
\end{equation}
\end{Proposition}
The relations (\ref{gArel_0}) and $(\ref{srpm})$ will be called {\it
Serre relations} for the Lie superalgebra $\fg(A)$. If $p=3$, then
the relation
\begin{equation}
\label{x^3} [X_{i}^{\pm}, [X_{i}^{\pm}, X_{i}^{\pm}]]=0~~
\;\;\text{for $X_{i}^{\pm}$ odd and $A_{ii} =1$}
\end{equation}
is not a consequence of the Jacobi identity; for simplicity,
however, we will include it in the set of Serre relations. {\bf
Observe that usually only the relations of the first line in
$(\ref{srpm})$ are said to be} {\it Serre relations} {\bf for the
Lie superalgebra $\fg(A)$}.

\ssec{Non-Serre relations} These are relations that correspond to
the third summand in~(\ref{3.4.2}). Let us consider the simplest
case: $\fsl (m|n)$ in the realization with the system of simple
roots
\begin{equation}
\label{circ}
 \xymatrix@C=1em{
 \ffcirc\ar@{-}[r]&\dots\ar@{-}[r]&\ffcirc\ar@{-}[r]&
 \ffotimes\ar@{-}[r]&\ffcirc\ar@{-}[r]&\dots\ar@{-}[r]&
 \ffcirc
 }
\end{equation}

Then $H_2(\fg_{\pm})$ from the third summand in (\ref{3.4.2}) is
just $\wedge^2(\fg_{\pm})$. For simplicity, we confine ourselves to
the positive roots. Let $X_{1}$, \dots , $X_{m-1}$ and
$Y_{1}$,~\dots,~$Y_{n-1}$ be the root vectors corresponding to even
roots separated by the root vector $Z$ corresponding to the root
$\otimes$.

If $n=1$ or $m=1$, then $\wedge^2(\fg)$ is an irreducible $\fg_{\bar
0}$-module and there are no non-Serre relations. If $n\neq 1$ and
$m\neq 1$, then $\wedge^2(\fg)$ splits into 2 irreducible $\fg_{\bar
0}$-modules. The lowest component of one of them corresponds to the
relation $[Z, Z]=0$, the other one corresponds to the non-Serre-type
relation
\begin{equation}
\label{*} [[X_{m-1}, Z], [Y_{1}, Z]] =0. \end{equation}

If, instead of $\fsl (m|n)$, we would have considered the Lie
algebra $\fsl(m+n)$, the same argument would have led us to the two
relations, both of Serre type:
\[
\ad_Z^2(X_{m-1})=0, \qquad \ad_Z^2(Y_{1})=0.
\]

In what follows we give an explicit description of the defining
relations in terms of the Chevalley generators of the Lie
(super)algebras of the form $\fg(A)$ or their simple subquotients
$\fg^{(1)}(A)/\fc$.

\section{The Lie (super)algebras of the
form $\fg(A)$ or their simple subquotients $\fg^{(1)}(A)/\fc$}\label{Ssteps}

\ssec{Over $\Cee$} Kaplansky was the first (see his newsletters in
\cite{Kapp}) to discover the exceptional algebras $\fag(2)$ and
$\fab(3)$ (he dubbed them $\Gamma_2$ and $\Gamma_3$, respectively)
and a parametric family $\fosp(4|2; \alpha)$ (he dubbed it
$\Gamma(A, B, C))$); our notation reflect the fact that
$\fag(2)_\ev=\fsl(2)\oplus\fg(2)$ and
$\fab(3)_\ev=\fsl(2)\oplus\fo(7)$ ($\fo(7)$ is $B_3$ in Cartan's
nomenclature).

\begin{figure}[ht]\centering
\includegraphics{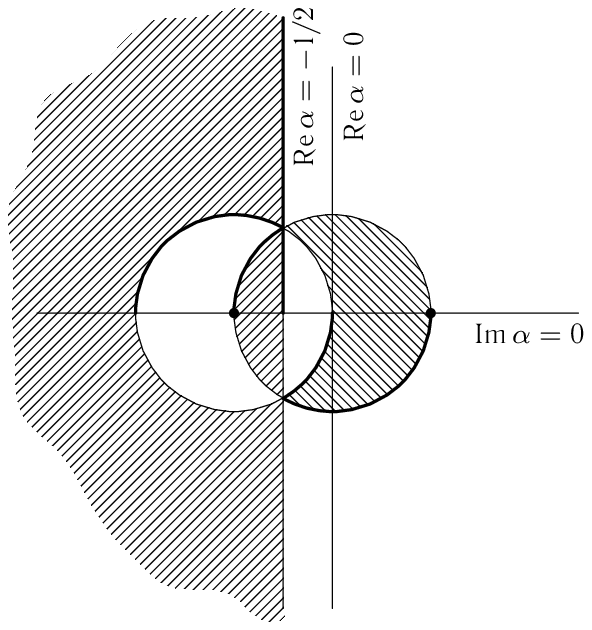}
\caption{}
\end{figure}\label{fig1}

\noindent Kaplansky's description (irrelevant to us at the moment
except for the fact that $A$, $B$ and $C$ are on equal footing) of
what we now identify as $\fosp(4|2; \alpha)$, a parametric family of
deforms of $\fosp(4|2)$, made
 an $S_3$-symmetry of the parameter manifest (to A.\ A.~Kirillov,
 and he informed us, in 1976).
Indeed, since $A+B+C=0$, and $\alpha\in \Cee\cup\infty$ is the ratio
of the two remaining parameters, we get an $S_3$-action on the plane
$A+B+C=0$ which in terms of $\alpha$ is generated by the
transformations:
\begin{equation}\label{osp42symm}
\alpha\longmapsto -1-\alpha, \qquad \alpha\longmapsto
\frac{1}{\alpha}.
\end{equation}
This symmetry should have immediately sprang to mind since
$\fosp(4|2; \alpha)$ is strikingly similar to $\fwk(3; a)$ found 5
years earlier, cf. (\ref{wkiso}), and since $S_3\simeq \SL(2;
\Zee/2)$.

Figure \ref{fig1} depicts the fundamental domains of the
$S_3$-action. The other transformations generated by
(\ref{osp42symm}) are
\[\alpha\longmapsto -\frac{1+\alpha}{\alpha},\quad\alpha\longmapsto
-\frac{1}{\alpha+1},\quad \alpha\longmapsto-\frac{\alpha}{\alpha+1}.
\]

\sssec{Notation: On matrices with a \lq\lq --" sign and other
notation in the lists of inequivalent Cartan matrices}\label{recmat}
The rectangular matrix at the beginning of each list of inequivalent
Cartan matrices for each Lie superalgebra shows the result of odd
reflections (the number of the row is the number of the Cartan
matrix in the list below, the number of the column is the number of
the root (given by small boxed number) in which the reflection is
made; the cells contain the results of reflections (the number of
the Cartan matrix obtained) or a \lq\lq --" if the reflection is not
appropriate because $A_{ii}\neq 0$. Some of the Cartan matrices thus
obtained are equivalent, as indicated.

The number of the matrix $A$ such that $\fg(A)$ has only one odd
simple root is \boxed{boxed}, that with all simple roots odd is {\bf
underlined}. The nodes are numbered by small boxed numbers; the
curly lines with arrows depict odd reflections.

\sssec{Cartan matrices}\label{agab} Recall that $\fag(2)$ of $\sdim
= 17|14$ has the following Cartan matrices\footnotesize
\begin{figure}[ht]\centering
\begin{minipage}[m]{0.2\linewidth}\centering
$
\begin{pmatrix}
 2 & - & - \\
 1 & 3 & - \\
 - & 2 & 4 \\
 - & - & 3
\end{pmatrix}
$
\end{minipage}
\begin{minipage}[m]{0.48\linewidth}\centering
\includegraphics{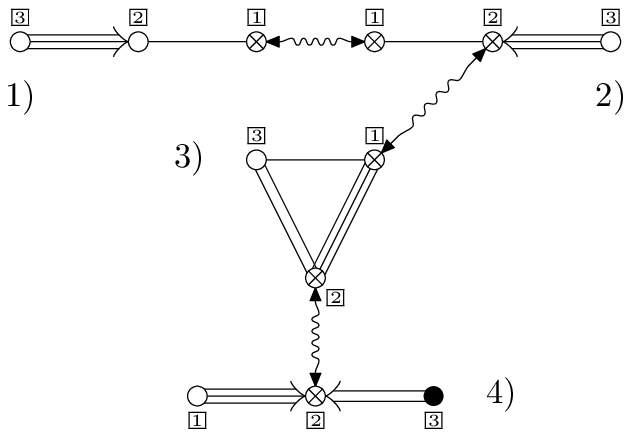}
\end{minipage}\hfill
\end{figure}
\footnotesize\begin{equation}\label{ag2cm}
\boxed{1)}\; \begin{pmatrix} 0 & -1 & 0 \\ -1 & 2 & -3 \\ 0 & -1 & 2
\end{pmatrix}\quad 2)\; \begin{pmatrix}
0 & -1 & 0 \\ -1 & 0 & 3 \\ 0 & -1 & 2
\end{pmatrix}\quad 3)\; \begin{pmatrix}
0 & -3 & 1 \\ -3 & 0 & 2 \\ -1 & -2 & 2
\end{pmatrix}\quad 4)\;
\begin{pmatrix}
2 & -1 & 0 \\ -3 & 0 & 2 \\ 0 & -1 & 1
\end{pmatrix}
\end{equation}\normalsize


Recall that $\fab(3)$ of $\sdim = 24|16$ has the following Cartan
matrices
\begin{figure}[ht]\centering
\begin{minipage}[m]{0.2\linewidth}\centering
$
\begin{pmatrix}
 - & 2 & - & - \\
 3 & 1 & 4 & - \\
 2 & - & - & - \\
 - & - & 2 & 5 \\
 - & 6 & - & 4 \\
 - & 5 & - & -
\end{pmatrix}
$
\end{minipage}
\begin{minipage}[m]{0.58\linewidth}\centering
\includegraphics{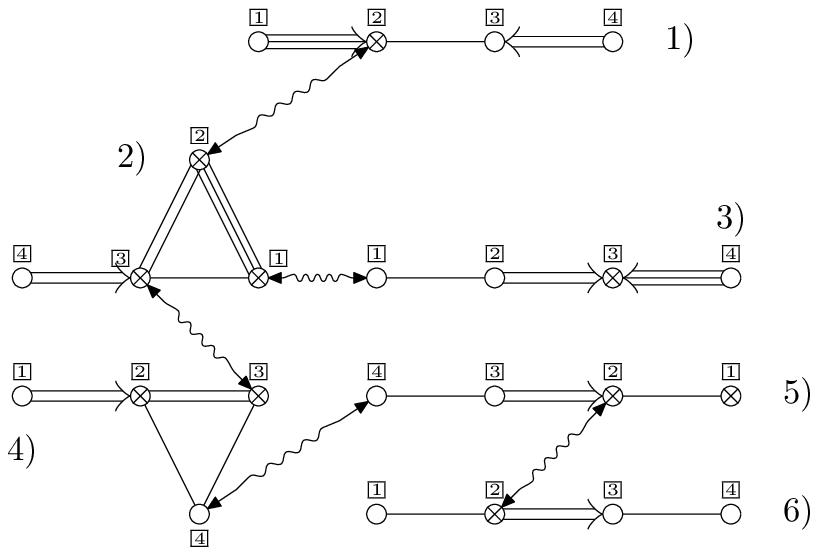}
\end{minipage}\hfill
\end{figure}
\footnotesize\begin{equation}\label{ab3cm}
\begin{matrix}
\boxed{1)}\; \begin{pmatrix}
2 & -1 & 0 & 0 \\ -3 & 0 & 1 & 0 \\ 0 & -1 & 2 & -2 \\
 0 & 0 & -1 & 2
\end{pmatrix}\quad 2)\;
\begin{pmatrix}
0 & -3 & 1 & 0 \\ -3 & 0 & 2 & 0 \\ 1 & 2 & 0 & -2 \\
0 & 0 & -1 & 2
\end{pmatrix}\quad \boxed{3)}\;
\begin{pmatrix}
2 & -1 & 0 & 0 \\ -1 & 2 & -1 & 0 \\ 0 & -2 & 0 & 3 \\
0 & 0 & -1 & 2
\end{pmatrix}\\
4)\; \begin{pmatrix}
2 & -1 & 0 & 0 \\ -2 & 0 & 2 & -1 \\ 0 & 2 & 0 & -1 \\
0 & -1 & -1 & 2
\end{pmatrix}\quad 5)\;
\begin{pmatrix}
0 & 1 & 0 & 0 \\ -1 & 0 & 2 & 0 \\ 0 & -1 & 2 & -1 \\
0 & 0 & -1 & 2
\end{pmatrix}\quad \boxed{6)}\;
\begin{pmatrix}
2 & -1 & 0 & 0 \\ -1 & 2 & -1 & 0 \\ 0 & -2 & 2 & -1 \\ 0 & 0 & -1 &
0
\end{pmatrix}
\end{matrix}
\end{equation}

\normalsize

\clearpage

\ssec{Modular Lie algebras and Lie superalgebras}

\sssec{$p=2$, Lie algebras} Weisfeiler and Kac \cite{WK}
discovered two new parametric families that we denote $\fwk(3;a)$
and $\fwk(4;a)$ ({\it Weisfeiler and Kac
algebras}).\index{$\fwk(3;a)$, Weisfeiler and Kac algebra}
\index{$\fwk(3;a)$, Weisfeiler and Kac algebra}\index{Weisfeiler and
Kac algebra}

$\fwk(3;a)$, where $a\neq 0, -1$, of dim 18 is a non-super version
of $\fosp(4|2; a)$ (although no $\fosp$ exists for $p=2$); the
dimension of its simple subquotient $\fwk(3;a)^{(1)}/\fc$ is equal
to 16; the inequivalent Cartan matrices are:
\[
1)\; \begin{pmatrix} \ev &a &0\\
a&\overline{0}&1\\0&1&\overline{0} \end{pmatrix} ,\;
2)\; \begin{pmatrix} \ev &1+a &a\\
1+a&\overline{0}&1\\
a&1&\overline{0} \end{pmatrix}
\]

$\fwk(4;a)$, where $a\neq 0, -1$, of $\dim=34$; the inequivalent
Cartan matrices are:
\[
1)\; \begin{pmatrix} \ev &a &0&0\\
a &\overline{0}&1&0\\
0&1&\overline{0}&1\\
0&0&1&\overline{0} \end{pmatrix} ,\;
2)\; \begin{pmatrix} \ev &1 &1+a&0\\
1 &\overline{0}& a & 0\\
a+1& a &\overline{0}&a\\
0&0&a&\overline{0} \end{pmatrix}
 ,\;
3)\;\begin{pmatrix} \ev &a & 0 &0\\
a &\overline{0}& a+1 & 0\\
0& a+1 &\overline{0}&1\\
0&0&1&\overline{0} \end{pmatrix}
\]

Weisfeiler and Kac investigated also which of these algebras are
isomorphic and the answer is as follows:
\begin{equation}\label{wkiso}
\renewcommand{\arraystretch}{1.4}\begin{array}{l}
\fwk(3;a)\simeq \fwk(3;a')\Longleftrightarrow
a'=\displaystyle\frac{\alpha a+\beta}{\gamma a+\delta},\text{where
$\begin{pmatrix}\alpha&\beta\\ \gamma &\delta\end{pmatrix}\in \SL(2; \Zee/2)$}\\
\fwk(4;a)\simeq \fwk(4;a')\Longleftrightarrow
a'=\displaystyle\frac{1}{a}.
\end{array}
\end{equation}

\sssec{$p=2$, Lie superalgebras} The same Cartan matrices as for
$\fwk$ algebras but with arbitrary distribution of 0's on the main
diagonal correspond to Lie superalgebras $\fbgl(3;a)$ and
$\fbgl(4;a)$ discovered in \cite{BGL}. The conditions when they are
isomorphic are the same as in \eqref{wkiso}, they have the same
inequivalent Cartan matrices, and are considered also only if $a\neq
0, 1$ (since otherwise they are not simple). We have $\sdim
\fbgl(3;a)=10/8|8$ and $\sdim \fbgl(4;a)=18|16$.

\begin{landscape}
{\bf $p=2$, Lie superalgebras. Dynkin diagrams for
$p=2$}\label{tbl}{}~{}


\begin{center}\label{tbll}

\extrarowheight=2pt \begin{tabular}%
{|>{\PBS\raggedright\hspace{0pt}}m{33mm}|
 >{\PBS\raggedright\hspace{0pt}}m{40mm}|
 >{\PBS\raggedright\hspace{0pt}}m{13mm}|
 >{\PBS\raggedright\hspace{0pt}}m{17mm}|
 >{\PBS\raggedright\hspace{0pt}}m{13mm}|
 >{\PBS\centering\hspace{0pt}}m{12mm}|
 >{\PBS\raggedright\hspace{0pt}}m{25mm}|}
\hline
 Diagrams\centering & $\fg$\centering & $v$\centering & $ev$ &
 $od$ &$png$ & $ng\le \min(* \ , *)$\\ \hline
 &&&
 $k_\ev-2$ & $k_\od$ &$\ev$ &$2k_\ev-4, 2k_\od$ \\
 &&&
 $k_\od$ &$k_\ev-2$ & $\od$ &$2k_\ev-3, 2k_\od-1$ \\
 &&&
 $k_\od-2$ & $k_\ev$ &$\ev$ &$2k_\ev, 2k_\od-4$ \\
 &&&
 $k_\ev$ &$k_\od-2$ & $\od$ &$2k_\ev-1,2k_\od-3$ \\ \cline{4-7}
 &&&
 $k_\ev-1$ & $k_\od-1$ & &$2k_\ev-2, 2k_\od-1$ \\
 \raisebox{3.5em}[0pt]{$\left.\arraycolsep=0pt\begin{array}{l}
 1)\;\includegraphics{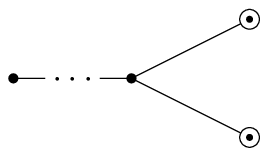}\\
 2)\;\includegraphics{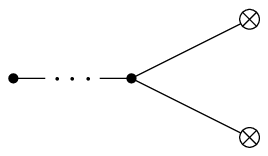}
 \end{array}\;\right\}$} &
 \raisebox{3.5em}[0pt]{$\begin{array}{l}\foo\fc(2;2k_\ev|2k_\od)\subplus\Kee I_0\\
 \text{if $k_\ev+k_\od$ is odd;}\\ \foo\fc(1;2k_\ev|2k_\od)\subplus\Kee I_0\\
 \text{if $k_\ev+k_\od$ is even.}\end{array}$}
 &
 \raisebox{3.5em}[0pt]{$k_{\ev}+k_{\od}$} &
 $k_\od-1$ &$k_\ev-1$ & &$2k_\ev-1, 2k_\od-2$ \\
\hline
 &&&
 $k_\ev-1$ & $k_\od$ &$\ev$ &$2k_\ev-2, 2k_\od$ \\
 &&&
 $k_\od$ &$k_\ev-1$ & $\od$ &$2k_\ev-1,2k_\od-1$ \\ \cline{4-7}
 &&&
 $k_\od-1$ & $k_\ev$ & $\ev$ &$2k_\ev, 2k_\od-2$ \\
 \raisebox{2.2em}[0pt]{$\left.\arraycolsep=0pt\begin{array}{l}
 3)\; \includegraphics{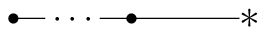}\\
 4)\;\includegraphics{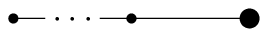}
 \end{array}\;\right\}$} &
 \raisebox{2.2em}[0pt]{$\foo^{(1)}_{I\Pi}(2k_\ev+1|2k_\od)$} &
 \raisebox{2.2em}[0pt]{$k_{\ev}+k_{\od}$} &
 $k_\ev$ &$k_\od-1$ & $\od$ &$2k_\ev-1, 2k_\od-1$ \\
 \hline
 5)\;\includegraphics{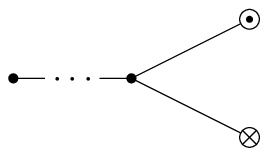}
 & $\begin{array}{l}\fpe\fc(2;m)\subplus\Kee I_0\\ \text{for $m$ odd;}\\
 \fpe\fc(1;m)\subplus\Kee I_0\\ \text{for $m$ even.}\end{array}$ & $m$&&&&\\ \hline
\end{tabular}

\end{center}

\ssec{Notation} The Dynkin diagrams in Table \ref{tbll} correspond
to Cartan matrix Lie superalgebras close to ortho-orthogonal and
periplectic Lie superalgebras. Each thin black dot may be $\motimes$
or $\odot$; the last five columns show conditions on the diagrams;
in the last four columns, it suffices to satisfy conditions in any
one row. Horizontal lines in the last four columns separate the
cases corresponding to different Dynkin diagrams. The notation are:
 $v$ is the total number of nodes in the diagram;  $ng$ is the
 number of \lq\lq grey'' nodes $\motimes$'s among the
thin black dots;  $png$ is the parity of this number;
 $ev$ and $od$ are the number of thin black dots such that the number
of $\otimes$'s to the left from them is even and odd, respectively.

\end{landscape}

\ssec{Systems of simple roots of the $\fe$-type Lie
superalgebras}\label{Sssr-e}

\sssbegin{Remark} Observe that if $p=2$ and the Cartan matrix has no
parameters, the reflections do not change the shape of the Dynkin
diagram. Therefore, for the $\fe$-superalgebras, it suffices to list
distributions of parities of the nodes in order to describe the
Dynkin diagrams. Since there are tens and even hundreds of diagrams
in these cases, this possibility saves a lot of space, see the lists
of all inequivalent Cartan matrices of the $\fe$-type Lie
superalgebras.
\end{Remark}

\sssec{$\fe(6, 1)\simeq \fe(6,5)$ of $\sdim 46|32$} All inequivalent
Cartan matrices are as follows (none of the matrices corresponding
to the symmetric pairs of Dynkin diagrams is excluded but are placed
one under the other for clarity, followed by three symmetric
diagrams):\tiny
\begin{equation*}\label{e65}
\begin{array}{llllllllllll}
\boxed{1)}&000010&3)&010001&5)&100110&7)&000011&9)&000110&11)&000111\\
\boxed{2)}&100000&4)&000101&6)&110010&8)&100001&10)&110000&12)&110001\\
13)&111001&15)&101001&17)&011000&19)&101100&21)&011001&23)&011110\\
14)&001111&16)&001011&18)&001100&20)&011010&22)&001101&24)&111100\\
25)&010100&26)&100010&27)&110110&&&&&&\\
 \end{array}
\end{equation*}

\normalsize

\sssec{$\fe(6, 6)$ of $\sdim = 38|40$} All inequivalent Cartan
matrices are as follows: \tiny
\begin{equation*}\label{e66}
\begin{array}{llllllllllllllll}
\boxed{1)}&000001&\boxed{2)}&000100&\boxed{3)}&001000
&\boxed{4)}&010000&5)&011011&6)&101110&7)&111110\\
8)&011100&9)&101111&10)&011101&11)&101010&12)&111101&
13)&010110&14)&101011\\
15)&110011&16)&001001&17)&011111&18)&110100&19)&010011&20)&101000&
21)&111011\\
22)&001010&23)&100011&24)&110101&25)&001110&26)&111000&27)&010010&
28)&100111\\
29)&100100&30)&110111&31)&100101&32)&111010&33)&010101&34)&010111&35)&101101\\
36)&111111\\
\end{array}
\end{equation*}

\normalsize

\sssec{$\fe(7,1)$ of $\sdim=80/78|54$} All inequivalent Cartan
matrices are as follows:\tiny
\begin{equation*}\label{e71}
\begin{array}{llllllllllllll}
\boxed{1)}&1000000&2)&1000010&3)&1000110&
4)&1001100&25)&0110000&26)&0110010&27)&0110110\\
5)&1010001&6)&1011001& 7)&1100000&8)&1100010&21)&0011010&
22)&0011110&23)&0100001\\
9)&1100110&
10)&1101100&11)&1110001&12)&1111001&17)&0001101&18)&0001111&
19)&0010100\\
13)&0000011&14)&0000101&15)&0000111&
16)&0001011&28)&0111100&24)&0101001&20)&0011000\\
\end{array}
\end{equation*}

\normalsize

\sssec{$\fe(7,6)$ of $\sdim=70/68|64$} All inequivalent Cartan
matrices are as follows:\tiny
\begin{equation*}\label{e76}
\begin{array}{llllllllllll}
\boxed{1)}&0000010&\boxed{2)}&0000100&3)&0000110&
\boxed{4)}&0001000&62)&1111100&63)&1111110\\
5)&0001010&6)&0001100&
7)&0001110&8)&0010001&60&1111000& 61)&1111010\\
9)&0010011&
10)&0010101&11)&0010111&12)&0011001& 58)&1110100&59)&1110110\\
13)&0011011&14)&0011101&15)&0011111&
\boxed{16)}&0100000&56)&1110000& 57)&1110010\\
17)&0100010&18)&0100100&
19)&0100110&20)&0101000&54)&1101101& 55)&1101111\\
21)&0101010&
22)&0101100&23)&0101110&24)&0110001& 52)&1101001& 53)&1101011\\
25)&0110011&26)&0110101&27)&0110111&
28)&0111001&50)&1100101&51)&1100111\\
29)&0111011&30)&0111101&
31)&0111111&32)&1000001&48)&1100001& 49)&1100011\\
33)&1000011&
34)&1000101&35)&1000111&36)&1001001& 46)&1011100&47)&1011110 \\
37)&1001011&38)&1001101&39)&1001111&
40)&1010000&44)&1011000& 45)&1011010\\
41)&1010010&42)&1010100& 43)&1010110&&
\end{array}
\end{equation*}

\normalsize

\sssec{$\fe(7,7)$ of $\sdim=64/62|70$} All inequivalent Cartan
matrices are as follows:\tiny
\begin{equation*}\label{e77}
\begin{array}{llllllllllll}
\boxed{1)}&0000001&2)& 0001001&\boxed{3)}&0010000&4)&0010010&
34)&1111011&35)&1111101\\
5)&0010110&6)&0011100&
7)&0100011&8)&0100101&32)&1110101& 33)&1110111\\
9)&0100111&
10)&0101011&11)&0101101&12)&0101111&30)&1101110& 31)&1110011\\
13)&0110100&14)&0111000&15)&0111010&
16)&0111110& 28)&1101000& 29)&1101010\\
17)&1000100&18)&1001000&
19)&1001010&20)&1001110&26)&1011111&27)&1100100\\
21)&1010011& 22)&1010101&23)&1010111&24)&1011011& 25)&1011101
\end{array}
\end{equation*}

\normalsize

 \sssec{$\fe(8,1)$ of $\sdim=136|112$} All inequivalent
Cartan matrices are as follows:\tiny
\begin{equation*}\label{81}
\begin{array}{llllllllll}
\boxed{1)}&10000000&2)&10000010&3)&10000011&4)&10000101&120)&01111110\\
5)&10000110&6)&10000111&7)&10001011&8)&10001100&119)&01111010\\
9)&10001101&10)&10001111&11)&10010001&12)&10010100&118)&01111001\\
13)&10011000&14)&10011001&15)&10011010&16)&10011110&117)&01111000\\
17)&10100000&18)&10100001&19)&10100010&20)&10100110&116)&01110100\\
21)&10101001&22)& 10101100&23)&10110000&24)&10110001&115)&01110001\\
25)&10110010&26)& 10110110&27&10111001&28)&10111100&114)&01101111\\
29)&11000000&30)&11000010&31)&11000011&32)&11000101&113)&01101101\\
33)&11000110&34)&11000111&35)&11001011&36)&11001100&112)&01101100\\
37)&11001101&38)&11001111&39)&11010001&40)&11010100&111)&01101011\\
41)&11011000&42)&11011001&43)&11011010&44)&11011110&110)&01100111\\
45)&11100000&46)&11100001&47)&11100010&48)&11100110&109)&01100110\\
49)&11101001&50)&11101100&51)&11110000&52)&11110001&108)&01100101\\
53)&11110010&54)& 11110110&55)&11111001&56)&11111100&107)&01100011\\
57)&00000011&\boxed{58)}&00000100&59)&00000101&60)&00000111&106)&01100010\\
\boxed{61)}&00001000&62)& 00001010&63)&00001011&64)&00001101&105)&01100000\\
65)&00001110&66)&00001111&67)&00010011&68)&00010100&104)&01011100\\
69)&00010101&70)& 00010111&71)&00011000&72)&00011010&103)&01011001\\
73)&00011011&74)&00011101&75)&00011110&76)&00011111&102)&01010110\\
77)&00100001&78)&00100100&79)&00101000&80)&00101001&101)&01010010\\
81)&00101010&82)&00101110&83)&00110000&84)&00110010&100)&01010001\\
85)&00110011&86)&00110101&87)&00110110&88)&00110111&99)&01010000\\
89)&00111011&90)&00111100&91)&00111101&92)&00111111&98)&01001100\\
\boxed{93)}&01000000&94)&01000001&95)&01000010&96)&01000110&
97)&01001001\\
\end{array}\end{equation*}

\normalsize

\sssec{$\fe(8,8)$ of $\sdim=120|128$} All inequivalent Cartan
matrices are as follows: \tiny
\begin{equation*}\label{e88}
\begin{array}{llllllllll}
\boxed{1)})&00000001&\boxed{2)}&00000010&\boxed{12)}&00100000&\boxed{6)}&00010000
&109)&11010101\\
5)&00001100&4)&00001001&7)&00010001&8)&00010010&110)&11010110\\
9)&00010110&10)&00011001&11)&00011100&3)&00000110&111)&11010111\\
13)&00100010&14)&00100011&15)&00100101&16)&00100110&112)&11011011\\
17)&00100111&18)&00101011&19)&00101100&20)&00101101&113)&11011100\\
21)&00101111&22)& 00110001&23)&00110100&24)&00111000&114)&11011101\\
25)&00111001&26)&00111010&27)&00111110&28)&01000011&115)&11011111\\
29)&01000100&30)& 01000101&31)&01000111&32)&01001000&116)&11100011\\
33)&01001010&34)&01001011&35)&01001101&36)&01001110&117)&11100100\\
37)&01001111&38)&01010011&39)&01010100&40)&01010101&118)&11100101\\
41)&01010111&42)&01011000&43)&01011010&44)&01011011&119)&11100111\\
45)&01011101&46)&01011110&47)&01011111&48)&01100001&120)&11101000\\
49)&01100100&50)&01101000&51)&01101001&52)&01101010&121)&11101010\\
53)&01101110&54)&01110000&55)&01110010&56)&01110011&122)&11101011\\
57)&01110101&58)&01110110&59)&01110111&60)&01111011&123)&11101101\\
61)&01111100&62)&01111101&63)&01111111&64)&10000001&124)&11101110\\
65)&10000100&66)&10001000&67)&10001001&68)&10001010&125)&11101111\\
69)&10001110&70)& 10010000&71)&10010010&72)&10010011&126)&11110011\\
73)&10010101&74)&10010110&75)&10010111&76)&10011011&127)&11110100\\
77)&10011100&78)&10011101&79)&10011111&80)&10100011&128)&11110101\\
81)&10100100&82)&10100101&83)&10100111&84)&10101000&129)&11110111\\
85)&10101010&86)&10101011&87)&10101101&88)&10101110&130)& 11111000\\
89)&10101111&90)&10110011&91)&10110100&92)&10110101&131)&11111010\\
93)&10110111&94)& 10111000&95)&10111010&96)&10111011&132)&11111011\\
97)&10111101&98)&10111110&99)&10111111&100)&11000001&133)&11111101\\
101)&11000100&102)&11001000&103)&11001001&104)&11001010&134)&11111110\\
105)&11001110&106)&11010000&107)&11010010&108)&11010011&\underline{135)}&
\underline{11111111}\\
 \end{array}
\end{equation*}

\normalsize

\sssec{$p=3$, Lie algebras} {\it Brown algebras} (CM=Cartan
matrix):\index{$\fbr(2, a)$, where $a\neq 0, -1$, Brown
algebra}\index{$\fbr(2)$, Brown algebra}\index{Brown
algebra}\index{$\fbr(3)$, Brown algebra}
\begin{equation}\label{br2a}
\begin{array}{l}
\fbr(2, a)\text{~with CM~}
\begin{pmatrix}2&-1\\a&2\end{pmatrix};\quad
\fbr(2)\text{~with CM~}
\begin{pmatrix}2&-1\\-1&0\end{pmatrix}\end{array}
\end{equation}
The reflections change the value of the parameter, so
\begin{equation}\label{brso}
\fbr(2,a)\simeq \fbr(2,a')\Longleftrightarrow a'=-(1+a).
\end{equation}

\begin{equation}\label{br3a}
\begin{array}{l}
1\fbr(3)\text{~with CM~}
\begin{pmatrix}2&-1&0\\-1&2&-1\\
0&-1&\ev\end{pmatrix};\quad 2\fbr(3)\text{~with CM~}
\begin{pmatrix}2&-1&0\\-2&2&-1\\
0&-1&\ev\end{pmatrix}\end{array}
\end{equation}

\sssec{$p=3$, Lie superalgebras}
{\it Brown superalgebra} $\fbrj(2;3)$\index{$\fbrj(2;3)$, Brown
superalgebra}\index{Brown superalgebra} of $\sdim = 10|8$ recently
discovered in \cite{El1} (Theorem 3.2(i); its Cartan matrices (first
listed in \cite{BGL1}) are as follows:
\[
1)\ \begin{pmatrix}0&-1\\
-2&1\end{pmatrix}, \quad 2)\ \begin{pmatrix}0&-1\\
-1&\ev\end{pmatrix}, \quad 3)\ \begin{pmatrix}1&-1\\
-1&\ev\end{pmatrix}.
\]
The Lie superalgebra $\fbrj(2;3)$ is a super analog of the Brown
algebra $\fbr(2)=\fbrj(2;3)_\ev$, its even part;
$\fbrj(2;3)_\od=R(2\pi_1)$ is irreducible $\fbrj(2;3)_\ev$-module.

Elduque \cite{El1, El2, CE, CE2} considered a particular case of the
classification problem of simple Lie superalgebras with Cartan
matrix and arranged the Lie (super)algebras he discovered in a
Supermagic Square all its entries being of the form $\fg(A)$. These
{\it Elduque and Cunha superalgebras}\index{Elduque and Cunha
superalgebra} are, indeed, exceptional ones.

\ssec{Elduque and Cunha superalgebras: Systems of simple
roots}\label{Secssr}

For de\-tails of description of Elduque and Cunha superalgebras in
terms of symmetric composition algebras, see \cite{El1,CE,CE2}. Here
we consider the simple Elduque and Cunha superalgebras with Cartan
matrix for $p=3$. In what follows, we list them using somewhat
shorter notation as compared with the original ones: Hereafter
$\fg(A,B)$ denotes the superalgebra occupying $(A, B)$th slot in the
Elduque Supermagic Square; the first Cartan matrix is usually the
one given in \cite{CE}, where only one Cartan matrix is given; the
other matrices are obtained from the first one by means of
reflections. Accordingly, $i\fg(A, B)$ is the shorthand for the
realization of $\fg(A, B)$ by means of the $i$th Cartan matrix.

There are no instances of isotropic even reflections.

\sssec{$\fg{}(1,6)$ of $\sdim =21|14$} We have $\fg(1,
6)_\ev=\fsp(6)$ and $\fg(1,6)_\od=R(\pi_3)$.
\begin{figure}[ht]\centering
\parbox{.5\linewidth}{\mbox{}\hfill$
\begin{pmatrix}
 -&-&2 \\
 -&-&1
 \end{pmatrix}
$\quad\mbox{}}\hfill
\parbox{.5\linewidth}{\includegraphics{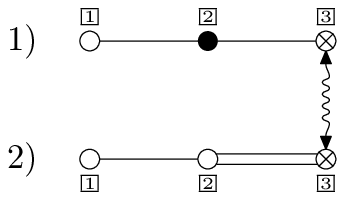}}
\end{figure}


{\footnotesize \[ 1) \begin{pmatrix}
 2&-1&0 \\
 -1&1&-1 \\
 0&-1&0
 \end{pmatrix}\quad
\boxed{2)} \begin{pmatrix}
 2&-1&0 \\
 -1&2&-2 \\
 0&-2&0
 \end{pmatrix}
\]
}

\clearpage

\sssec{$\fg(2,3)$ of $\sdim =12/10|14$} We have
$\fg(2,3)_\ev=\fgl(3)\oplus \fsl(2)$ and
$\fg(2,3)_\od=\fpsl(3)\otimes \id$.
\begin{figure}[ht]\centering
\parbox{.28\linewidth}{$
\begin{pmatrix}
 -&-&2 \\
 3&4&1 \\
 2&5&- \\
 5&2&- \\
 4&3&-
 \end{pmatrix}
$} \hfill
\parbox{.7\linewidth}{\includegraphics{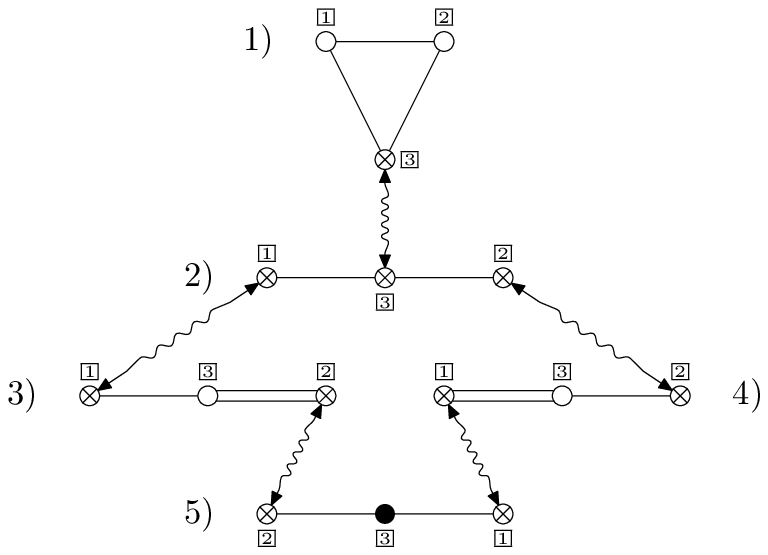}}
\end{figure}

{\footnotesize
\[
\begin{matrix}
\boxed{1)} \begin{pmatrix}
 2&-1&-1 \\
 -1&2&-1 \\
 -1&-1&0
 \end{pmatrix}\quad
\underline{2)} \begin{pmatrix}
 0&0&-1 \\
 0&0&-1 \\
 -1&-1&0
 \end{pmatrix}\quad
3) \begin{pmatrix}
 0&0&-1 \\
 0&0&-2 \\
 -1&-2&2
 \end{pmatrix}\\
\noalign{\vspace{2ex}} \hspace{1.em} 4) \begin{pmatrix}
 0&0&-2 \\
 0&0&-1 \\
 -2&-1&2
 \end{pmatrix}\quad
\underline{5)} \begin{pmatrix}
 0&0&-1 \\
 0&0&-1 \\
 -1&-1&1
 \end{pmatrix}\\
\end{matrix}
\]
}

\clearpage

\sssec{$\fg{}(3,6)$ of $\sdim =36|40$} We have $\fg(3,
6)_\ev=\fsp(8)$ and $\fg(3,6)_\od=R(\pi_3)$.
\[
\footnotesize{
\parbox{.28\linewidth}{$\begin{pmatrix}
 2&-&-&3 \\
 1&4&-&5 \\
 5&-&-&1 \\
 -&2&-&6 \\
 3&6&-&2 \\
 -&5&7&4 \\
 -&-&6&-
 \end{pmatrix} $}\quad
1) \begin{pmatrix}
 0&-1&0&0 \\
 -1&2&-1&0 \\
 0&-1&1&-1 \\
 0&0&-1&0
 \end{pmatrix}\quad
\underline{2)} \begin{pmatrix}
 0&-1&0&0 \\
 -1&0&-1&0 \\
 0&-1&1&-1 \\
 0&0&-1&0
 \end{pmatrix}
 }
\]


\begin{figure}[ht]\centering
\includegraphics{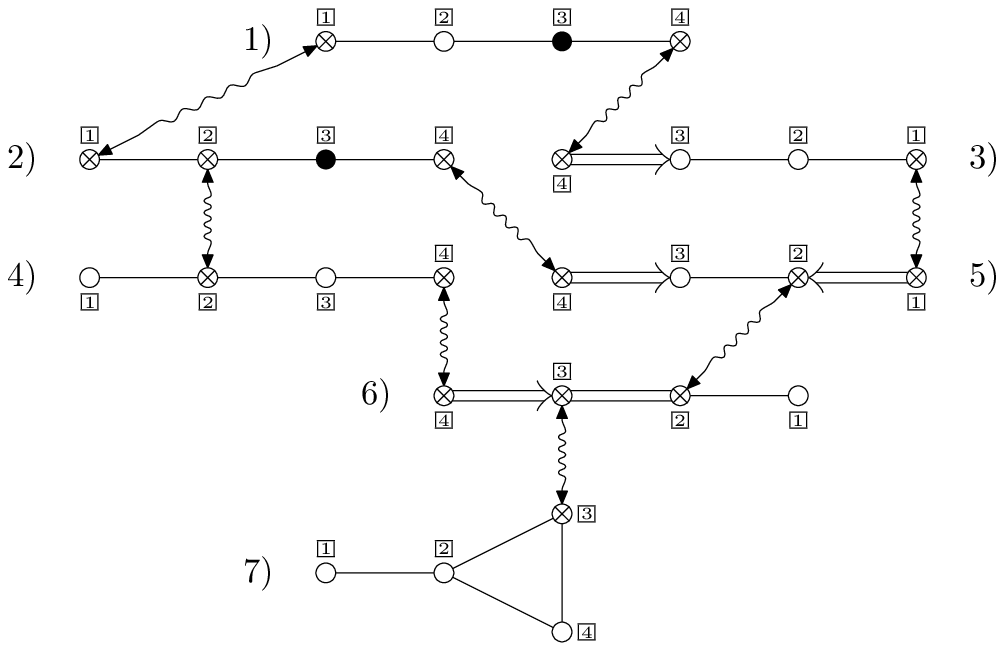}
\end{figure}


\footnotesize
\[
 \begin{matrix}
3) \begin{pmatrix}
 0&-1&0&0 \\
 -1&2&-1&0 \\
 0&-1&2&-2 \\
 0&0&-1&0
 \end{pmatrix} \quad
4)
\begin{pmatrix}
 2&-1&0&0 \\
 -1&0&-2&0 \\
 0&-2&2&-1 \\
 0&0&-1&0
 \end{pmatrix}\quad
5) \begin{pmatrix}
 0&-1&0&0 \\
 -2&0&-1&0 \\
 0&-1&2&-2 \\
 0&0&-1&0
 \end{pmatrix} \\
6) \begin{pmatrix}
 2&-1&0&0 \\
 -1&0&-2&0 \\
 0&-2&0&-2 \\
 0&0&-1&0
 \end{pmatrix}\quad
\boxed{7)} \begin{pmatrix}
 2&-1&0&0 \\
 -1&2&-1&-1 \\
 0&-1&0&-1 \\
 0&-1&-1&2
 \end{pmatrix}\\
\end{matrix}
\]
\normalsize

\clearpage

\sssec{$\fg{}(3,3)$ of $\sdim =23/21|16$} We have
\[\fg(3,3)_\ev=(\fo(7)\oplus \Kee z)\oplus \Kee d\text{~ and
~}\fg(3,3)_\od=(\spin_7)_+\oplus (\spin_7)_-;\] the action of $d$
separates the summands --- identical $\fo(7)$-modules $\spin_7$ ---
acting on one as the scalar multiplication by 1, on the other one by
$-1$. \[\tiny
\begin{pmatrix}
 -&-&-&2 \\
 -&-&3&1 \\
 -&4&2&- \\
 5&3&-&6 \\
 4&-&-&7 \\
 7&-&-&4 \\
 6&8&-&5 \\
 -&7&9&- \\
 10&-&8&- \\
 9&-&-&-
 \end{pmatrix} \quad \begin{matrix}\boxed{1)} \begin{pmatrix}
 2&-1&0&0 \\
 -1&2&-1&0 \\
 0&-2&2&-1 \\
 0&0&-1&0
 \end{pmatrix}\quad
2) \begin{pmatrix}
 2&-1&0&0 \\
 -1&2&-1&0 \\
 0&-1&0&-1 \\
 0&0&-1&0
 \end{pmatrix}\\
 3) \begin{pmatrix}
 2&-1&0&0 \\
 -1&0&-2&-2 \\
 0&-2&0&-2 \\
 0&-1&-1&2
 \end{pmatrix} \quad
4) \begin{pmatrix}
 0&-1&0&0 \\
 -2&0&-1&-1 \\
 0&-1&2&0 \\
 0&-1&0&0
 \end{pmatrix}
 \end{matrix}
\]

\begin{figure}[ht]\centering
\includegraphics{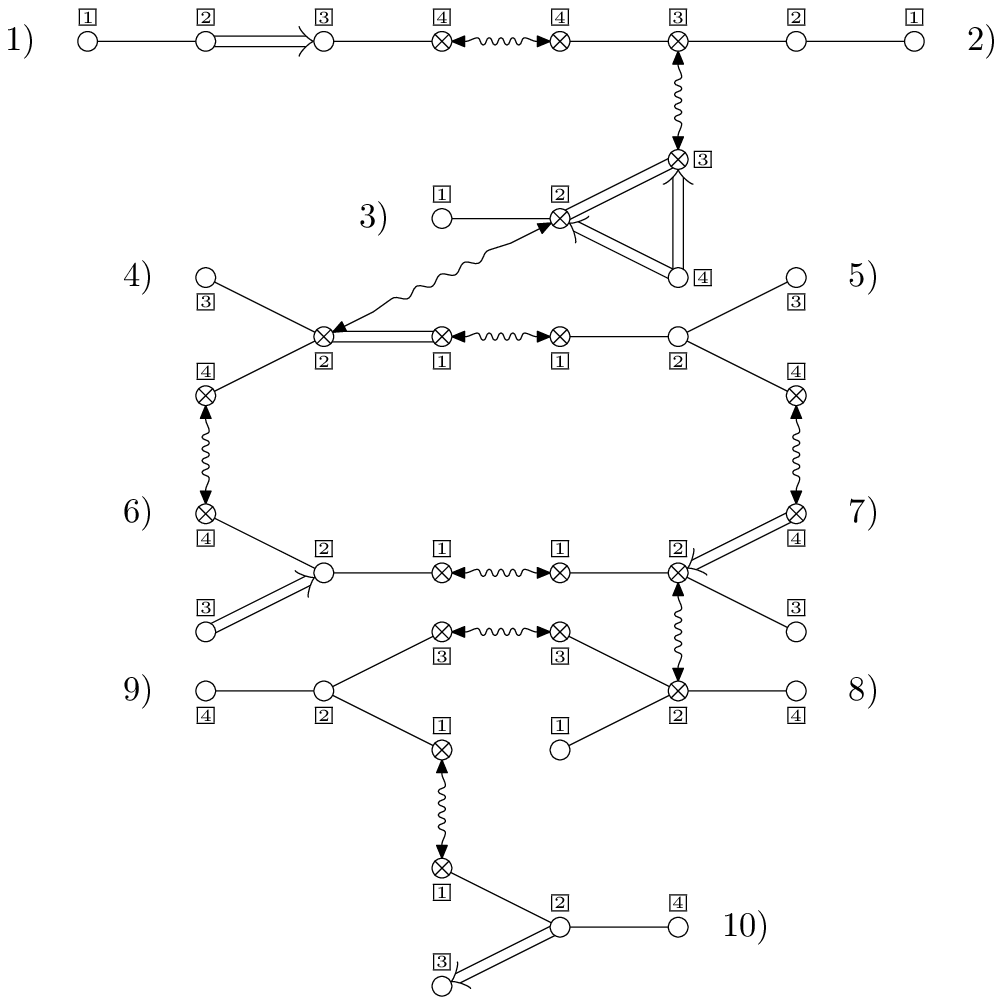}
\end{figure}

\tiny

\[
\begin{matrix}
\quad 5) \begin{pmatrix}
 0&-1&0&0 \\
 -1&2&-1&-1 \\
 0&-1&2&0 \\
 0&-1&0&0
 \end{pmatrix}\quad
6) \begin{pmatrix}
 0&-1&0&0 \\
 -1&2&-2&-1 \\
 0&-1&2&0 \\
 0&-1&0&0
 \end{pmatrix} \quad
7)\begin{pmatrix}
 0&-1&0&0 \\
 -1&0&-1&-2 \\
 0&-1&2&0 \\
 0&-1&0&0
 \end{pmatrix}\\
8) \begin{pmatrix}
 2&-1&-1&0 \\
 -2&0&-2&-1 \\
 -1&-1&0&0 \\
 0&-1&0&2
 \end{pmatrix}\quad
9)
\begin{pmatrix}
 0&0&-1&0 \\
 0&2&-1&-1 \\
 -1&-1&0&0 \\
 0&-1&0&2
 \end{pmatrix}\quad
\boxed{10)} \begin{pmatrix}
 0&0&-1&0 \\
 0&2&-1&-1 \\
 -1&-2&2&0 \\
 0&-1&0&2
 \end{pmatrix}\\
\end{matrix}
\]

\normalsize

\sssec{$\fg{}(4,3)$ of $\sdim =24|26$} We have
$\fg(4,3)_\ev=\fsp(6)\oplus \fsl(2)$ and
$\fg(4,3)_\od=R(\pi_2)\otimes \id$.
\[\tiny{ \;\begin{pmatrix}
 -&-&-&2 \\
 -&3&-&1 \\
 4&2&5&- \\
 3&-&6&- \\
 6&-&3&7 \\
 5&8&4&9 \\
 9&-&-&5 \\
 -&6&-&10 \\
 7&10&-&6 \\
 -&9&-&8
 \end{pmatrix}\quad \begin{matrix}\boxed{1)} \begin{pmatrix}
 2&-1&0&0 \\
 -1&2&-2&-1 \\
 0&-1&2&0 \\
 0&-1&0&0
 \end{pmatrix}\quad
2) \begin{pmatrix}
 2&-1&0&0 \\
 -1&0&-2&-2 \\
 0&-1&2&0 \\
 0&-1&0&0
 \end{pmatrix}\\
3) \begin{pmatrix}
 0&-1&0&0 \\
 -2&0&-1&-1 \\
 0&-1&0&-1 \\
 0&-1&-1&2
 \end{pmatrix} \quad
4) \begin{pmatrix}
 0&-1&0&0 \\
 -1&2&-1&-1 \\
 0&-1&0&-1 \\
 0&-1&-1&2
 \end{pmatrix}\end{matrix} }
\]


\begin{figure}[ht]\centering
\includegraphics{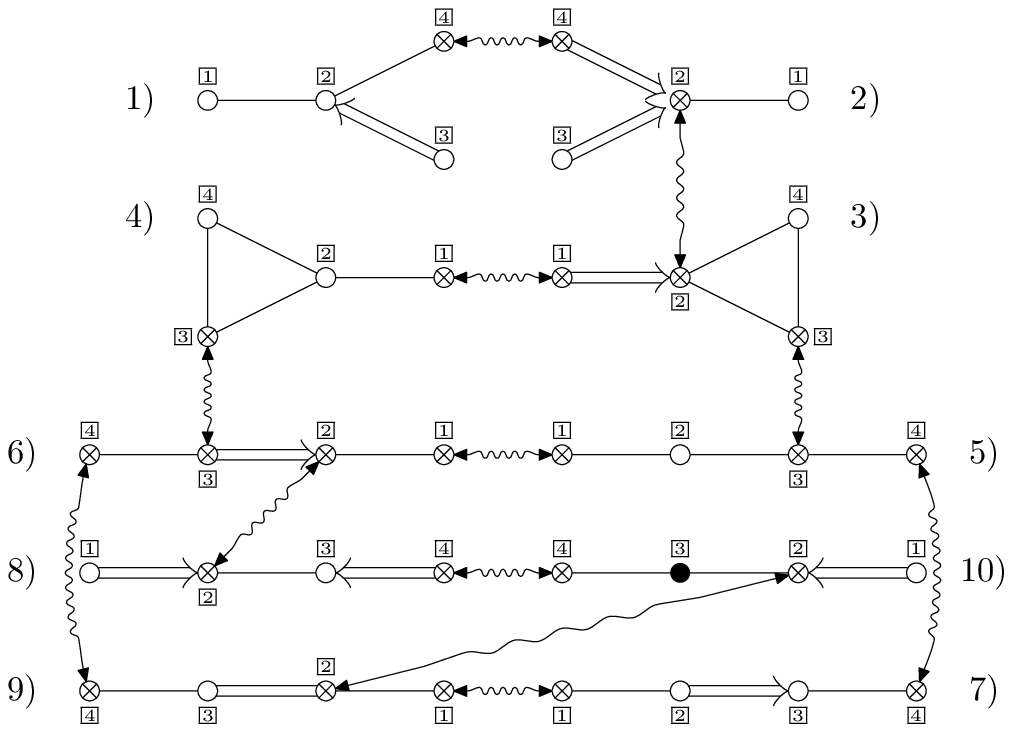}
\end{figure}


\tiny
\[
\begin{matrix}
5) \begin{pmatrix}
 0&-1&0&0 \\
 -1&2&-1&0 \\
 0&-1&0&-1 \\
 0&0&-1&0
 \end{pmatrix}\quad
\underline{6)} \begin{pmatrix}
 0&-1&0&0 \\
 -1&0&-2&0 \\
 0&-1&0&-1 \\
 0&0&-1&0
 \end{pmatrix} \quad
7) \begin{pmatrix}
 0&-1&0&0 \\
 -1&2&-1&0 \\
 0&-2&2&-1 \\
 0&0&-1&0
 \end{pmatrix}\\
8) \begin{pmatrix}
 2&-1&0&0 \\
 -2&0&-1&0 \\
 0&-1&2&-2 \\
 0&0&-1&0
 \end{pmatrix}\quad
9) \begin{pmatrix}
 0&-1&0&0 \\
 -1&0&-2&0 \\
 0&-2&2&-1 \\
 0&0&-1&0
 \end{pmatrix}\quad
10) \begin{pmatrix}
 2&-1&0&0 \\
 -2&0&-1&0 \\
 0&-1&1&-1 \\
 0&0&-1&0
 \end{pmatrix}\\
\end{matrix}
\]
\normalsize

\sssec{$\fg{}(2,6)$ of $\sdim =36/34|20$} We have
$\fg(2,6)_\ev=\fgl(6)$ and $\fg(2,6)_\od=R(\pi_3)$.

\begin{figure}[ht]\centering
\includegraphics{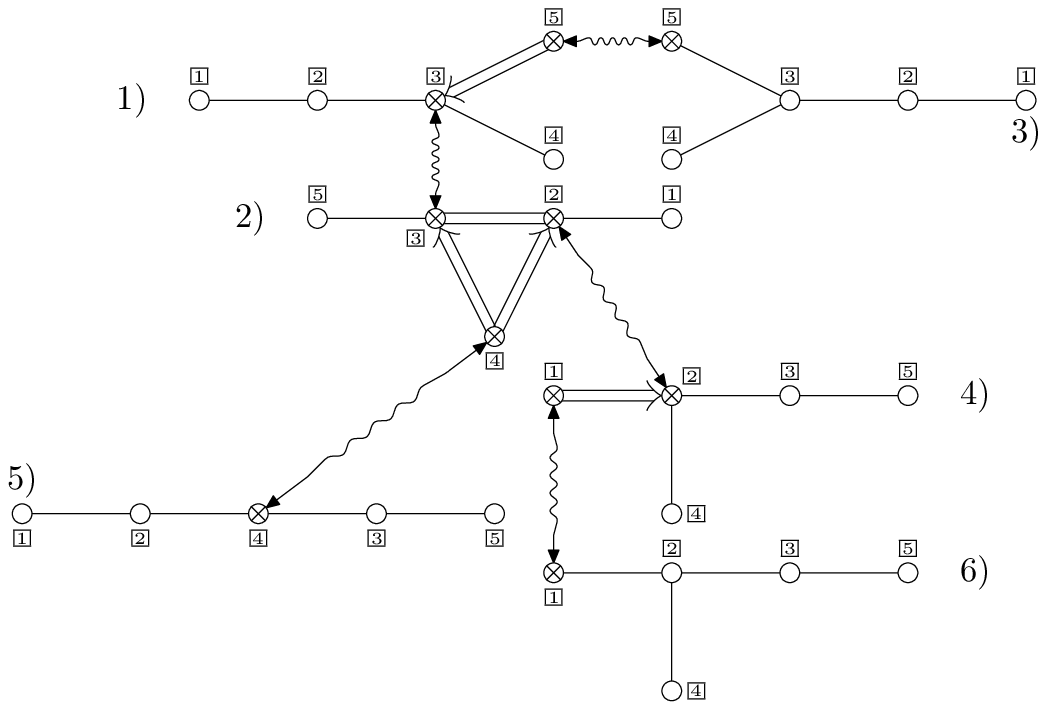}
\end{figure}

\[
\tiny{
\begin{pmatrix}
 -&-&2&-&3 \\
 -&4&1&5&- \\
 -&-&-&-&1 \\
 6&2&-&-&- \\
 -&-&-&2&- \\
 4&-&-&-&-
 \end{pmatrix} \quad
 1) \begin{pmatrix}
 2&-1&0&0&0 \\
 -1&2&-1&0&0 \\
 0&-1&0&-1&-2 \\
 0&0&-1&2&0 \\
 0&0&-1&0&0
 \end{pmatrix}\quad
2) \begin{pmatrix}
 2&-1&0&0&0 \\
 -1&0&-2&-2&0 \\
 0&-2&0&-2&-1 \\
 0&-1&-1&0&0 \\
 0&0&-1&0&2
 \end{pmatrix}
 }
\]

\tiny

\[\begin{matrix} \arraycolsep=3pt \boxed{3)} \begin{pmatrix}
 2&-1&0&0&0 \\
 -1&2&-1&0&0 \\
 0&-1&2&-1&-1 \\
 0&0&-1&2&0 \\
 0&0&-1&0&0
 \end{pmatrix} \sim
\boxed{6)} \begin{pmatrix}
 0&-1&0&0&0 \\
 -1&2&-1&-1&0 \\
 0&-1&2&0&-1 \\
 0&-1&0&2&0 \\
 0&0&-1&0&2
 \end{pmatrix}\\
4)
 \begin{pmatrix}
 0&-1&0&0&0 \\
 -2&0&-1&-1&0 \\
 0&-1&2&0&-1 \\
 0&-1&0&2&0 \\
 0&0&-1&0&2
 \end{pmatrix}\quad\boxed{5)}\begin{pmatrix}
 2&-1&0&0&0 \\
 -1&2&0&-1&0 \\
 0&0&2&-1&-1 \\
 0&-1&-1&0&0 \\
 0&0&-1&0&2
 \end{pmatrix}
\\
\end{matrix}
\]

\normalsize

\clearpage

\sssec{$\fg{}(8,3)$ of $\sdim =55|50$} We have
$\fg(8,3)_\ev=\ff(4)\oplus \fsl(2)$ and 
$\fg(8,3)_\od=R(\pi_4)\otimes \id$. \tiny
\[
\begin{pmatrix}
 -&-&-&-&2 \\
 -&-&-&3&1 \\
 -&-&4&2&- \\
 -&5&3&-&- \\
 6&4&-&7&- \\
 5&-&-&8&- \\
 8&-&-&5&9 \\
 7&10&-&6&11 \\
 11&-&-&-&7 \\
 -&8&12&-&13 \\
 9&13&-&-&8 \\
 14&-&10&-&15 \\
 -&11&15&16&10 \\
 12&-&-&-&17 \\
 17&-&13&18&12 \\
 -&-&18&13&- \\
 15&-&-&19&14 \\
 19&20&16&15&- \\
 18&21&-&17&- \\
 21&18&-&-&- \\
 20&19&-&-&-
 \end{pmatrix}
\begin{matrix}
\boxed{1)} 
\begin{pmatrix}
 2&-1&0&0&0 \\
 -1&2&-1&0&0 \\
 0&-2&2&-1&0 \\
 0&0&-1&2&-1 \\
 0&0&0&1&0
 \end{pmatrix} \\
2) \begin{pmatrix}
 2&-1&0&0&0 \\
 -1&2&-1&0&0 \\
 0&-2&2&-1&0 \\
 0&0&-2&0&-1 \\
 0&0&0&-1&0
 \end{pmatrix}\\
3) \begin{pmatrix}
 2&-1&0&0&0 \\
 -1&2&-1&0&0 \\
 0&-1&0&-1&0 \\
 0&0&-1&0&-2 \\
 0&0&0&-1&2
 \end{pmatrix} \\
4) \begin{pmatrix}
 2&-1&0&0&0 \\
 -1&0&-2&-2&0 \\
 0&-1&0&-1&0 \\
 0&-1&-1&2&-1 \\
 0&0&0&-1&2
 \end{pmatrix}
 \end{matrix}
\begin{matrix}
5) \begin{pmatrix}
 0&-1&0&0&0 \\
 -2&0&-1&-1&0 \\
 0&-1&2&0&0 \\
 0&-1&0&0&-2 \\
 0&0&0&-1&2
 \end{pmatrix}\\
6) \begin{pmatrix}
 0&-1&0&0&0 \\
 -1&2&-1&-1&0 \\
 0&-1&2&0&0 \\
 0&-1&0&0&-2 \\
 0&0&0&-1&2
 \end{pmatrix} \\
7) \begin{pmatrix}
 0&-1&0&0&0 \\
 -1&2&-2&-1&0 \\
 0&-1&2&0&0 \\
 0&-2&0&0&-1 \\
 0&0&0&-1&0
 \end{pmatrix}\\
8) \begin{pmatrix}
 0&-1&0&0&0 \\
 -1&0&-1&-2&0 \\
 0&-1&2&0&0 \\
 0&-2&0&0&-1 \\
 0&0&0&-1&0
 \end{pmatrix}
 \end{matrix}
\]

\[
\begin{matrix}
9) \begin{pmatrix}
 0&-1&0&0&0 \\
 -1&2&-2&-1&0 \\
 0&-1&2&0&0 \\
 0&-1&0&2&-1 \\
 0&0&0&-1&0
 \end{pmatrix} \quad
10) \arraycolsep=2.8pt
 \begin{pmatrix}
 2&-1&-1&0&0 \\
 1&0&1&2&0 \\
 1&1&0&0&0 \\
 0&-1&0&2&-1 \\
 0&0&0&-1&0
 \end{pmatrix}\quad
11) \begin{pmatrix}
 0&-1&0&0&0 \\
 -1&0&-1&-2&0 \\
 0&-1&2&0&0 \\
 0&-1&0&2&-1 \\
 0&0&0&-1&0
 \end{pmatrix}\\
 \end{matrix}
\]
\[
\begin{matrix}
12) \arraycolsep=2.8pt \begin{pmatrix}
 0&0&-1&0&0 \\
 0&2&-1&-1&0 \\
 -1&-1&0&0&0 \\
 0&-1&0&2&-1 \\
 0&0&0&-1&0
 \end{pmatrix} \quad 13)
 \begin{pmatrix}
 2&-1&-1&0&0 \\
 -1&0&-1&-2&0 \\
 -1&-1&0&0&0 \\
 0&-2&0&0&-1 \\
 0&0&0&-1&0
 \end{pmatrix}\quad
14) \begin{pmatrix}
 0&0&-1&0&0 \\
 0&2&-1&-1&0 \\
 -1&-2&2&0&0 \\
 0&-1&0&2&-1 \\
 0&0&0&-1&0
 \end{pmatrix} \end{matrix}
\]
\[
\begin{matrix}
15) \arraycolsep=2.8pt\begin{pmatrix}
 0&0&-1&0&0 \\
 0&2&-1&-1&0 \\
 -1&-1&0&0&0 \\
 0&-2&0&0&-1 \\
 0&0&0&-1&0
 \end{pmatrix} \quad
16) 
 \begin{pmatrix}
 2&-1&-1&0&0 \\
 -1&2&-1&-1&0 \\
 -1&-1&0&0&0 \\
 0&-1&0&0&-2 \\
 0&0&0&-1&2
 \end{pmatrix}\quad
17) \begin{pmatrix}
 0&0&-1&0&0 \\
 0&2&-1&-1&0 \\
 -1&-2&2&0&0 \\
 0&-2&0&0&-1 \\
 0&0&0&-1&0
 \end{pmatrix} \end{matrix}
\]
\[
\begin{matrix}
18) \begin{pmatrix}
 0&0&-1&0&0 \\
 0&0&-2&-1&0 \\
 -1&-1&0&0&0 \\
 0&-1&0&0&-2 \\
 0&0&0&-1&2
 \end{pmatrix} \quad
19) 
 \begin{pmatrix}
 0&0&-1&0&0 \\
 0&0&-2&-1&0 \\
 -1&-2&2&0&0 \\
 0&-1&0&0&-2 \\
 0&0&0&-1&2
 \end{pmatrix}\\
20) \begin{pmatrix}
 0&0&-1&0&0 \\
 0&0&-1&-2&0 \\
 -2&-1&2&0&0 \\
 0&-1&0&2&-1 \\
 0&0&0&-1&2
 \end{pmatrix}\quad
21) \begin{pmatrix}
 0&0&-1&0&0 \\
 0&0&-1&-2&0 \\
 -1&-1&1&0&0 \\
 0&-1&0&2&-1 \\
 0&0&0&-1&2
 \end{pmatrix}\\
\end{matrix}
\]

\begin{center}
\includegraphics{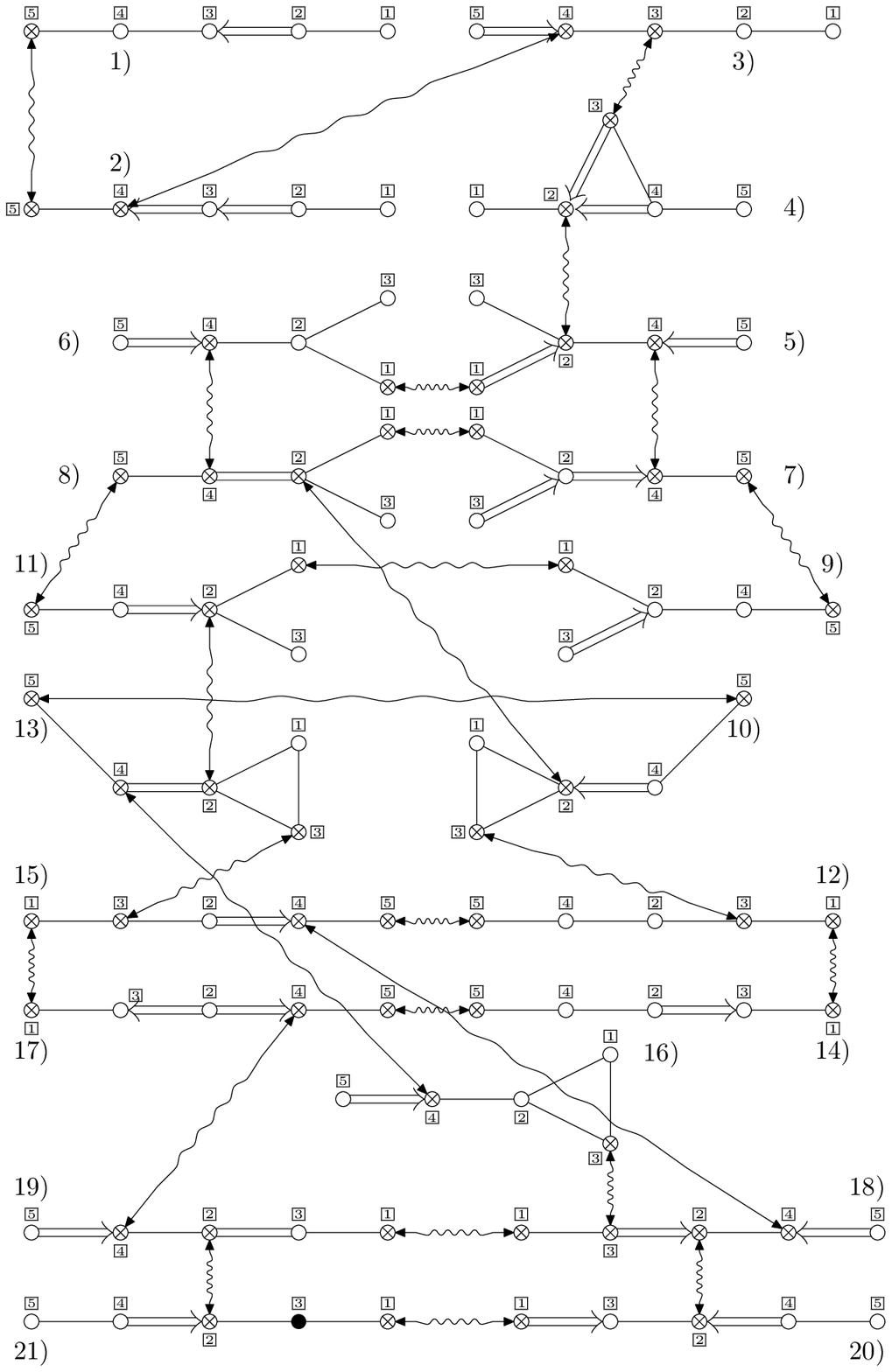}
\end{center}

\normalsize

\clearpage

\sssec{$\fg{}(4,6)$ of $\sdim =66|32$} We have
$\fg(4,6)_\ev=\fo(12)$ and $\fg(4,6)_\od=R(\pi_{5})$.
\[
\tiny \begin{pmatrix}
 -&-&-&2&-&3 \\[-1pt]
 -&-&4&1&5&- \\[-1pt]
 -&-&-&-&-&1 \\[-1pt]
 -&6&2&-&-&- \\[-1pt]
 -&-&-&-&2&- \\[-1pt]
 7&4&-&-&-&- \\[-1pt]
 6&-&-&-&-&-
 \end{pmatrix}\quad 1)
\begin{pmatrix}
 2&-1&0&0&0&0 \\
 -1&2&-1&0&0&0 \\
 0&-1&2&-1&0&0 \\
 0&0&-2&0&-2&-1 \\
 0&0&0&-1&2&0 \\
 0&0&0&-1&0&0
 \end{pmatrix}
\]

\begin{figure}[ht]\centering
\includegraphics{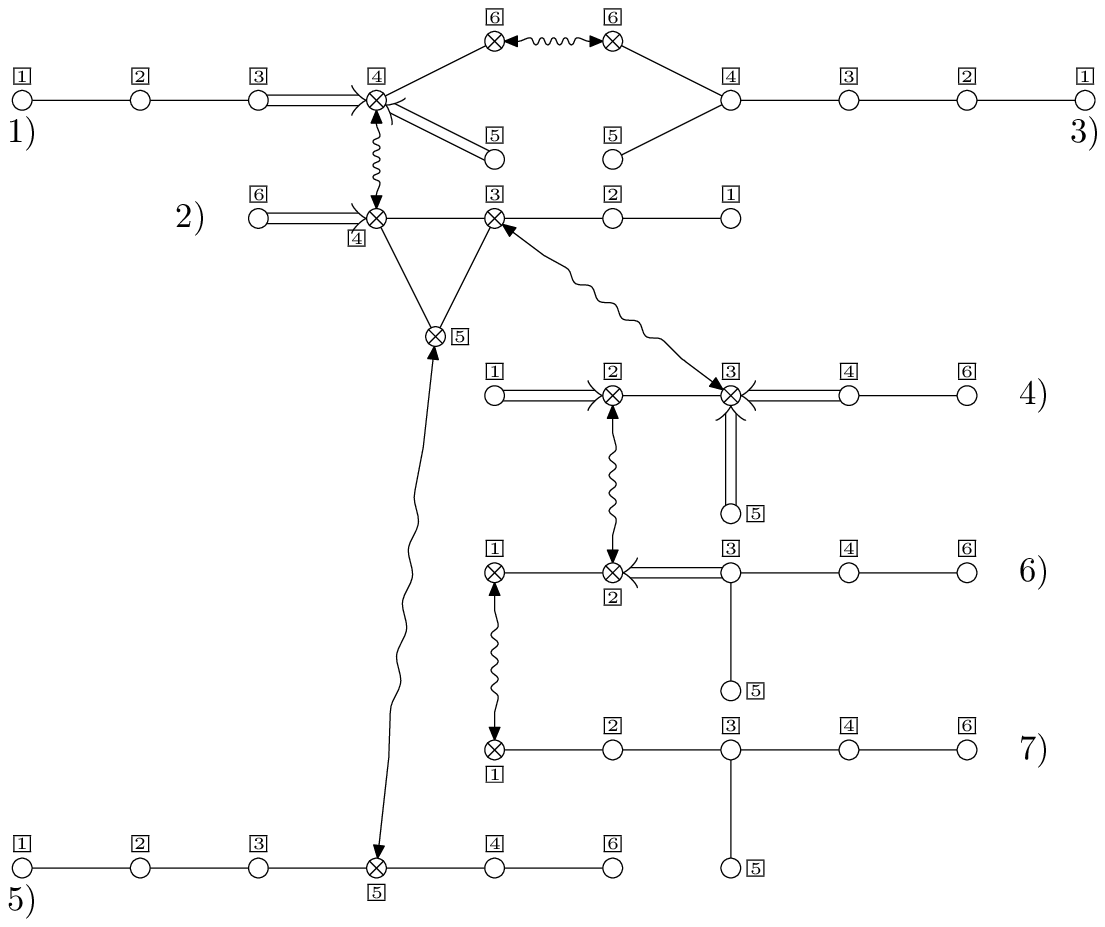}
\end{figure}

\vspace*{-5mm}

\tiny
\[
\begin{matrix}
2) \begin{pmatrix}
 2&-1&0&0&0&0 \\
 -1&2&-1&0&0&0 \\
 0&-2&0&-1&-1&0 \\
 0&0&-1&0&-1&-2 \\
 0&0&-1&-1&0&0 \\
 0&0&0&-1&0&2
 \end{pmatrix}
\quad \boxed{3)} \begin{pmatrix}
 2&-1&0&0&0&0 \\
 -1&2&-1&0&0&0 \\
 0&-1&2&-1&0&0 \\
 0&0&-1&2&-1&-1 \\
 0&0&0&-1&2&0 \\
 0&0&0&-1&0&0
 \end{pmatrix} \\
\noalign{\vspace{1ex}} 
4) \arraycolsep=3pt
 \begin{pmatrix}
 2&-1&0&0&0&0 \\
 -2&0&-1&0&0&0 \\
 0&-1&0&-2&-2&0 \\
 0&0&-1&2&0&-1 \\
 0&0&-1&0&2&0 \\
 0&0&0&-1&0&2
 \end{pmatrix}\quad
\boxed{5)} \begin{pmatrix}
 2&-1&0&0&0&0 \\
 -1&2&-1&0&0&0 \\
 0&-1&2&0&-1&0 \\
 0&0&0&2&-1&-1 \\
 0&0&-1&-1&0&0 \\
 0&0&0&-1&0&2
 \end{pmatrix} \\
\noalign{\vspace{1ex}} 
6) \arraycolsep=3pt
 \begin{pmatrix}
 0&-1&0&0&0&0 \\
 -1&0&-2&0&0&0 \\
 0&-1&2&-1&-1&0 \\
 0&0&-1&2&0&-1 \\
 0&0&-1&0&2&0 \\
 0&0&0&-1&0&2
 \end{pmatrix}\quad
\boxed{7)} \begin{pmatrix}
 0&-1&0&0&0&0 \\
 -1&2&-1&0&0&0 \\
 0&-1&2&-1&-1&0 \\
 0&0&-1&2&0&-1 \\
 0&0&-1&0&2&0 \\
 0&0&0&-1&0&2
 \end{pmatrix}\\
\end{matrix}
\]


\normalsize

\sssec{$\fg{}(6,6)$ of $\sdim =78|64$} We have
$\fg(6,6)_\ev=\fo(13)$ and $\fg(6,6)_\od=\spin_{13}$.
\[
\tiny \begin{pmatrix}
 2&3&-&4&-&5 \\[-1.5pt]
 1&-&-&6&-&7 \\[-1.5pt]
 -&1&8&9&-&10 \\[-1.5pt]
 6&9&11&1&12&- \\[-1.5pt]
 7&10&-&-&-&1 \\[-1.5pt]
 4&-&13&2&14&- \\[-1.5pt]
 5&-&-&-&-&2 \\[-1.5pt]
 -&-&3&-&-&15 \\[-1.5pt]
 -&4&-&3&16&- \\[-1.5pt]
 -&5&15&-&-&3 \\[-1.5pt]
 13&-&4&-&-&- \\[-1.5pt]
 14&16&-&-&4&- \\[-1.5pt]
 11&17&6&-&-&- \\[-1.5pt]
 12&-&-&-&6&- \\[-1.5pt]
 -&-&10&18&-&8 \\[-1.5pt]
 -&12&19&-&9&- \\[-1.5pt]
 -&13&-&-&-&- \\[-1.5pt]
 -&-&-&15&20&- \\[-1.5pt]
 -&-&16&-&-&- \\[-1.5pt]
 -&-&-&-&18&21 \\[-1.5pt]
 -&-&-&-&-&20
 \end{pmatrix} \begin{matrix}
1) 
\begin{pmatrix}
 0&-1&0&0&0&0 \\
 -1&0&-2&0&0&0 \\
 0&-1&2&-1&0&0 \\
 0&0&-2&0&-2&-1 \\
 0&0&0&-1&2&0 \\
 0&0&0&-1&0&0
 \end{pmatrix}\\
2) \begin{pmatrix}
 0&-2&0&0&0&0 \\
 -1&2&-1&0&0&0 \\
 0&-1&2&-1&0&0 \\
 0&0&-2&0&-2&-1 \\
 0&0&0&-1&2&0 \\
 0&0&0&-1&0&0
 \end{pmatrix}\\
3) \begin{pmatrix}
 2&-1&0&0&0&0 \\
 -2&0&-1&0&0&0 \\
 0&-1&0&-2&0&0 \\
 0&0&-2&0&-2&-1 \\
 0&0&0&-1&2&0 \\
 0&0&0&-1&0&0
 \end{pmatrix}
 \end{matrix}
\]

\tiny
\[
\begin{matrix}
4) \arraycolsep=1pt
\begin{pmatrix}
 0&-1&0&0&0&0 \\
 -1&0&-2&0&0&0 \\
 0&-2&0&-1&-1&0 \\
 0&0&-1&0&-1&-2 \\
 0&0&-1&-1&0&0 \\
 0&0&0&-1&0&2
 \end{pmatrix}\quad
5) \begin{pmatrix}
 0&-1&0&0&0&0 \\
 -1&0&-2&0&0&0 \\
 0&-1&2&-1&0&0 \\
 0&0&-1&2&-1&-1 \\
 0&0&0&-1&2&0 \\
 0&0&0&-1&0&0
 \end{pmatrix}\quad
6) \begin{pmatrix}
 0&-1&0&0&0&0 \\
 -1&2&-1&0&0&0 \\
 0&-2&0&-1&-1&0 \\
 0&0&-1&0&-1&-2 \\
 0&0&-1&-1&0&0 \\
 0&0&0&-1&0&2
 \end{pmatrix} \\
\noalign{\vspace{3ex}} 
7) \arraycolsep=1pt
 \begin{pmatrix}
 0&-1&0&0&0&0 \\
 -1&2&-1&0&0&0 \\
 0&-1&2&-1&0&0 \\
 0&0&-1&2&-1&-1 \\
 0&0&0&-1&2&0 \\
 0&0&0&-1&0&0
 \end{pmatrix}\quad
8) \begin{pmatrix}
 2&-1&0&0&0&0 \\
 -1&2&-1&0&0&0 \\
 0&-2&0&-1&0&0 \\
 0&0&-1&2&-2&-1 \\
 0&0&0&-1&2&0 \\
 0&0&0&-1&0&0
 \end{pmatrix} \quad
9) \arraycolsep=1pt
\begin{pmatrix}
 2&-1&0&0&0&0 \\
 -2&0&-1&0&0&0 \\
 0&-1&2&-1&-1&0 \\
 0&0&-1&0&-1&-2 \\
 0&0&-1&-1&0&0 \\
 0&0&0&-1&0&2
 \end{pmatrix}\\
\noalign{\vspace{3ex}} 
10) \arraycolsep=1pt
 \begin{pmatrix}
 2&-1&0&0&0&0 \\
 -2&0&-1&0&0&0 \\
 0&-1&0&-2&0&0 \\
 0&0&-1&2&-1&-1 \\
 0&0&0&-1&2&0 \\
 0&0&0&-1&0&0
 \end{pmatrix}\quad
11) \begin{pmatrix}
 0&-1&0&0&0&0 \\
 -1&2&-1&0&0&0 \\
 0&-1&0&-2&-2&0 \\
 0&0&-1&2&0&-1 \\
 0&0&-1&0&2&0 \\
 0&0&0&-1&0&2
 \end{pmatrix} \quad
12) \arraycolsep=1pt
 \begin{pmatrix}
 0&-1&0&0&0&0 \\
 -1&0&-2&0&0&0 \\
 0&-1&2&0&-1&0 \\
 0&0&0&2&-1&-1 \\
 0&0&-1&-1&0&0 \\
 0&0&0&-1&0&2
 \end{pmatrix}\\
\noalign{\vspace{3ex}} 
13) \arraycolsep=1pt
 \begin{pmatrix}
 0&-1&0&0&0&0 \\
 -2&0&-1&0&0&0 \\
 0&-1&0&-2&-2&0 \\
 0&0&-1&2&0&-1 \\
 0&0&-1&0&2&0 \\
 0&0&0&-1&0&2
 \end{pmatrix}\quad
14) \begin{pmatrix}
 0&-1&0&0&0&0 \\
 -1&2&-1&0&0&0 \\
 0&-1&2&0&-1&0 \\
 0&0&0&2&-1&-1 \\
 0&0&-1&-1&0&0 \\
 0&0&0&-1&0&2
 \end{pmatrix} \quad
15) \begin{pmatrix}
 2&-1&0&0&0&0 \\
 -1&2&-1&0&0&0 \\
 0&-2&0&-1&0&0 \\
 0&0&-1&0&-2&-2 \\
 0&0&0&-1&2&0 \\
 0&0&0&-1&0&0
 \end{pmatrix}\\
\noalign{\vspace{3ex}} 
16) \arraycolsep=1pt
 \begin{pmatrix}
 2&-1&0&0&0&0 \\
 -2&0&-1&0&0&0 \\
 0&-1&0&0&-2&0 \\
 0&0&0&2&-1&-1 \\
 0&0&-1&-1&0&0 \\
 0&0&0&-1&0&2
 \end{pmatrix}\quad
\boxed{17)} \begin{pmatrix}
 2&-1&0&0&0&0 \\
 -1&0&-2&0&0&0 \\
 0&-1&2&-1&-1&0 \\
 0&0&-1&2&0&-1 \\
 0&0&-1&0&2&0 \\
 0&0&0&-1&0&2
 \end{pmatrix} \quad
18) \begin{pmatrix}
 2&-1&0&0&0&0 \\
 -1&2&-1&0&0&0 \\
 0&-1&2&-1&0&0 \\
 0&0&-2&0&-1&-1 \\
 0&0&0&-1&0&-1 \\
 0&0&0&-1&-1&2
 \end{pmatrix}\\
\noalign{\vspace{3ex}} 
\boxed{19)} \arraycolsep=1pt
 \begin{pmatrix}
 2&-1&0&0&0&0 \\
 -1&2&-1&0&0&0 \\
 0&-2&0&0&-1&0 \\
 0&0&0&2&-1&-1 \\
 0&0&-1&-2&2&0 \\
 0&0&0&-1&0&2
 \end{pmatrix} \quad
20) \begin{pmatrix}
 2&-1&0&0&0&0 \\
 -1&2&-1&0&0&0 \\
 0&-1&2&-1&0&0 \\
 0&0&-1&2&-1&0 \\
 0&0&0&-1&0&-1 \\
 0&0&0&0&-1&0
 \end{pmatrix}\quad
\boxed{21)} \begin{pmatrix}
 2&-1&0&0&0&0 \\
 -1&2&-1&0&0&0 \\
 0&-1&2&-1&0&0 \\
 0&0&-1&2&-1&0 \\
 0&0&0&-2&2&-1 \\
 0&0&0&0&-1&0
 \end{pmatrix}\\
\end{matrix}
\]

\begin{figure}[ht]\centering
\includegraphics{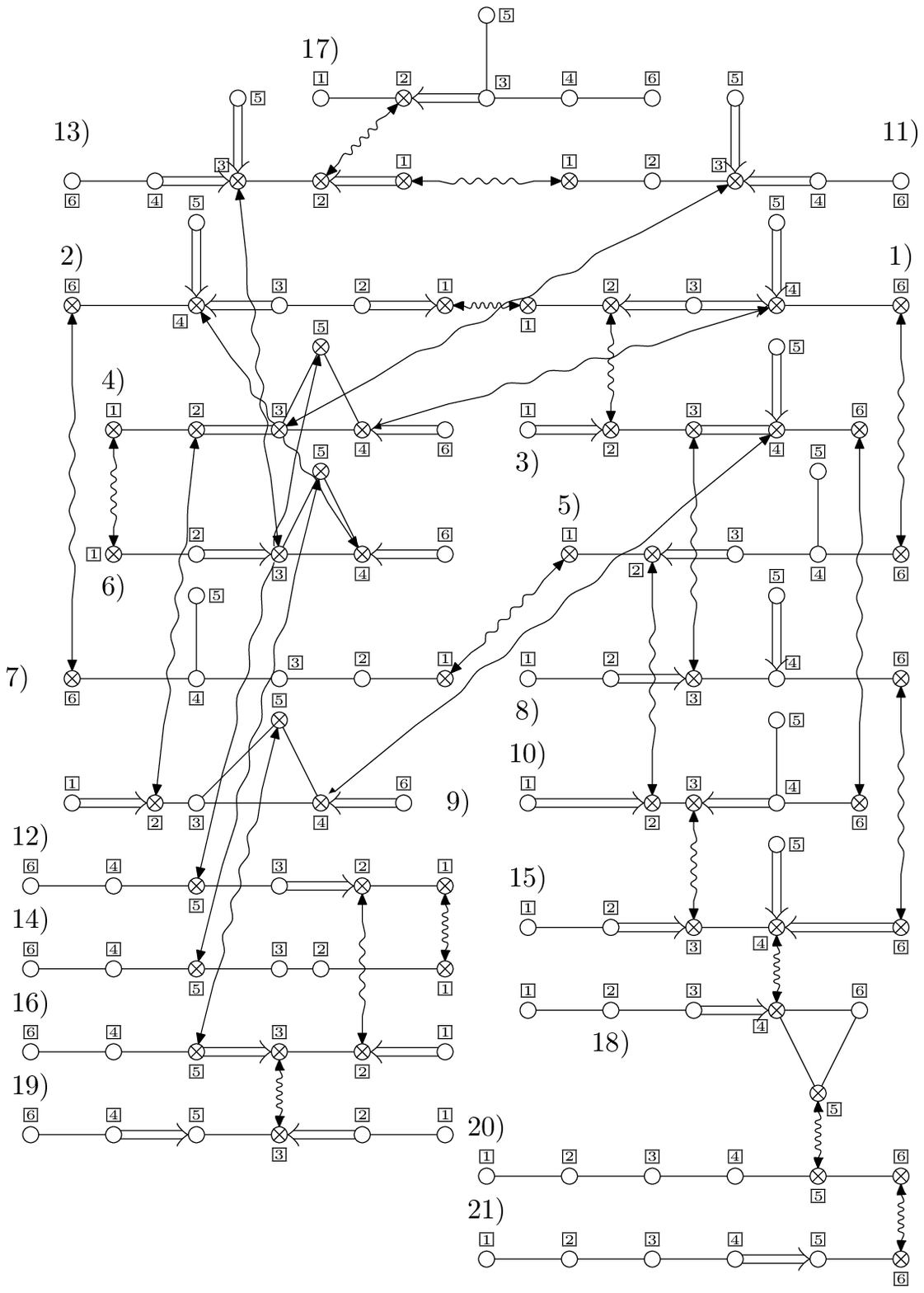}
\end{figure}

\clearpage


\normalsize

\sssec{$\fg{}(8,6)$ of $\sdim =133|56$} We have
$\fg(8,6)_\ev=\fe(7)$ and $\fg(8,6)_\od=R(\pi_1)$. \tiny
\[
\begin{pmatrix}
 -&-&-&-&-&2&3 \\
 -&-&-&-&4&1&- \\
 -&-&-&-&-&-&1 \\
 -&-&-&5&2&-&- \\
 -&6&7&4&-&-&- \\
 -&5&-&-&-&-&- \\
 8&-&5&-&-&-&- \\
 7&-&-&-&-&-&-
 \end{pmatrix}
\]

\[\begin{matrix} 1)
\begin{pmatrix}
 2&0&-1&0&0&0&0 \\
 0&2&0&-1&0&0&0 \\
 -1&0&2&-1&0&0&0 \\
 0&-1&-1&2&-1&0&0 \\
 0&0&0&-1&2&-1&0 \\
 0&0&0&0&-2&0&-1 \\
 0&0&0&0&0&-1&0
 \end{pmatrix}\quad
2) \begin{pmatrix}
 2&0&-1&0&0&0&0 \\
 0&2&0&-1&0&0&0 \\
 -1&0&2&-1&0&0&0 \\
 0&-1&-1&2&-1&0&0 \\
 0&0&0&-2&0&-1&0 \\
 0&0&0&0&-1&0&-2 \\
 0&0&0&0&0&-1&2
 \end{pmatrix} \\
\noalign{\vspace{3.5ex}} 
\boxed{3)} \arraycolsep=5pt
 \begin{pmatrix}
 2&0&-1&0&0&0&0 \\
 0&2&0&-1&0&0&0 \\
 -1&0&2&-1&0&0&0 \\
 0&-1&-1&2&-1&0&0 \\
 0&0&0&-1&2&-1&0 \\
 0&0&0&0&-1&2&-1 \\
 0&0&0&0&0&-1&0
 \end{pmatrix}\quad
4) \begin{pmatrix}
 2&0&-1&0&0&0&0 \\
 0&2&0&-1&0&0&0 \\
 -1&0&2&-1&0&0&0 \\
 0&-2&-2&0&-1&0&0 \\
 0&0&0&-1&0&-2&0 \\
 0&0&0&0&-1&2&-1 \\
 0&0&0&0&0&-1&2
 \end{pmatrix}\\
\noalign{\vspace{3.5ex}} 
5) \arraycolsep=5pt
 \begin{pmatrix}
 2&0&-1&0&0&0&0 \\
 0&0&-1&-1&0&0&0 \\
 -2&-1&0&-1&0&0&0 \\
 0&-1&-1&0&-2&0&0 \\
 0&0&0&-1&2&-1&0 \\
 0&0&0&0&-1&2&-1 \\
 0&0&0&0&0&-1&2
 \end{pmatrix}\quad
\boxed{6)} \begin{pmatrix}
 2&0&-1&0&0&0&0 \\
 0&0&-2&-2&0&0&0 \\
 -1&-1&2&0&0&0&0 \\
 0&-1&0&2&-1&0&0 \\
 0&0&0&-1&2&-1&0 \\
 0&0&0&0&-1&2&-1 \\
 0&0&0&0&0&-1&2
 \end{pmatrix} \\
\noalign{\vspace{3.5ex}} 
7) \arraycolsep=5pt
 \begin{pmatrix}
 0&0&-1&0&0&0&0 \\
 0&2&-1&0&0&0&0 \\
 -1&-2&0&-2&0&0&0 \\
 0&0&-1&2&-1&0&0 \\
 0&0&0&-1&2&-1&0 \\
 0&0&0&0&-1&2&-1 \\
 0&0&0&0&0&-1&2
 \end{pmatrix}\quad
\boxed{8)}\begin{pmatrix}
 0&0&-1&0&0&0&0 \\
 0&2&-1&0&0&0&0 \\
 -1&-1&2&-1&0&0&0 \\
 0&0&-1&2&-1&0&0 \\
 0&0&0&-1&2&-1&0 \\
 0&0&0&0&-1&2&-1 \\
 0&0&0&0&0&-1&2
 \end{pmatrix}\\
\end{matrix}
\]

\begin{landscape}
\begin{figure}[ht]\centering
\includegraphics{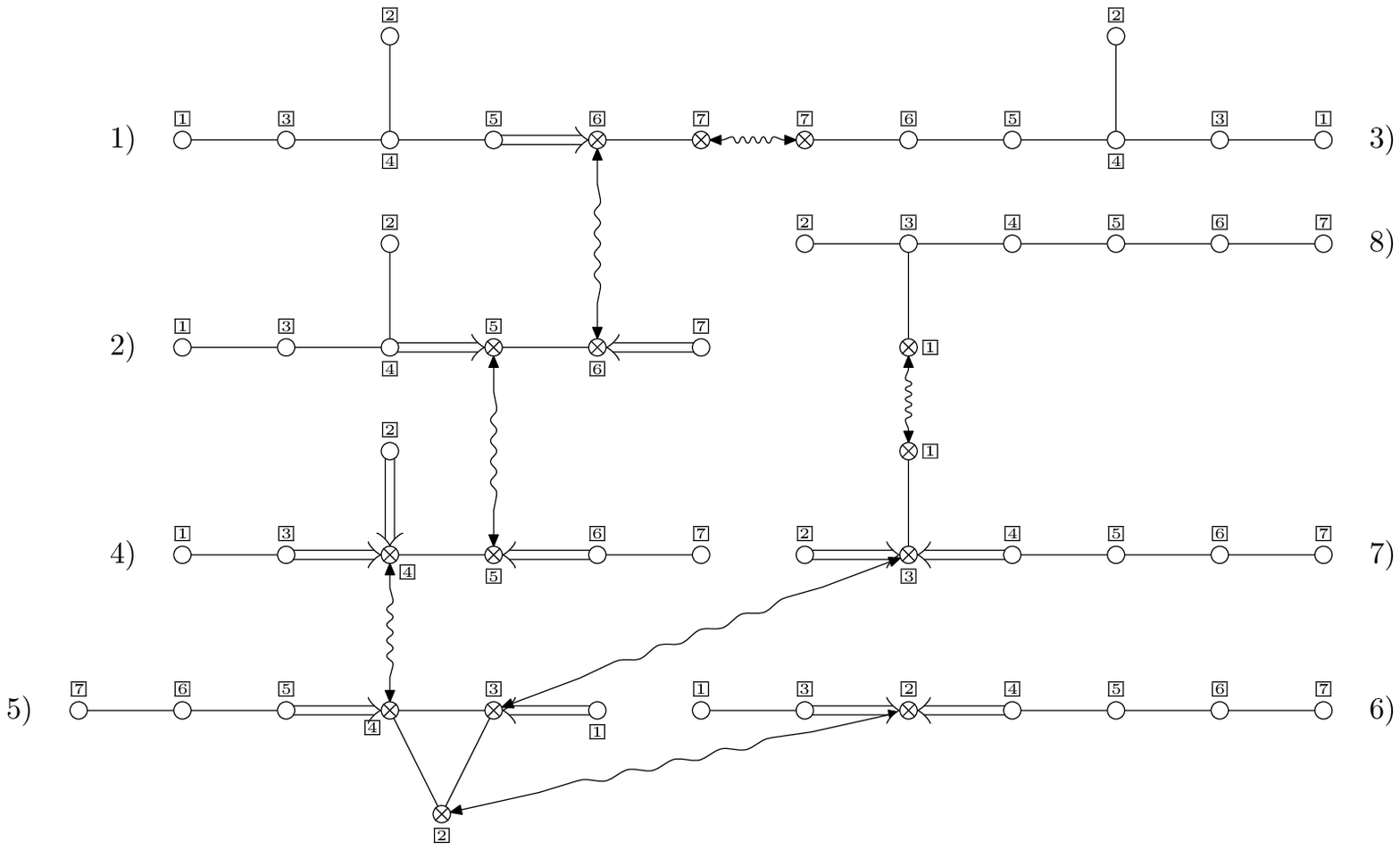}
\end{figure}
\end{landscape}



\normalsize

 {\it The Elduque superalgebra $\fel(5; 3)$: Systems of simple
roots}\label{SSel3ssr} Its superdimension is $39|32$; the even part
is $\fel(5; 3)_\ev=\fo(9)\oplus\fsl(2)$ and its odd part is
irreducible: $\fel(5; 3)_\od=R(\pi_4)\otimes \id$.

The following are all its Cartan matrices:
\begin{figure}[ht]
\parbox{.7\linewidth}{\includegraphics{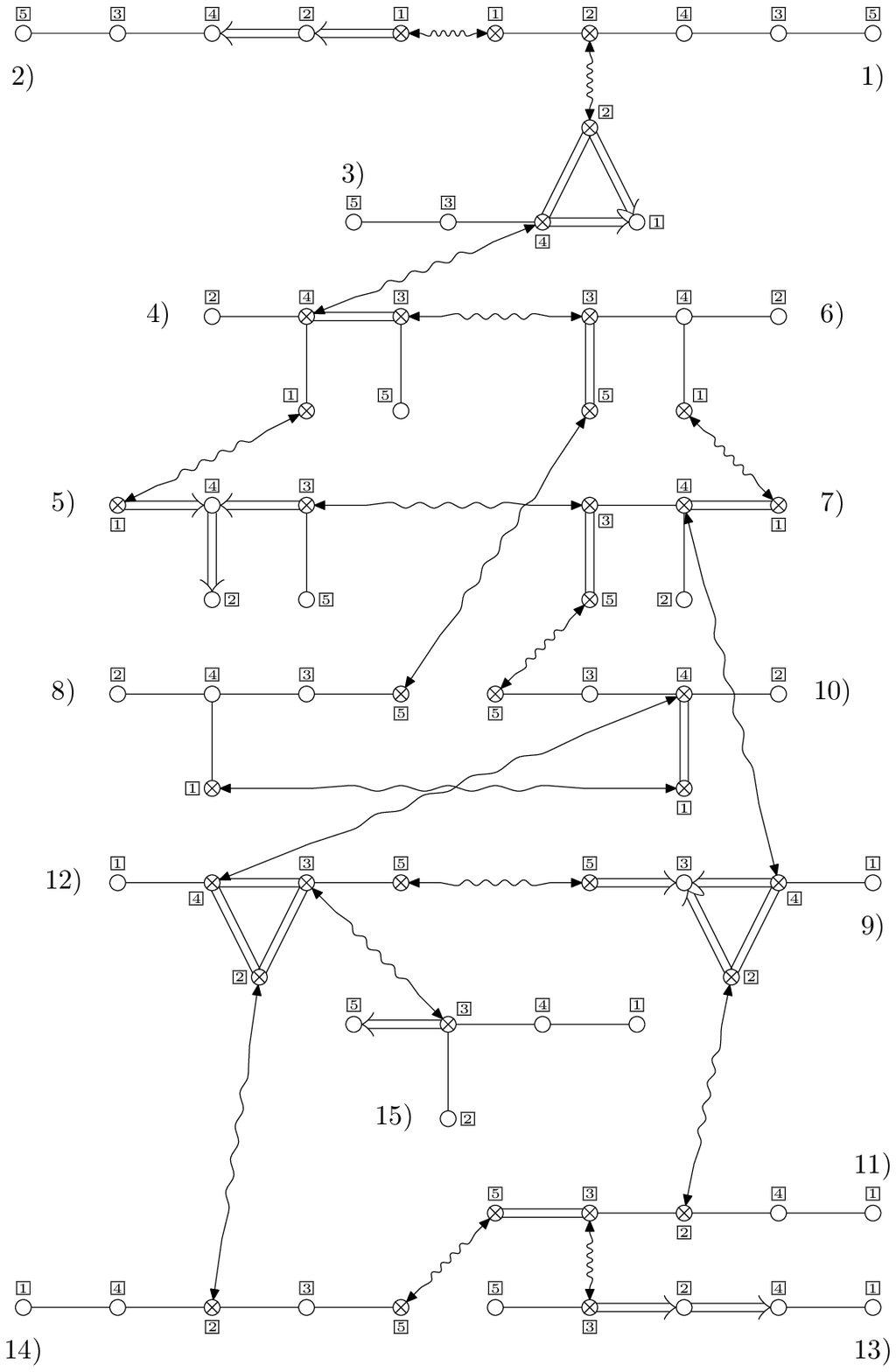}}
\end{figure}

\tiny
{
\[ \begin{matrix}
\begin{pmatrix}
2 & 3 & - & -& -\\
1 & - & - & -& -\\
- & 1 & - & 4& -\\
5 & - & 6 & 3& -\\
4 & - & 7 & -& -\\
7 & - & 4 & -& 8\\
6 & - & 5 & 9 & 10\\
10 & - & - & -& 6\\
- & 11& - &7& 12\\
8 & - & - & 12& 7\\
- & 9 & 13 & -& 14\\
- &14 & 15 & 10& 9\\
- & - & 11 & -& -\\
- & 12 & - & -& 11\\
- & - & 12 & -& -\\
\end{pmatrix}\quad \begin{matrix}
1) \;
\begin{pmatrix}
 0 & -1 & 0 & 0 & 0 \\
 -1 & 0 & 0 & -1 & 0 \\
 0 & 0 & 2 & -1 & -1 \\
 0 & -1 & -1 & 2 & 0 \\
 0 & 0 & -1 & 0 & 2
\end{pmatrix}
\quad \boxed{2)} \;
\begin{pmatrix}
 0 &-2 & 0 & 0 & 0 \\
 -1 & 2 & 0 & -2 & 0 \\
 0 & 0 & 2 & -1 & -1 \\
 0 & -1 & -1 & 2 & 0 \\
 0 & 0 & -1 & 0 & 2
\end{pmatrix}
\\
 3) \;
\begin{pmatrix}
 2 & -1 & 0 & -1 & 0 \\
-2 & 0 & 0 &-2 & 0 \\
 0 & 0 & 2 & -1 & -1 \\
-2 &-2 & -1 & 0 & 0 \\
 0 & 0 & -1 & 0 & 2
\end{pmatrix}
\quad 4) \;
\begin{pmatrix}
 0 & 0 & 0 & -1 & 0 \\
 0 & 2 & 0 & -1 & 0 \\
 0 & 0 & 0 &-2 & -1 \\
 -1 & -1 &-2 & 0 & 0 \\
 0 & 0 & -1 & 0 & 2
\end{pmatrix}
\end{matrix}
\\
 5) \;
\begin{pmatrix}
 0 & 0 & 0 &-2 & 0 \\
 0 & 2 & 0 & -1 & 0 \\
 0 & 0 & 0 &-2 & -1 \\
 -1 & -2 & -1 & 2 & 0 \\
 0 & 0 & -1 & 0 & 2
\end{pmatrix}
\quad 6) \;
\begin{pmatrix}
 0 & 0 & 0 & -1 & 0 \\
 0 & 2 & 0 & -1 & 0 \\
 0 & 0 & 0 & -1 &-2 \\
 -1 & -1 & -1 & 2 & 0 \\
 0 & 0 &-2 & 0 & 0
\end{pmatrix}
\\
7) \;
\begin{pmatrix}
 0 & 0 & 0 &-2 & 0 \\
 0 & 2 & 0 & -1 & 0 \\
 0 & 0 & 0 & -1 &-2 \\
-2 & -1 & -1 & 0 & 0 \\
 0 & 0 &-2 & 0 & 0
\end{pmatrix}
\quad 8) \;
\begin{pmatrix}
 0 & 0 & 0 & -1 & 0 \\
 0 & 2 & 0 & -1 & 0 \\
 0 & 0 & 2 & -1 & -1 \\
 -1 & -1 & -1 & 2 & 0 \\
 0 & 0 & -1 & 0 & 0
\end{pmatrix}
\quad 9) \;
\begin{pmatrix}
 2 & 0 & 0 & -1 & 0 \\
 0 & 0 &-2 &-2 & 0 \\
 0 & -1 & 2 & -1 & -1 \\
 -1 &-2 &-2 & 0 & 0 \\
 0 & 0 &-2 & 0 & 0
\end{pmatrix}
\\ 10)\;
\begin{pmatrix}
 0 & 0 & 0 &-2 & 0 \\
 0 & 2 & 0 & -1 & 0 \\
 0 & 0 & 2 & -1 & -1 \\
-2 & -1 & -1 & 0 & 0 \\
 0 & 0 & -1 & 0 & 0
\end{pmatrix}
\quad 11)\;
\begin{pmatrix}
 2 & 0 & 0 & -1 & 0 \\
 0 & 0 & -1 & -1 & 0 \\
 0 & -1 & 0 & 0 &-2 \\
 -1 & -1 & 0 & 2 & 0 \\
 0 & 0 &-2 & 0 & 0
\end{pmatrix}
\quad 12)\;
\begin{pmatrix}
 2 & 0 & 0 & -1 & 0 \\
 0 & 0 &-2 &-2 & 0 \\
 0 &-2 & 0 &-2 & -1 \\
 -1 &-2 &-2 & 0 & 0 \\
 0 & 0 & -1 & 0 & 0
\end{pmatrix}
\\ \boxed{13)}\;
\begin{pmatrix}
 2 & 0 & 0 & -1 & 0 \\
 0 & 2 & -1 & -2 & 0 \\
 0 &-2 & 0 & 0 & -1 \\
 -1 & -1 & 0 & 2 & 0 \\
 0 & 0 & -1 & 0 & 2
\end{pmatrix}
\quad 14)\;
\begin{pmatrix}
 2 & 0 & 0 & -1 & 0 \\
 0 & 0 & -1 & -1 & 0 \\
 0 & -1 & 2 & 0 & -1 \\
 -1 & -1 & 0 & 2 & 0 \\
 0 & 0 & -1 & 0 & 0
\end{pmatrix}
\quad \boxed{15)}\;
\begin{pmatrix}
 2 & 0 & 0 & -1 & 0 \\
 0 & 2 & -1 & 0 & 0 \\
 0 & -1 & 0 & -1 &-2 \\
 -1 & 0 & -1 & 2 & 0 \\
 0 & 0 & -1 & 0 & 2
\end{pmatrix}
\end{matrix}
\]
}

\clearpage

\normalsize

\sssec{$p=5$, Lie superalgebras} {\it Brown
superalgebra}\index{$\fbrj(2;5)$, Brown superalgebra}\index{Brown
superalgebra} $\fbrj(2;5)$ of $\sdim = 10|12$, recently discovered
in \cite{BGL1}, such that $\fbrj(2;5)_\ev=\fsp(4)$ and
$\fbrj(2;5)_\od=R(\pi_1+\pi_2)$ is an irreducible
$\fbrj(2;5)_\ev$-module.\footnote{To the incredulous reader: The
Cartan subalgebra of $\fsp(4)$ is generated by $h_2$ and $2h_1+h_2$.
The highest weight vector is
$x_{10}=[[x_2,\;[x_2,\;[x_1,\;x_2]]],\;[[x_1,\;x_2],\;[x_1,\;x_2]]]$
and its weight is not a multiple of a fundamental weight, but
$(1,1)$. We encounter several more instances of non-fundamental
weights in descriptions of exceptions for $p=2$.} The Lie
superalgebra $\fbrj(2;5)$ has the following Cartan matrices:
\[
\begin{pmatrix}
 2 & - \\
 1 & -
\end{pmatrix}
 \qquad 1)\ \begin{pmatrix}0&-1\\
-2&1\end{pmatrix}, \quad 2)\ \begin{pmatrix}0&-1\\
-3&2\end{pmatrix}.
\]

{\it Elduque superalgebra} $\fel(5; 5)$\index{$\fel(5; 5)$, Elduque
superalgebra}\index{Elduque superalgebra} of $\sdim=55|32$, where
$\fel(5; 5)_\ev=\fo(11)$ and $\fel(5; 5)_\od=\spin_{11}$. Its
inequivalent Cartan matrices, first described in \cite{BGL1}, are as
follows: \label{SSel5ssr}

Instead of joining nodes by four segments in the cases where
$A_{ij}=A_{ji}=1\equiv -4\mod 5$ we use one dotted segment.

\begin{figure}[ht]
\parbox{.7\linewidth}{\includegraphics{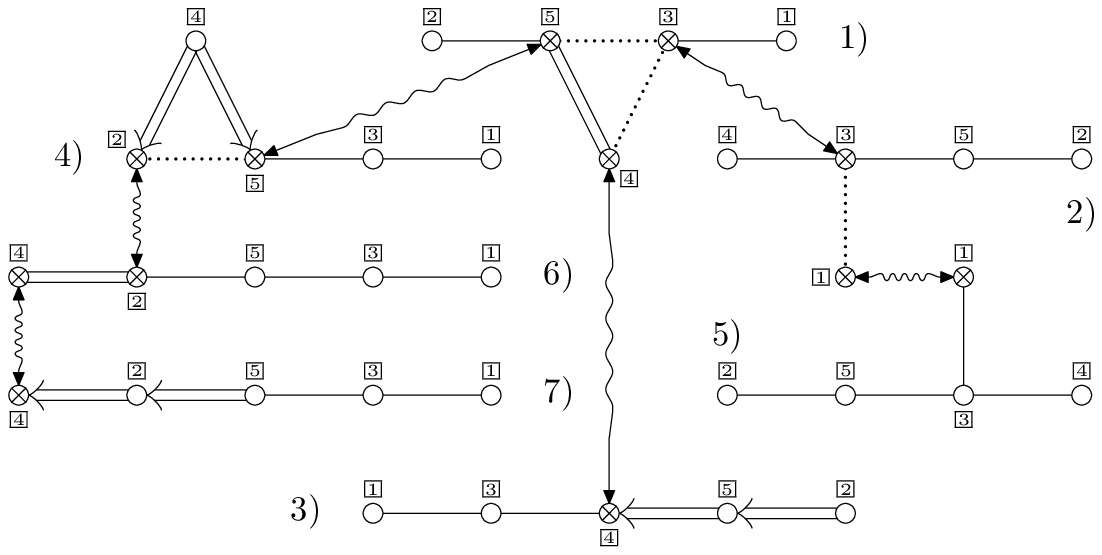}}
\end{figure}

{\tiny 
\[
\begin{matrix}
1) \begin{pmatrix}
2&0&-1&0&0 \\
 0&2&0&0&-1 \\
 -1&0&0&-4&-4 \\
 0&0&-4&0&-2 \\
 0&-1&-4&-2&0
\end{pmatrix}\quad
2) \begin{pmatrix}
 0&0&-4&0&0 \\
 0&2&0&0&-1 \\
 -4&0&0&-1&-1 \\
 0&0&-1&2&0 \\
 0&-1&-1&0&2
 \end{pmatrix}
\quad
 \boxed{3)} \begin{pmatrix}
 2&0&-1&0&0 \\
 0&2&0&0&-1 \\
 -1&0&2&-1&0 \\
 0&0&-1&0&2 \\
 0&-2&0&-1&2
 \end{pmatrix}\end{matrix}
\]
\[
\begin{matrix}
4) \begin{pmatrix}
 2&0&-1&0&0 \\
 0&0&0&2&-4 \\
 -1&0&2&0&-1 \\
 0&-1&0&2&-1 \\
 0&-4&-1&2&0
 \end{pmatrix}
\quad \boxed{5)} \begin{pmatrix}
 0&0&-1&0&0 \\
 0&2&0&0&-1 \\
 -1&0&2&-1&-1 \\
 0&0&-1&2&0 \\
 0&-1&-1&0&2
 \end{pmatrix}\quad
6) \begin{pmatrix}
 2&0&-1&0&0 \\
 0&0&0&-2&-1 \\
 -1&0&2&0&-1 \\
 0&-2&0&0&0 \\
 0&-1&-1&0&2
 \end{pmatrix}
\end{matrix}
\]
\[
\begin{matrix}
 \boxed{7)} \begin{pmatrix}
 2&0&-1&0&0 \\
 0&2&0&-1&-2 \\
 -1&0&2&0&-1 \\
 0&2&0&0&0 \\
 0&-1&-1&0&2
 \end{pmatrix}\quad
 8) \begin{pmatrix} 
 -&-&2&3&4 \\
 5&-&1&-&- \\
 -&-&-&1&- \\
 -&6&-&-&1 \\
 2&-&-&-&- \\
 -&4&-&7&- \\
 -&-&-&6&-
 \end{pmatrix}
\end{matrix}
\]
}

\clearpage

\section{Defining relations in characteristic $2$}
To save space, in what follows we omit indicating the Serre
relations; their fulfilment is assumed. Additionally there appear
relations of a new type (non-Serre relations). Here we describe
them. We have proved them analytically only for Lie (super)algebras
of $\fsl$ type and their relatives. Relations for the rest of the
(super)algebras are results of computations with \texttt{SuperLie}.

For serial Lie (super)algebras (like $\fo$, $\fo\fo$, $\fosp$,
$\fspe$), they are conjectural.

\ssec{Results}\label{non-Serre rel} Here we consider the classical
Lie algebras and superalgebras as preserving the volume element or a
non-degenerate bilinear form. We usually interpret the exceptional
Lie (super)algebras as preserving a non-integrable distribution, cf.
\cite{Shch} but here we just construct them from their Cartan
matrices.

For subalgebras of $\fgl$, we set $x_i=E_{i,i+1}, y_i=E_{i+1,i},
h_i=E_{i,i}-E_{i+1,i+1}$; the Lie sub(super)algebra $\fn$ consists
of upper-triangular (super)matrices.

\ssbegin{Theorem}\label{sl-rel} For $\fg=\fsl(n+1)$ or $\fsl(a|b)$,
where $a+b=n+1$: In characteristic $>2$, the Serre relations
$(\ref{srpm})$ define $\fn$; in characteristic $2$, the following
additional relations are required:
\begin{equation}
\label{sl} [[x_{i-1},x_i],[x_i,x_{i+1}]]=0\quad \text{ for $1<i<n$}.
\end{equation}
\end{Theorem}

\parbegin{Remark} In characteristic $p>0$, the Lie algebra $\fsl(pk)$
is not simple, since it contains the center $\fc=\{\lambda\cdot
1_{pk}\mid \lambda\in\Kee\}$. The corresponding simple Lie algebra
$\fsl(pk)/\fc$ is denoted by $\fpsl(pk)$. Since the reduction from
$\fsl(pk)$ to $\fpsl(pk)$ does not affect the structure of $\fn$,
its presentation is the same for $\fgl(pk)$, $\fsl(pk)$, and
$\fpsl(pk)$. Same applies to any other Lie (super)algebra with
non-invertible Cartan matrix. \end{Remark}

\sssbegin{Theorem}\label{e-rel} Let the nodes of the Dynkin diagram
of $\fe(8)$ be numbered as usual:
$\begin{array}{ccccccc}1&2&3&4&5&6&7\\&&&&8&&\end{array}$

For $\fg=\fe(n)$ or $\fg=\fe(n;i)$: In characteristic $2$, in the
case of $\fg=\fe(8)$, the following list of relations must be added
to the Serre relations:
\begin{equation}
\label{erel}
\begin{array}{l}
~ [[x_1, x_2], [x_2, x_3]]= 0; \\
~ [[x_2, x_3], [x_3, x_4]]= 0; \\
~ [[x_3, x_4], [x_4, x_5]]= 0; \\
~ [[x_4, x_5], [x_5, x_6]]= 0; \\
~ [[x_5, x_6], [x_6, x_7]]= 0; \\
~ [[x_4, x_5], [x_5, x_8]]= 0; \\
~ [[x_5, x_6], [x_5, x_8]]= 0; \\
~ [[x_4, [x_5, x_6]], [x_4, [x_5, x_8]]]= 0; \\
~ [[x_4, [x_5, x_6]], [x_8, [x_5, x_6]]]= 0; \\
~ [[x_4, [x_5, x_8]], [x_8, [x_5, x_6]]]= 0; \\
~ [[x_3, [x_4, [x_5, x_6]]], [x_3, [x_4, [x_5, x_8]]]]= 0; \\
~ [[x_4, [x_5, [x_6, x_7]]], [x_8, [x_5, [x_6, x_7]]]]= 0; \\
~ [[x_2, [x_3, [x_4, [x_5, x_6]]]], [x_2, [x_3, [x_4, [x_5, x_8]]]]]= 0; \\
~ [[x_1, [x_2, [x_3, [x_4, [x_5, x_6]]]]], [x_1, [x_2, [x_3, [x_4,
[x_5, x_8]]]]]]= 0.
\end{array}
\end{equation}

To obtain the corresponding lists of relations for $\fe(6)$ or
$\fe(7)$, one should delete the relations containing the \lq\lq
extra'' $x_i$ and renumber the rest of the $x_i$, i.e:

1) delete the relations containing $x_1$ for $\fe(7)$, $x_1$ and
$x_2$ for $\fe(6)$;

2) decrease all indices of the $x_i$ by $1$ for $\fe(7)$, by $2$ for
$\fe(6)$.\end{Theorem}

Proof: Direct computer calculations.

\parbegin{Remark} Here is a shorter way to describe these relations. Let {\it a
chain of nodes} for a Dynkin diagram with $n$ nodes be a sequence
$i_1,\dots,i_k$, where $k\geq 2$ and

1) $i_j\in\overline{1,n}$ for all $j=1,\dots,k$;

2) $i_j\neq i_{j'}$ for $j\neq j'$;

3) nodes with numbers $i_j$ and $i_{j+1}$ are connected for all
$j=1,\dots,k-1$.

The above non-Serre relations (both for $\fsl(n+1)$ and $\fe(n)$ can
be represented in the form
\begin{equation}\label{node_chain}
\begin{array}{l}[~[x_{i_1},[\dots,[x_{i_{k-1}},x_{i_k}]\dots]],~
[x_{i_1},[\dots,[x_{i_{k-1}},x_{i'_k}]\dots]]~]=0,\\
\text{where~} i_1,\dots,i_{k-1},i_k\text{~and~}
i_1,\dots,i_{k-1},i'_k \text{~are
two chains of nodes}\\
 \text{that differ only in the last
element.}\end{array}
\end{equation}
All the relations that can be represented in the form
(\ref{node_chain}) are necessary.
\end{Remark}

In what follows we only consider the Lie algebras $\fg(A)$; {\bf the
non-Serre relations of Lie superalgebras} $s(\fg(A))$ {\bf from
which} $\fg(A)$ {\bf can be obtained by means of forgetful functor
are the same as those of} $\fg(A)$. Theoretically, there could be
redundant ones among them, we can only conjecture (by analogy with
$\fsl$ and $\fe$ types) that no redundances occur.

\sssec{$\fg=\fo_B(2n)$} The orthogonal algebra is, by definition,
the Lie algebra of linear transformations preserving a given
non-degenerate symmetric bilinear form $B$. The bilinear form is
usually taken with the Gram matrix $1_{2n}$ or $\Pi_{2n}$. In
characteristic $>2$, these two forms are equivalent over any perfect
field. The corresponding Lie algebra has the same defining relations
as in characteristic $0$, so in this subsection we only consider
$p=2$.

It turns out (\cite{Le1}) that these two forms are not equivalent
over any ground field $\Kee$ of characteristic $2$. If $\Kee$ is
perfect, then any non-degenerate symmetric bilinear form is
equivalent to one of these two forms: It is equivalent to $\Pi_n$,
if it is zero-diagonal; otherwise, it is equivalent to $1_n$.

The orthogonal Lie algebras corresponding to these two forms (we
denote them $\fo_I(n)$ and $\fo_\Pi(n)$, respectively) are not
isomorphic and have different properties. In particular, {\bf only
$\fo_\Pi(2n)$ for $n\geq 3$ is close to an algebra with a Cartan
matrix} (same as in characteristic $0$). The corresponding algebra
$\fg^{(1)}(A)$ is $\fo\fc(2;2n)$ (i.e., the central extension of
$\fo^{(2)}_\Pi(2n)$, given by the formula \eqref{cocycle}).

\sssec{$\fo\fc(2;2n)$} The algebra $\fo^{(2)}_\Pi(2n)$ (whose
central extension is $\fo\fc(2;2n)$) consists of matrices of the
following form (where $ZD(n)$ denotes the space (Lie algebra if
$p=2$) of symmetric zero-diagonal $n\times n$\defis matrices):
\[
\left(\begin{array}{cc}A&B\\C&A^T\end{array}\right),
\quad\begin{array}{l}\text{where $A\in\fsl(n)$;}\\ B,C\in
ZD(n).\end{array}
\]
The Chevalley generators of $\fo\fc(2;2n)$ are:
\[
\begin{array}{l}
x_i=E_{i,i+1}+E_{n+i+1,n+i}\quad \text{for~}1\leq i\leq n-1; \\
x_n=E_{n-1,2n}+E_{n,2n-1};\\
y_i=x_i^T\quad\text{for~}1\leq i\leq n; \\
h_i=E_{i,i}+E_{i+1,i+1}+E_{n+i,n+i}+E_{n+i+1,n+i+1}
\quad \text{for~}1\leq i\leq n-1; \\
h_{n}=h_{n-1}+z,
\end{array}
\]
where $z$ is central element.

\parbegin{Theorem}\label{o_P_2n} In characteristic $2$, for
$\fo\fc(2;2n)$, where $n\geq 4$, the defining relations for $\fn$
are Serre relations plus the following ones:
\[
{}[[x_{i-1},x_i],[x_i,x_{i+1}]]=0\quad \text{for~}2\leq i\leq n-2;
\]
\[
\begin{array}{ll}[[x_{n-3},x_{n-2}],[x_{n-2},x_n]]=0;&\\
{}[[x_{n-2},x_{n-1}],[x_{n-2},x_n]]=0;&\\
{}[[x_{n-3},[x_{n-2},x_{n-1}]],[x_{n},[x_{n-1},x_{n-2}]]=0;&\\
{}[[x_{n-3},[x_{n-2},x_{n}]],[x_{n},[x_{n-1},x_{n-2}]]=0;&\\
\end{array}
\]
and, for $1\leq i\leq n-3$,
\[
[[x_{n-1},[x_i,[x_{i+1},\dots,[x_{n-3},x_{n-2}]\dots]]],
[x_{n},[x_i,[x_{i+1},\dots,[x_{n-3},x_{n-2}]\dots]]]]=0.
\]
\end{Theorem}
(We don't consider the case of $n=3$ in the theorem because
$\fo\fc(2;6)$ is isomorphic to $\fsl(4)$.)

\sssec{$\fg=\fo^{(1)}_I(2n)$} As shown above,  
if $n\geq 2$, then
$\fo_I(2n)\not\simeq\fo_I^{(1)}(2n)\simeq\fo_I^{(2)}(2n)$ (and if
$n=1$, then the algebra $\fo_I(2n)$ is nilpotent). So any set of
generators of $\fo_I(2n)$ contains \lq\lq extra" (as compared with
generators of $\fo_I^{(1)}(2n)$) generators $a_1,\dots,a_{2n}$. The
relations containing these generators say nothing new about the
structure of the simple (and, thus, more interesting) algebra
$\fo_I^{(1)}(2n)$. Because of this and because we want to make the
set of generators we use as small as possible, we consider the
algebra $\fo_I^{(1)}(2n)$. It consists of symmetric zero-diagonal
$2n\times 2n$-matrices. We can choose the following generators (for
the whole algebra since in this case there is no $\fn$):
\[
X_i=E_{i,i+1}+E_{i+1,i}\quad\text{for~}1\leq i\leq 2n-1.
\]

\parbegin{Theorem} The following are the defining relations for
$\fo^{(1)}_I(2n)$, $n\geq 2$:
\[
\begin{array}{ll}
~[X_i, X_j]=0&\text{for~}1\leq i,j\leq 2n-1, |i-j|\geq 2;\\
\left.\begin{array}{l}~[X_i, [X_i, X_{i+1}]]=x_{i+1} \\
~[X_{i+1}, [X_i,
X_{i+1}]]=x_{i}\end{array}\right\}&\text{for~}1\leq i\leq 2n-2;\\
~[[X_{i-1},X_i],[X_i,X_{i+1}]]=0&\text{for~} 2\leq i\leq 2n-2.
\end{array}
\]
\end{Theorem}

\begin{proof} (Sketch of.) The algebra $\fo^{(1)}_I(2n)$ is filtered:
\[
0=L_0\subset...\subset L_{2n-1},
\]
where $L_k$ consists of all symmetric zero-diagonal matrices $M$
such that $M_{ij}=0$ for all $i,j$ such that $|i-j|>k$. The
associated graded algebra is isomorphic to the algebra of
upper-triangular matrices, i.e., a maximal nilpotent subalgebra of
$\fsl(2n)$. So we can use Theorem \ref{sl-rel}.
\end{proof}

\parbegin{Remark} Presentations of the Lie algebra
$\fo^{(1)}_I(2n+1)$, where $n\geq 1$, are similar in
shape.\end{Remark}

\sssec{$\fg=\fo_B(2n+1)$}  For this algebra, again, the case of
characteristic $>2$ does not differ from the case of characteristic
$0$, so we only consider the case of characteristic $2$. Then, if
the ground field is perfect, all the non-degenerate symmetric
bilinear form over a linear space of dimension $2n+1$ are
equivalent. We choose the form $\Pi_{2n+1}$.

\paragraph{$\fg=\fo_\Pi(2n+1)$} It is easy to see that
\[\fo_\Pi(2n+1)\not \simeq \fo^{(1)}_\Pi(2n+1)\text{~ and~ $
\fo^{(1)}_\Pi(2n+1)\simeq\fo^{(2)}_\Pi(2n+1)$ for $n\geq 1$.}\]
So,
as for $\fo_I(2n)$, we consider the first derived algebra
$\fo^{(1)}_\Pi(2n+1)$. The algebra $\fo^{(1)}_\Pi(2n+1)$ consists of
matrices of the following form:
\[
\left(\begin{array}{ccc}
A&X&B\\Y^T&0&X^T\\C&Y&A^T\end{array}\right),\quad\begin{array}{l}
\text{where~}A\in\fgl(n); B,C\in ZD(n);\\
X,Y \text{~are $n$-vectors.}\end{array}
\]

This algebra has a Cartan matrix. The Chevalley generators are:
\[
\begin{array}{l}
x_i=E_{i,i+1}+E_{n+i+2,n+i+1}\quad \text{for~}1\leq i\leq n-1; \\
x_n=E_{n,n+2}+E_{n+1,2n+1};\\
y_i=x_i^T\quad\text{for~}1\leq i\leq n; \\
h_i=E_{i,i}+E_{i+1,i+1}+E_{n+i+1,n+i+1}+E_{n+i+2,n+i+2}\quad
\text{for~}1\leq i\leq n-1; \\
h_{n}=E_{n,n}+E_{2n+1,2n+1}.
\end{array}
\]

\parbegin{Theorem}\label{o_P_2n+1} In characteristic $2$, for
$\fg=\fo^{(1)}_\Pi(2n+1)$, the defining relations for $\fn$ are the
Serre relations plus the following ones:
\[
~[[x_{i-1},x_i],[x_i,x_{i+1}]]=0\quad\text{for~} 2\leq i\leq n-2.
\]
\end{Theorem}

\sssec{$\fg(2)$}  The Cartan matrix of $\fg(2)$ reduced modulo 2
coincides with Cartan matrix of $\fsl(3)$. There is, however,
another approach: Select the Chevalley basis in the Lie algebra
$\fg(2)$ as explicitly described in \cite{FH}, p. 346. Reducing the
integer structure constants reduced modulo 2 we get a {\bf simple}
Lie algebra $\fg(2)_\Kee$ (its basis is that of $\fg(2)$). This Lie
algebra is isomorphic to $\fpsl(4)$.

\sssec{$\ff(4)$} There is no $\Zee$-form of $\ff(4)$ such that the
algebra $\ff(4)_\Kee$ is still simple.

\sssec{$\fwk(3; a)$ and $\fbgl(3; a)$} The non-Serre relations are:
\[
\begin{array}{ll}
\text{For the first Cartan matrix:}&\qquad\text{For the second Cartan matrix:}\\
{}[[x_1,x_2], [x_3,
[x_1,x_2]]]=0&\qquad{}[x_2,[x_1,x_3]]=a[x_3,[x_1,x_2]]\\
{}[[x_2,x_3], [x_3, [x_1,x_2]]]=0&
\end{array}
\]

\sssec{$\fwk(4; a)$ and $\fbgl(4; a)$} The non-Serre relations are:
\[
%
\]



\sssec{$\fg=\foo(n|m)$, $p=2$} Here we consider Cartan matrix Lie
superalgebras close to some of the ortho-orthogonal algebras. There
are two kinds of such Cartan matrices:

1) The Cartan matrix
\[
%
\]
generates $\fo\fo_{I\Pi}^{(1)}(2k_\ev|2k_\od+1)$, $k_\ev+k_\od=n$
(parities of the rows of the matrix may be different; the connection
between these parities and $k_\ev,k_\od$ is described in Table
\ref{tbl}). The corresponding non-Serre relations are as in Theorem
\ref{o_P_2n+1}.

2) The Cartan matrix
\[
%
\]
generates an algebra close to $\fo\fo_{\Pi\Pi}(2k_\ev|2k_\od)$,
$k_\ev+k_\od=n$ (parities of the rows of the matrix may be
different; the connection between these parities and $k_\ev,k_\od$
is described in Table \ref{tbl}; the exact description of the Cartan
matrix Lie superalgebra is in subsection \ref{oo-oo_PP} ). The
corresponding non-Serre relations are as in Theorem \ref{o_P_2n}.

\sssec{$\fg=\fag(2)$, $p=2$} The Cartan matrices for $p=0$ are
\[
1)\; %
\]

If $p=2$, these Cartan matrices do not produce anything \lq\lq
resembling'' $\fag(2)$  since they contain $-3\equiv -1\pmod 2$ .
(In particular, the Lie superalgebra that corresponds to the
matrices 1) and 2) is isomorphic to $\fsl(1|3)$.)

We do not know an integer basis of $\fag(2)$ in which the
corresponding Lie superalgebra in characteristics $p=2$ or 3 is
simple. Elduque suggested, nevertheless,  its $p=3$ analog, see
\cite{CE}.

\sssec{$\fg=\fab(3)$, $p=2$} The Cartan matrices for $p=0$ are
\[
1)\;
%
\]
This algebra can not be constructed in $p=2$ as a Cartan matrix Lie
superalgebra. We do not know any integer basis of $\fab(3)$ such
that the corresponding Lie superalgebra in characteristics $p=2$ or
3 is simple. Elduque suggested, nevertheless, its $p=3$ analog, see
\cite{CE}.

\ssec{Defining relations: Brown (super)algebras}

\footnotesize{
\begin{arab}

\item[]


$\underline{\fbr(2;\alpha), \;\dim =10}$

\item[]
$ \arraycolsep=0pt
\renewcommand{\arraystretch}{1.2}
\begin{array}{ll}
1)\ & \ad_{x_1}^2(x_2)= 0,\quad \ad_{x_2}^3(x_1)= 0
      \end{array}
$

\end{arab}
}

\footnotesize{
\begin{arab}

\item[]


$\underline{\fbr(3), \;\dim =29}$

\item[]
$ \arraycolsep=0pt
\renewcommand{\arraystretch}{1.2}
\begin{array}{ll}
1)\ & [x_1,\;x_3]=0,\quad
(\ad_{[x_3,[x_2,x_3]]})^2([x_3,[x_1,x_2]])=0
      \end{array}
$

$ \arraycolsep=0pt
\renewcommand{\arraystretch}{1.2}
\begin{array}{ll}
2)\ & [[x_2,\;x_3],\;[x_2,\;[x_1,\;x_2]]]= 0,\quad
[[x_1,\;x_2],\;[x_3,\;[x_2,\;x_3]]]= 2\,
[[x_2,\;x_3],\;[x_3,\;[x_1,\;x_2]]]
      \end{array}
$

\end{arab}
}

\footnotesize{
\begin{arab}

\item[]


$\underline{\fbrj(2; 3), \;\dim =(10|8)}$

\item[]
$ \arraycolsep=0pt
\renewcommand{\arraystretch}{1.2}
\begin{array}{ll}
1)\ & (\mathrm{ad}_{[x_1,x_2]})^2(x_2)=0,\quad
(\mathrm{ad}_{[x_2,x_2]})^2([x_1,x_2])=0
      \end{array}
$

$ \arraycolsep=0pt
\renewcommand{\arraystretch}{1.2}
\begin{array}{ll}
2)\ & (\mathrm{ad}_{[x_2,[x_1,x_2]]})^2([x_1,x_2])=0
      \end{array}
$

$ \arraycolsep=0pt
\renewcommand{\arraystretch}{1.2}
\begin{array}{ll}
3)\ & [[x_1,\;x_1],\;[x_1,\;x_2]]=
0,\quad[x_2,\;[x_2,\;[x_1,\;x_1]]]= [[x_1,\;x_2],\;[x_1,\;x_2]]
\end{array}
$
\end{arab}
}

\footnotesize{
\begin{arab}

\item[]


$\underline{\fbrj(2; 5), \;\dim =(10|12)}$

\item[]
$ \arraycolsep=0pt
\renewcommand{\arraystretch}{1.2}
\begin{array}{ll}
1)\ & [[x_1,\;x_2],\;[[x_1,\;x_2],\;[x_2,\;x_2]]]= 3\,
[[x_2,\;[x_1,\;x_2]],\;[x_2,\;[x_1,\;x_2]]];\\
   & [[x_2,\;x_2],\;[x_1,\;x_2],\;[x_2,\;x_2]]]= 0,\quad
    [[x_2,\;[x_1,\;x_2]],\;[[x_1,\;x_2],\;[x_2,\;[x_1,\;x_2]]]]= 0
      \end{array}
$

$ \arraycolsep=0pt
\renewcommand{\arraystretch}{1.2}
\begin{array}{ll}
2)\ & [[x_2,\;[x_1,\;x_2]],\;[x_2,\;[x_2,\;[x_1,\;x_2]]]]=
0,\quad[[[x_1,\;x_2],\;[x_1,\;x_2]],\;[[x_1,\;x_2],\;[x_2,\;[x_1,\;x_2]]]]=
0
      \end{array}
$

\end{arab}
}

\section{Defining relations: Elduque superalgebras}
\footnotesize{
\begin{arab}

\item[]


$\underline{\fel(5; 5), \;\dim =(55|32)}$

\item[]
$ \arraycolsep=0pt
\renewcommand{\arraystretch}{1.2}
\begin{array}{ll}
1)\ 
   & [x_4,\;[x_3,\;x_5]]- [x_5,\;[x_3,\;x_4]]=0,\quad[[x_1,\;x_3],\;[x_3,\;x_4]]= 0,
   \quad[[x_1,\;x_3],\;[x_3,\;x_5]]= 0,\\
   &[[x_2,\;x_5],\;[x_3,\;x_5]]=
   0,\quad[[x_4,\;x_5],\;[[x_2,\;x_5],\;[x_4,\;x_5]]]=0
      \end{array}
$

$ \arraycolsep=0pt
\renewcommand{\arraystretch}{1.2}
\begin{array}{ll}
2)\ 
   & [[x_1,\;x_3],\;[x_3,\;x_4]]= 0,\quad[[x_1,\;x_3],\;[x_3,\;x_5]]= 0,\\
   &[[[x_3,\;x_4],\;[x_3,\;x_5]],\;[[x_3,\;[x_2,\;x_5]],\;
   [[x_3,\;x_4],\;[x_3,\;x_5]]]]= 0
      \end{array}
$

$ \arraycolsep=0pt
\renewcommand{\arraystretch}{1.2}
\begin{array}{ll}
3)\ 
   & [[x_3,\;x_4],\;[[x_2,\;x_5],\;[x_4,\;x_5]]]- 4\, [[x_4,\;[x_2,\;x_5]],\;
   [x_5,\;[x_3,\;x_4]]]=0,\\
   &[[x_4,\;[x_1,\;x_3]],\;[[x_3,\;x_4],\;[x_4,\;x_5]]]=
   0
\end{array}
$

$ \arraycolsep=0pt
\renewcommand{\arraystretch}{1.2}
\begin{array}{ll}
4)\ 
   &[x_4,\;
  [x_2,\;x_5]]- 3\, [x_5,\;[x_2,\;x_4]]=0,
   \quad [[x_2,\;x_5],\;[x_3\;x_5]]= 0,\\
   & [[x_5,\;[x_1,\;x_3]],\;[[x_3\;
   x_5],\;[x_4,\;x_5]]]= 0
      \end{array}
$

$ \arraycolsep=0pt
\renewcommand{\arraystretch}{1.2}
\begin{array}{ll}
5)\ 
   &  [[[x_4,\;[x_1\;
   \;x_3]],\;[x_5,\;[x_1,\;x_3]]]\;
   \;[[[x_1,\;x_3],\;[x_2,\;x_5]]\;
   \;[[x_4,\;[x_1,\;x_3]],\;[x_5,\;[x_1\;
   \;x_3]]]]]= 0
      \end{array}
$

$ \arraycolsep=0pt
\renewcommand{\arraystretch}{1.2}
\begin{array}{ll}
6)\ 
   & [[x_2,\;x_4]\;
   \;[[x_2,\;x_4],\;[x_2,\;x_5]]]= 0,\quad
   [[[x_1\;
   \;x_3],\;[x_2,\;x_5]],\;[[x_3,\;[x_2\;
   \;x_5]],\;[[x_2,\;x_4],\;[x_2,\;x_5]]]]=
   0
      \end{array}
$

$ \arraycolsep=0pt
\renewcommand{\arraystretch}{1.2}
\begin{array}{ll}
7)\ 
   & [[[x_2\;
   \;x_4],\;[x_5,\;[x_1,\;x_3]]],\;[[x_5\;
   \;[x_2,\;x_4]],\;[x_3,\;[x_2,\;[x_2\;
   \;x_5]]]]]- \\
   &2\, [[[x_2,\;x_4],\;[[x_1,\;x_3]\;
   \;[x_2,\;x_5]]],\;[[x_3,\;[x_2\;
   \;x_5]],\;[x_5,\;[x_2,\;x_4]]]]=0
      \end{array}
$

\end{arab}
}

\footnotesize{
\begin{arab}

\item[]

$\underline{\fel(5; 3), \;\dim =(39|32)}$

\item[]
$ \arraycolsep=0pt
\renewcommand{\arraystretch}{1.2}
\begin{array}{ll}
1)\ 
   & [[x_1,\;x_2],\;[[x_1,\;x_2]\;
   \;[x_2,\;x_4]]]= 0,\quad
   [[x_3,\;[x_2,\;x_4]]\;
   \;[[x_1,\;x_2],\;[x_2,\;x_4]]]=
   0,\\
&   [[[x_1\;
   \;x_2],\;[x_3,\;x_4]],\;[[x_1,\;x_2]\;
   \;[x_4,\;[x_1,\;x_2]]]]= 0
      \end{array}
$

$ \arraycolsep=0pt
\renewcommand{\arraystretch}{1.2}
\begin{array}{ll}
2)\ 
&[[x_1,\;x_2],\;[x_2,\;[x_2,\;x_4]]]=
   0,\quad [[x_2,\;[x_2,\;x_4]],\;[x_3,\;[x_2\;
   \;x_4]]]= 0,\\
&
   [[[x_1,\;x_2],\;[x_2,\;x_4]]\;
   \;[[x_1,\;x_2],\;[x_3,\;x_4]]]= 0
      \end{array}
$

$ \arraycolsep=0pt
\renewcommand{\arraystretch}{1.2}
\begin{array}{ll}
3)\ 
   &[x_2,\;[x_1,\;x_4]]- [x_4\;
   \;[x_1,\;x_2]]=0,\quad
   [[x_1\;
   \;x_4],\;[x_3,\;x_4]]= 0,\quad
   [[x_2,\;x_4]\;
   \;[x_3,\;x_4]]= 0,\\
&   [[x_3,\;[x_1,\;x_4]]\;
   \;[[x_1,\;x_2],\;[x_1,\;x_4]]]= 0
\end{array}
$

$ \arraycolsep=0pt
\renewcommand{\arraystretch}{1.2}
\begin{array}{ll}
4)\ 
   &[[x_1,\;x_4],\;[x_3,\;x_4]]= 0,\quad
   [[x_2\;
   \;x_4],\;[x_3,\;x_4]]= 0,[[x_3,\;x_4]\;
   \;[x_3,\;x_5]]= 0,\\
&
   [[x_1,\;x_4],\;[[x_1\;
   \;x_4],\;[x_2,\;x_4]]]= 0,\\
   &[[x_3,\;[x_1\;
   \;x_4]],\;[[x_1,\;x_4],\;[x_2,\;[x_1\;
   \;x_4]]]]= 0
      \end{array}
$

$ \arraycolsep=0pt
\renewcommand{\arraystretch}{1.2}
\begin{array}{ll}
5)\ 
   &[[x_3,\;x_4],\;[x_3,\;x_5]]=
   0,\quad
   [[x_1,\;x_4],\;[x_4,\;[x_2,\;x_4]]]=
   0,\\
&
   [[x_2,\;x_4],\;[x_3,\;[x_1,\;x_4]]]-2\,
   [[x_3,\;x_4],\;[x_2,\;[x_1,\;x_4]]]=0,\\
   &[[x_3\;
   \;x_4],\;[x_4,\;[x_2,\;x_4]]]= 0
      \end{array}
$

$ \arraycolsep=0pt
\renewcommand{\arraystretch}{1.2}
\begin{array}{ll}
6)\ 
   &[[x_3,\;x_4],\;[x_3,\;x_5]]= 0,\quad
   [[x_3\;
   \;[x_1,\;x_4]],\;[[x_2,\;[x_1,\;x_4]]\;
   \;[x_3,\;[x_1,\;x_4]]]]= 0,\\
  & [[x_3,\;[x_1\;
   \;x_4]],\;[[x_3,\;[x_1,\;x_4]],\;[x_3\;
   \;[x_2,\;x_4]]]]= 0
      \end{array}
$

$ \arraycolsep=0pt
\renewcommand{\arraystretch}{1.2}
\begin{array}{ll}
7)\ 
   & [[x_1\;
   \;x_4],\;[x_2,\;x_4]]= 0,\quad
   [[x_1,\;x_4]\;
   \;[x_3,\;x_4]]= 0,\quad
   [[x_3,\;x_4],\;[x_3,\;x_5]]=
   0,\\
&   [[x_3,\;x_4],\;[[x_2,\;x_4],\;[x_3\;
   \;x_4]]]= 0,\quad
   [[x_3,\;[x_1,\;x_4]],\;[[x_3\;
   \;x_4],\;[x_3,\;[x_2,\;x_4]]]]= 0
      \end{array}
$

$ \arraycolsep=0pt
\renewcommand{\arraystretch}{1.2}
\begin{array}{ll}
8)\ 
   & [[[x_1,\;x_4]\;
   \;[x_3,\;x_5]],\;[[x_2,\;[x_1,\;x_4]]\;
   \;[x_3,\;[x_1,\;x_4]]]]= 0,\\
&[[[x_1,\;x_4]\;
   \;[x_3,\;x_5]],\;[[[x_1,\;x_4],\;[x_3\;
   \;x_5]],\;[[x_2,\;x_4],\;[x_3,\;x_5]]]]= 0
      \end{array}
$

$ \arraycolsep=0pt
\renewcommand{\arraystretch}{1.2}
\begin{array}{ll}
9)\ 
   & [x_3,\;[x_2\;
   \;x_4]]- [x_4,\;[x_2,\;x_3]]=0,\quad
   [[x_1,\;x_4]\;
   \;[x_2,\;x_4]]= 0,\\
&[[x_1,\;x_4],\;[x_3,\;x_4]]=
   0,\quad
   [[x_3,\;[x_1,\;x_4]],\;[[x_2,\;x_3]\;
   \;[x_3,\;x_4]]]= 0
      \end{array}
$

$ \arraycolsep=0pt
\renewcommand{\arraystretch}{1.2}
\begin{array}{ll}
10)\ 
   & [[x_1,\;x_4],\;[x_2,\;x_4]]= 0,\quad
   [[x_1,\;x_4],\;[x_3,\;x_4]]=
   0,\\
 &  [[x_5,\;[x_3,\;x_4]],\;[[x_2,\;x_4]\;
   \;[x_3,\;x_4]]]= 0,\\
 &[[[x_1,\;x_4],\;[x_3\;
   \;x_5]],\;[[x_5,\;[x_3,\;x_4]],\;[[x_2\;
   \;x_4],\;[x_3,\;x_5]]]]= 0
      \end{array}
$

$ \arraycolsep=0pt
\renewcommand{\arraystretch}{1.2}
\begin{array}{ll}
11)\ 
   & [[x_2,\;x_3]\;
   \;[x_3,\;x_5]]= 0,\quad
   [[x_2,\;x_3],\;[[x_2\;
   \;x_3],\;[x_2,\;x_4]]]= 0,\\
&
   [[x_2,\;[x_1\;
   \;x_4]],\;[[x_2,\;x_3],\;[x_2,\;x_4]]]=
   0,\\
&[[[x_1,\;x_4],\;[x_2,\;x_3]],\;[[x_2\;
   \;x_3],\;[x_4,\;[x_2,\;x_3]]]]= 0
      \end{array}
$

$ \arraycolsep=0pt
\renewcommand{\arraystretch}{1.2}
\begin{array}{ll}
12)\ 
   & [x_3,\;[x_2,\;x_4]]-
   [x_4,\;[x_2,\;x_3]]=0,\quad
   [[x_1,\;x_4],\;[x_2\;
   \;x_4]]= 0,\quad
   [[x_1,\;x_4],\;[x_3,\;x_4]]= 0,\\
&[[x_2\;
   \;x_3],\;[x_3,\;x_5]]= 0,\quad [[x_3,\;x_4]\;
   \;[x_3,\;x_5]]= 0,\\
&
   [[[x_1,\;x_4],\;[x_3,\;x_5]],\;[[x_5,\;[x_2,\;x_3]],\;[x_5,\;[x_3,\;x_4]]]]= 0
      \end{array}
$

$ \arraycolsep=0pt
\renewcommand{\arraystretch}{1.2}
\begin{array}{ll}
13)\ 
&[x_1,\;[x_1,\;x_4]]= 0,\quad [x_2,\;x_3],\;[x_3,\;x_5]]= 0,\\
&[[x_2,\;x_3],\;[x_2,\;[x_2,\;x_4]]]= 0,\quad
   [[x_2,\;[x_1,\;x_4]],\;[x_2,\;[x_2,\;x_4]]]= 0,\\
 &  [[[x_1,\;x_4],\;[x_2,\;x_3]],\;[[x_2,\;x_3],\;[x_2,\;x_4]]]= 0
    \end{array}
$

$ \arraycolsep=0pt
\renewcommand{\arraystretch}{1.2}
\begin{array}{ll}
14)\ 
   & [[x_2,\;[x_1,\;x_4]],\;[[x_2,\;x_3],\;[x_2,\;x_4]]]= 0,\quad
[[x_5,\;[x_2\;
   \;x_3]],\;[[x_2,\;x_3],\;[x_2,\;x_4]]]=
   0,\\
&[[[x_3,\;x_5],\;[x_2,\;[x_1,\;x_4]]],\;[[x_5,\;[x_2,\;x_3]],\;[[x_2,\;x_4],\;[x_3,\;x_5]]]]=
0
    \end{array}
$

$ \arraycolsep=0pt
\renewcommand{\arraystretch}{1.2}
\begin{array}{ll}
15)\ 
&[[x_2,\;x_3],\;[x_3,\;x_5]]= 0,\quad [[x_3,\;x_4],\;[x_3,\;x_5]]=
   0,\\
&[[x_3,\;[x_1,\;x_4]]\;
   \;[[x_2,\;x_3],\;[x_3,\;x_4]]]= 0,\\
&  [[[x_1\;
   \;x_4],\;[x_3,\;x_5]],\;[[x_5,\;[x_2\;
   \;x_3]],\;[x_5,\;[x_3,\;x_4]]]]= 0
    \end{array}
$

\end{arab}
}

\section{Defining relations: Elduque and Cuhna superalgebras}

\footnotesize{
\begin{arab}

\item[]


$\underline{\fg(2, 3), \;\sdim =(11|14)}$

\item[]
$
\arraycolsep=0pt
\renewcommand{\arraystretch}{1.2}
\begin{array}{ll}
1)\ &{{[x_{2},\,[x_{1},\,x_{3}]]}-{2\,
[x_{3},\,[x_{1},\,x_{2}]]}}=0;\\ &
\left(\ad_{[x_{1},\,x_{2}]}\right)^2([x_{1},\,x_{3}])=0;\quad
\left(\ad_{[x_{1},\,x_{2}]}\right)^2([x_{2},\,x_{3}])=0.\\
      \end{array}
$

$
\arraycolsep=0pt
\renewcommand{\arraystretch}{1.2}
\begin{array}{ll}
2)\ &   \left (\ad_{[x_1,x_3] } \right )^2([x_2,x_3])=0;\quad
 \left(\ad_{[x_2,x_3]} \right )^2([x_1,x_3])=0;\\
& {{[[x_{2},\,[x_{1},\,x_{3}]],\,
      [[x_{1},\,x_{3}],\,[x_{2},
      \,[x_{1},\,x_{3}]]]]}= 0}; \\
&  {{[[x_{2},\,[x_{1},\,x_{3}]],\,
      [[x_{2},\,x_{3}],\,[x_{2},\,
      [x_{1},\,x_{3}]]]]}=0.}\\
      \end{array}
$

$
\arraycolsep=0pt
\renewcommand{\arraystretch}{1.2}
\begin{array}{ll}
3)\ &
{{[[x_{1},\,x_{3}],\,[x_{3},\,[x_{2},\,x_{3}]]]}=0};\quad
  {\left ( \mathrm{ad}_{[x_2,x_3]}\right )^2(x_3) =0}; \\
&  {{[[x_{2},\,[x_{1},\,x_{3}]]\!,\,
    [[x_{1},\,x_{3}]\!,\,[x_{2},\,x_{3}]]]}\,{=}\,0};\\
&  {{[[x_{2},\,[x_{1},\,x_{3}]]\!,\,
    [[x_{2},\,x_{3}]\!,\,[x_{2},\,x_{3}]]]}\,{=}\,0.}
\end{array}
$

$
\arraycolsep=0pt
\renewcommand{\arraystretch}{1.2}
\begin{array}{ll}
4) \ &  \left ( \mathrm{ad}_{[x_1,x_3]}\right )^2(x_3)=0; \quad
  {{[[x_{2},\,x_{3}],\,[x_{3},\,[x_{1},\,x_{3}]]]}=0};\\
&  {{[[x_{2},\,[x_{1},\,x_{3}]]\!,\,
    [[x_{1},\,x_{3}]\!,\,[x_{1},\,x_{3}]]]}\,{=}\,0};\\
&  {{[[x_{2},\,[x_{1},\,x_{3}]]\!,\,
    [[x_{1},\,x_{3}]\!,\,[x_{2},\,x_{3}]]]}\,{=}\,0.}
\end{array}
 $

 $
\arraycolsep=0pt
\renewcommand{\arraystretch}{1.2}
\begin{array}{ll}
5)\ & {{[x_{3},\,[x_{2},\,[x_{1},\,x_{3}]]]}-
    {[[x_{1},\,x_{3}],\,[x_{2},\,x_{3}]]}=0};\\
&
{{[[x_{1},\,x_{3}],\,[x_{3},\,x_{3}]]}=0};\quad
  {{[[x_{2},\,x_{3}],\,[x_{3},\,x_{3}]]}=0};\\
& {{[[x_{2},\,[x_{1},\,x_{3}]],\,[x_{3},\,
  [x_{1},\,x_{3}]]]}=0};\quad
  {{[[x_{2},\,[x_{1},\,x_{3}]],\,
  [x_{3},\,[x_{2},\,x_{3}]]]}=0.}
    \end{array}
$

\item[]

$\underline{\fg(1,6), \;\sdim =(21|14)}$

\item[]

$
\arraycolsep=0pt
\renewcommand{\arraystretch}{1.2}
\begin{array}{ll}
1)\ &
{{[x_{3},\,[x_{2},\,[x_{1},\,x_{2}]]]}-
    {[[x_{1},\,x_{2}],\,[x_{2},\,x_{3}]]}}=0;\\
&
{{[[x_{1},\,x_{2}],\,[x_{2},\,x_{2}]]}=0};
  \;\;
  {{[[x_{2},\,x_{2}],\,[x_{2},\,x_{3}]]}=
  0},\\
&
 {{[[x_{3},\,[x_{1},\,x_{2}]],\,[x_{3},\,[x_{2},\,x_{2}]]]}= 0};\quad
 {{[[x_{3},\,[x_{1},\,x_{2}]],\,[[x_{1},\,x_{2}],\,[x_{1},\,x_{2}]]]}=0.}
\end{array}
$

$
\arraycolsep=0pt
\renewcommand{\arraystretch}{1.2}
\begin{array}{ll}
2)\ &  {{[[x_{1},\,x_{2}],\,[x_{2},
  \,[x_{2},\,x_{3}]]]}= 0};\quad
  {{[[x_{2},\,x_{3}],\,[x_{2},
  \,[x_{2},\,x_{3}]]]}= 0};\\
&  {{[[x_{3},\,[x_{1},\,x_{2}]],\,
      [x_{3},\,[x_{2},\,[x_{2},\,x_{3}]]]]}=0};\\
&  {{[[x_{3},\,[x_{1},\,x_{2}]],\,
      [[x_{1},\,x_{2}],\,[x_{3},
      \,[x_{1},\,x_{2}]]]]}=0.}
      \end{array}
$

\item[]

$\underline{\fg(3, 6), \;\sdim =(36|40)}$

\item[]

$
\arraycolsep=0pt
\renewcommand{\arraystretch}{1.2}
\begin{array}{ll}
1)\ & {{[x_{4},\,[x_{3},\,[x_{2},\,x_{3}]]]}-
   {[[x_{2},\,x_{3}],\,[x_{3},\,x_{4}]]}}=0;\quad
 {{[[x_{2},\,x_{3}],\,[x_{3},\,x_{3}]]}=0};\\
 &
 {{[[x_{3},\,x_{3}],\,[x_{3},\,x_{4}]]}=0};\quad
 {{[[x_{4},\,[x_{2},\,x_{3}]],\,[x_{4},\,[x_{3},\,x_{3}]]]}=0};\\
&
 {{[[x_{4},\,[x_{2},\,x_{3}]],\,
   [[x_{2},\,x_{3}],\,[x_{2},\,x_{3}]]]}= 0};\\
& [[x_{3},\,[x_{2},\,x_{3}]],\,
   [[x_{1},\,x_{2}],\,[x_{3},\,x_{4}]]]{-}[[x_{4},\,[x_{2},\,x_{3}]],\,
   [[x_{1},\,x_{2}],\,[x_{3},\,x_{3}]]]{=}0.
\end{array}
$

$
\arraycolsep=0pt
\renewcommand{\arraystretch}{1.2}
\begin{array}{ll}
2)\ &
{{[x_{4},\,[x_{3},\,[x_{2},\,x_{3}]]]}-
   {[[x_{2},\,x_{3}],\,[x_{3},\,x_{4}]]}}=0;\\
&
  {{[[x_{1},\,x_{2}],\,[x_{2},\,x_{3}]]}=0};\quad
  {{[[x_{2},\,x_{3}],\,[x_{3},\,x_{3}]]}=0};\quad
  {{[[x_{3},\,x_{3}],\,[x_{3},\,x_{4}]]}=0};\\
&  {{[[x_{3},\,[x_{2},\,x_{3}]],\,
    [x_{4},\,[x_{2},\,x_{3}]]]}=0};\quad
  {{[[x_{4},\,[x_{2},\,x_{3}]],\,
    [x_{4},\,[x_{3},\,x_{3}]]]}=0};\\
&
{{[[[x_{1},\,x_{2}],\,[x_{3},\,x_{4}]],\,
    [[x_{1},\,x_{2}],\,[[x_{1},\,x_{2}],\,
    [x_{3},\,x_{3}]]]]}= 0.}
      \end{array}
$

$
\arraycolsep=0pt
\renewcommand{\arraystretch}{1.2}
\begin{array}{ll}
3)\ &  {{[[x_{2},\,x_{3}],\,
  [x_{3},\,[x_{3},\,x_{4}]]]}=0};\quad
  {{\left(\ad_{[x_{3},\,x_{4}]}\right)^2(x_{3})}=0}; \\
&  {{[[x_{4},\,[x_{2},\,x_{3}]],\,
      [x_{4},\,[x_{3},\,[x_{3},\,x_{4}]]]]}=0};\\
&  {{[[x_{4},\,[x_{2},\,x_{3}]],\,
      [[x_{2},\,x_{3}],\,
      [x_{4},\,[x_{2},\,x_{3}]]]]}=0};\\
&  {[[x_{4},\,[x_{2},\,x_{3}]],\,
      [[x_{3},\,x_{4}],\,[x_{3},\,
      [x_{1},\,x_{2}]]]]}-\\
&   {2\,
[[[x_{1},\,x_{2}],\,[x_{3},\,x_{4}]],\,
        [[x_{2},\,x_{3}],\,[x_{3},\,x_{4}]]]}=0.
        \end{array}
$

$
\arraycolsep=0pt
\renewcommand{\arraystretch}{1.2}
\begin{array}{ll}
5)\ & {{[[x_{1},\,x_{2}],\,[x_{2},\,x_{3}]]}=0};\quad
  {{[[x_{2},\,x_{3}],\,[x_{3},\,[x_{2},\,x_{3}]]]}=0};\\
& {{[[x_{3},\,x_{4}],\,[x_{3},\,[x_{2},\,x_{3}]]]}=0};\quad
  {{[[x_{4},\,[x_{2},\,x_{3}]],\,
    [[x_{2},\,x_{3}],\,[x_{2},\,x_{3}]]]}=0;}\\
&  {{[[x_{4},\,[x_{2},\,x_{3}]],\,
    [[x_{2},\,x_{3}],\,[x_{3},\,x_{4}]]]}=0};\\
& {{[[[x_{1},\,x_{2}],\,[x_{3},\,x_{4}]],\,
    [[x_{1},\,x_{2}],\,
    [x_{3},\,[x_{3},\,[x_{1},\,x_{2}]]]]]}=0.}
\end{array}
$

$ \arraycolsep=0pt
\renewcommand{\arraystretch}{1.2}
\begin{array}{ll}
4)\ & {{\left(\ad_{[x_{3},\,x_{4}]}\right)^2(x_{3})}\,{=}\,0};\;\;
 {{[[x_{1},\,x_{2}],\,[x_{2},\,x_{3}]]}\,{=}\,0};\;\;
 {{[[x_{2},\,x_{3}],\,[x_{3},\,[x_{3},\,x_{4}]]]}\,{=}\,0};\\
&{{[[x_{4},\,[x_{2},\,x_{3}]],\,[x_{4},\,
   [x_{3},\,[x_{3},\,x_{4}]]]]}\,{=}\,0};\\
& {{[[x_{4},\,[x_{2},\,x_{3}]],\,
   [[x_{2},\,x_{3}],\,[x_{3},\,x_{4}]]]}\,{=}\,0};\\
&
{{[[[x_{1},\,x_{2}],\,[x_{3},\,x_{4}]],\,
   [[x_{3},\,[x_{1},\,x_{2}]],\,
   [[x_{1},\,x_{2}],\,[x_{3},\,x_{4}]]]]}=0.}
\end{array}
$

$
\arraycolsep=0pt
\renewcommand{\arraystretch}{1.2}
\begin{array}{ll}
6)\ &
{{[[x_{1},\,x_{2}],\,[x_{2},\,x_{3}]]}=0};\quad
  {{[[x_{2},\,x_{3}],\,
  [[x_{2},\,x_{3}],\,[x_{3},\,x_{4}]]]}=0};\\
&
{{[[x_{3},\,x_{4}],\,[[x_{2},\,x_{3}],\,
  [x_{3},\,x_{4}]]]}\,{=}\,0};\\
& {{[[x_{4},\,[x_{2},\,x_{3}]],\,
      [[x_{2},\,x_{3}],\,[x_{4},\,
      [x_{2},\,x_{3}]]]]}\,{=}\,0};\\
& {{[[x_{4},\,[x_{2},\,x_{3}]],\,
      [[x_{3},\,x_{4}],\,[x_{4},\,
      [x_{2},\,x_{3}]]]]}=0};\\
&
{{[[[x_{1},\,x_{2}],\,[x_{3},\,x_{4}]],\,
      [[x_{3},\,[x_{1},\,x_{2}]],\,
      [[x_{1},\,x_{2}],\,[x_{3},\,x_{4}]]
       ]]}=0.}
       \end{array}
$

$
\arraycolsep=0pt
\renewcommand{\arraystretch}{1.2}
\begin{array}{ll}
7)\ & \left (\mathrm{ad}_{[x_2,x_4]}\right )^2([x_2,x_3])=0;\quad
\left (\mathrm{ad}_{[x_2,x_4]}\right )^2([x_3,x_4])=0;\\
&{{[x_{3},\,[x_{2},\,x_{4}]]}-{2\,
  [x_{4},\,[x_{2},\,x_{3}]]}=0};\\
& {{[[x_{4},\,[x_{1},\,x_{2}]],\,
      [[x_{3},\,[x_{1},\,x_{2}]],\,
      [x_{4},\,[x_{1},\,x_{2}]]]]}\,{=}\,0.}
       \end{array}
$

\end{arab}
}


$\underline{\fg(3, 3), \;\sdim=23/21|16}$

\footnotesize{
\begin{arab}

\item[]

$
\arraycolsep=0pt
\renewcommand{\arraystretch}{1.2}
\begin{array}{ll}
1)\ &
{{[[x_{3},\,x_{4}],\,[x_{3},\,[x_{2},\,x_{3}]]]}=0};\quad
  {{[[x_{3},\,[x_{1},\,x_{2}]],\,[x_{3},\,[x_{2},\,x_{3}]]]}=0};\\
&
{{[[[x_{1},\,x_{2}],\,[x_{3},\,x_{4}]],\,
    [[x_{2},\,x_{3}],\,[x_{3},\,x_{4}]]]}=0.}
\end{array}
$

$
\arraycolsep=0pt
\renewcommand{\arraystretch}{1.2}
\begin{array}{ll}
2)\ &
{{\left(\ad_{[x_{3},\,x_{4}]}\right)^2([x_{2},\,x_{3}])}=0};\quad
  {{[[x_{3},\,[x_{1},\,x_{2}]],\,
    [[x_{2},\,x_{3}],\,[x_{3},\,x_{4}]]]}=0};\\
&
{{[[[x_{1},\,x_{2}],\,[x_{3},\,x_{4}]],\,
    [[x_{3},\,x_{4}],\,[x_{4},\,[x_{2},\,x_{3}]]]]}=0.}
\end{array}
$

$
\arraycolsep=0pt
\renewcommand{\arraystretch}{1.2}
\begin{array}{ll}
3)\ &
{{[x_{3},\,[x_{2},\,x_{4}]]}-{[x_{4},\,[x_{2},\,x_{3}]]}=0};\quad
  {{[[x_{1},\,x_{2}],\,[x_{2},\,x_{3}]]}=0};\\
&
{{[[x_{1},\,x_{2}],\,[x_{2},\,x_{4}]]}=0};\quad
  {{[[x_{4},\,[x_{1},\,x_{2}]],\,
    [[x_{2},\,x_{4}],\,[x_{3},\,x_{4}]]]}=0.}
\end{array}
$

$
\arraycolsep=0pt
\renewcommand{\arraystretch}{1.2}
 \begin{array}{ll}
4)\ &
  {{[[x_{1},\,x_{2}],\,[x_{2},\,x_{3}]]}=0};\quad
  {{[[x_{1},\,x_{2}],\,[x_{2},\,x_{4}]]}=0};\\
& {{[[x_{2},\,x_{4}],\,[[x_{2},\,x_{3}],\,[x_{2},\,x_{4}]]]}\,{=}\,0};\\
& {{[[x_{4},\,[x_{1},\,x_{2}]],\,[[x_{2},\,x_{4}],\,
    [x_{4},\,[x_{2},\,x_{3}]]]]}\,{=}\,0.}
\end{array}
$

$
\arraycolsep=0pt
\renewcommand{\arraystretch}{1.2}
\begin{array}{ll}
5)\ &
  {{[[x_{4},\,[x_{1},\,x_{2}]],\,
      [[x_{3},\,[x_{1},\,x_{2}]],
      \,[x_{4},\,[x_{1},\,x_{2}]]]]}
    =0;}\\
 &   {{[[x_{4},\,[x_{1},\,x_{2}]],\,
      [[x_{4},\,[x_{1},\,x_{2}]],\,
      [x_{4},\,[x_{2},\,x_{3}]]]]}=0}.

\end{array}
$

$ \arraycolsep=0pt
\renewcommand{\arraystretch}{1.2}
\begin{array}{ll}
6)\ &
  {{[[x_{1},\,x_{2}],\,[x_{2},\,[x_{2},\,x_{3}]]]}=0};\quad
   {{[[x_{2},\,x_{4}],\,[x_{2},\,[x_{2},\,x_{3}]]]}=0;}\\
 & {{[[x_{2},\,x_{3}],\,[x_{4},\,[x_{1},\,x_{2}]]]}-
   {2\,
    [[x_{2},\,x_{4}],\,[x_{3},\,[x_{1},\,x_{2}]]]}=0.}
  \end{array}
$

$ \arraycolsep=0pt
\renewcommand{\arraystretch}{1.2}
\begin{array}{ll}
7)\ &
  {{[[x_{1},\,x_{2}],\,[x_{2},\,x_{4}]]}=0};\quad
  {{[[x_{4},\,[x_{1},\,x_{2}]],\,
      [[x_{1},\,x_{2}],\,[x_{3},\,[x_{1},\,x_{2}]]]]}=
      0};\\
  &{{[[x_{2},\,x_{3}],\,[x_{2},\,x_{4}]]}=0};
  \quad
  \left ( \mathrm{ad}_{[x_1,x_2]}\right )^2([x_2,x_3])=0; \quad
\end{array}
$

$ \arraycolsep=0pt
\renewcommand{\arraystretch}{1.2}
\begin{array}{ll}
8)\ &
  {{[x_{2},\,[x_{1},\,x_{3}]]}-{[x_{3},\,[x_{1},\,x_{2}]]}=0};\quad
  {{[[x_{1},\,x_{2}],\,[x_{2},\,x_{4}]]}=0};\\
  &
  {{[[x_{2},\,x_{3}],\,[x_{2},\,x_{4}]]}=0};\quad
  {{[[x_{4},\,[x_{1},\,x_{2}]],\,[[x_{1},\,x_{2}],\,
  [x_{1},\,x_{3}]]]}= 0.}
\end{array}
$

$ \arraycolsep=0pt
\renewcommand{\arraystretch}{1.2}
\begin{array}{ll}
9)\ &
  {{[[x_{1},\,x_{3}],\,[[x_{1},\,x_{3}],\,
    [x_{2},\,x_{3}]]]}=0};\quad
  {{[[x_{4},\,[x_{2},\,x_{3}]],\,
    [[x_{1},\,x_{3}],\,[x_{2},\,x_{3}]]]}=0};\\
 &{{[[[x_{1},\,x_{3}],\,[x_{2},\,x_{4}]],\,
    [[x_{1},\,x_{3}],\,[x_{2},\,[x_{1},\,x_{3}]]]]}=0.}
\end{array}
$

$ \arraycolsep=0pt
\renewcommand{\arraystretch}{1.2}
\begin{array}{ll}
10)\ &
   {{[[x_{1},\,x_{3}],\,[x_{3},\,[x_{2},\,x_{3}]]]}=0},
   {{[[x_{3},\,[x_{2},\,x_{3}]],\,
   [x_{4},\,[x_{2},\,x_{3}]]]}=0};\\
   & {{[[[x_{1},\,x_{3}],\,[x_{2},\,x_{3}]],\,
      [[x_{1},\,x_{3}],\,[x_{2},\,x_{4}]]]}=0}
      \end{array}
$
\item[]

$\underline{\fg(4, 3), \;\sdim =(24|26)}$

\item[]

$
\arraycolsep=0pt
\renewcommand{\arraystretch}{1.2}
\begin{array}{ll}
1)\ &
{{[[x_{1},\,x_{2}],\,[x_{2},\,[x_{2},\,x_{3}]]]}=0};\\
&
{{[[x_{2},\,x_{3}],\,[x_{4},\,[x_{1},\,x_{2}]]]}-
  {2\,[[x_{2},\,x_{4}],\,[x_{3},\,[x_{1},\,x_{2}]]]}=0};\\
&
{{[[x_{2},\,x_{4}],\,[x_{2},\,[x_{2},\,x_{3}]]]}=0};\\
&  {{[[x_{4},\,[x_{1},\,x_{2}]],\,
   [[x_{1},\,x_{2}],\,[x_{3},\,[x_{1},\,x_{2}]]]]}=0.}
\end{array}
$

$
\arraycolsep=0pt
\renewcommand{\arraystretch}{1.2}
\begin{array}{ll}
2)\ &
{{[[x_{1},\,x_{2}],\,[x_{2},\,x_{3}]]}=0};\quad
  {{[[x_{2},\,x_{3}],\,[x_{2},\,x_{4}]]}=0};\\
&
{{[[x_{2},\,x_{4}],\,[[x_{1},\,x_{2}],\,
   [x_{2},\,x_{4}]]]}\,{=}\,0};\\
 & {{[[x_{4},\,[x_{2},\,x_{3}]],\,
   [[x_{2},\,x_{4}],\,[x_{4},\,
   [x_{1},\,x_{2}]]]]}\,{=}\,0};\\
& {{[[x_{4},\,[x_{2},\,x_{3}]],\,
   [[x_{3},\,[x_{1},\,x_{2}]],\,
   [x_{4},\,[x_{2},\,x_{3}]]]]}=0.}
\end{array}
$

$ \arraycolsep=0pt
\renewcommand{\arraystretch}{1.2}
\begin{array}{ll}
3)\ &
  {{[x_{3},\,[x_{2},\,x_{4}]]}-
    {[x_{4},\,[x_{2},\,x_{3}]]}=0};\quad
  {{[[x_{1},\,x_{2}],\,[x_{2},\,x_{3}]]}=0};\\
  &
  {{[[x_{1},\,x_{2}],\,[x_{2},\,x_{4}]]}=0};\quad
  {{[[x_{4},\,[x_{1},\,x_{2}]],\,
      [[x_{2},\,x_{4}],\,[x_{3},\,x_{4}]]]}=0};\\
 & {{[[x_{4},\,[x_{1},\,x_{2}]],\,
      [[x_{3},\,[x_{1},\,x_{2}]],\,
      [x_{4},\,[x_{1},\,x_{2}]]]]}=0}.
\end{array}
$

$ \arraycolsep=0pt
\renewcommand{\arraystretch}{1.2}
\begin{array}{ll}
4)\ &
  {{[x_{3},\,[x_{2},\,x_{4}]]}-
  {[x_{4},\,[x_{2},\,x_{3}]]}=0}; \quad
  {{[[x_{1},\,x_{2}],\,[x_{2},\,x_{3}]]}=0};\\
  &{{[[x_{1},\,x_{2}],\,[x_{2},\,x_{4}]]}=
  0}; \quad
  {{[[x_{4},\,[x_{1},\,x_{2}]],\,
      [[x_{2},\,x_{4}],\,[x_{3},\,x_{4}]]]}=0};\\
 & {{[[x_{4},\,[x_{1},\,x_{2}]],\,
      [[x_{3},\,[x_{1},\,x_{2}]],\,[x_{4},\,
      [x_{1},\,x_{2}]]]]}=0}.

\end{array}
$

$ \arraycolsep=0pt
\renewcommand{\arraystretch}{1.2}
\begin{array}{ll}
5)\ &
   {{[[x_{3},\,[x_{1},\,x_{2}]],\,
      [[x_{2},\,x_{3}],\,[x_{3},\,x_{4}]]]}=0};\\
& {{[[[x_{1},\,x_{2}],\,[x_{3},\,x_{4}]],\,
      [[x_{3},\,x_{4}],\,[x_{4},\,[x_{2},\,x_{3}]]]]}=0};\\
 & {{[[[x_{1},\,x_{2}],\,[x_{3},\,x_{4}]],\,
      [[x_{3},\,[x_{1},\,x_{2}]],\,
       [[x_{1},\,x_{2}],\,[x_{3},\,x_{4}]]]]}=0.}

\end{array}
$

$ \arraycolsep=0pt
\renewcommand{\arraystretch}{1.2}
\begin{array}{ll}
6)\ & %
  {{[[x_{1},\,x_{2}],\,[x_{2},\,x_{3}]]}= 0};\quad
   \left (\mathrm{ad}_{[x_2,x_3]} \right )^2([x_3,x_4])=0; \\
   &{{[[x_{3},\,x_{4}],\,
      [[x_{2},\,x_{3}],\,[x_{3},\,x_{4}]]]}=0};\\
 &  {{[[x_{4},\,[x_{2},\,x_{3}]],\,
      [[x_{3},\,x_{4}],\,[x_{4},\,[x_{2},\,x_{3}]]]]}=0;}\\
      &
  {{[[x_{4},\,[x_{2},\,x_{3}]],\,
      [[x_{2},\,x_{3}],\,[x_{4},\,[x_{2},\,x_{3}]]]]}=0.}

\end{array}
$

$ \arraycolsep=0pt
\renewcommand{\arraystretch}{1.2}
\begin{array}{ll}
7)\ &
  {{[[x_{3},\,x_{4}],\,[x_{3},\,[x_{2},\,x_{3}]]]}=0};\quad
  {{[[x_{3},\,[x_{1},\,x_{2}]],\,
  [x_{3},\,[x_{2},\,x_{3}]]]}=0};\\
  &
  {{[[[x_{1},\,x_{2}],\,[x_{3},\,x_{4}]],\,
      [[x_{2},\,x_{3}],\,[x_{3},\,x_{4}]]]}=0};\\
      &
  {{[[[x_{1},\,x_{2}],\,[x_{3},\,x_{4}]],\,
    [[x_{1},\,x_{2}],\,[x_{3},\,[x_{3},\,[x_{1},\,x_{2}]]]]]}= 0.}

\end{array}
$

$ \arraycolsep=0pt
\renewcommand{\arraystretch}{1.2}
\begin{array}{ll}
8)\ &
  {{[[x_{1},\,x_{2}],\,[x_{2},\,x_{3}]]}=0};\quad
  {{[[x_{2},\,x_{3}],\,[x_{3},\,[x_{3},\,x_{4}]]]}=0;}\\
  &
  \left (\mathrm{ad}_{[x_3,x_4]}\right )^2 (x_3)=0; \quad
  {{[[x_{4},\,[x_{2},\,x_{3}]],\,
      [x_{4},\,[x_{3},\,[x_{3},\,x_{4}]]]]}=0};\\
      &
  {{[[x_{4},\,[x_{2},\,x_{3}]],\,
      [[x_{2},\,x_{3}],\,[x_{3},\,x_{4}]]]}=0.}

\end{array}
$

$ \arraycolsep=0pt
\renewcommand{\arraystretch}{1.2}
\begin{array}{ll}
9)\ & 
  {{[[x_{1},\,x_{2}],\,[x_{2},\,x_{3}]]}=0};\quad
  \left ( \mathrm{ad}_{[x_2,x_3]}\right )^2(x_3)=0;\\
  &
  {{[[x_{3},\,x_{4}],\,[x_{3},\,[x_{2},\,x_{3}]]]}=0};\quad
  {{[[x_{4},\,[x_{2},\,x_{3}]],\,
      [[x_{2},\,x_{3}],\,[x_{2},\,x_{3}]]]}=0};\\
      &
  {{[[x_{4},\,[x_{2},\,x_{3}]],\,
      [[x_{2},\,x_{3}],\,[x_{3},\,x_{4}]]]}=0.}
\end{array}
$

$
\arraycolsep=0pt
\renewcommand{\arraystretch}{1.2}
\begin{array}{ll}
10)\ &
  {{[x_{3},\,[x_{3},\,x_{3}]]}=0};\quad
  {{[x_{4},\,[x_{3},\,[x_{2},\,x_{3}]]]}-
  {[[x_{2},\,x_{3}],\,[x_{3},\,x_{4}]]}=0};\\
  &
  {{[[x_{1},\,x_{2}],\,[x_{2},\,x_{3}]]}=
  0};\quad
  {{[[x_{2},\,x_{3}],\,[x_{3},\,x_{3}]]}=0};\\
  &
  {{[[x_{3},\,x_{3}],\,[x_{3},\,x_{4}]]}=0};\quad
  {{[[x_{3},\,[x_{2},\,x_{3}]],\,[x_{4},\,
  [x_{2},\,x_{3}]]]}=0};\\
  &{{[[x_{4},\,[x_{2},\,x_{3}]],\,[x_{4},\,[x_{3},\,x_{3}]]]}
     =0.}

\end{array}
$

\item[]

$\underline{\fg(2, 6), \; \sdim =(35|20)}$


\item[]

$
\arraycolsep=0pt
\renewcommand{\arraystretch}{1.2}
\begin{array}{ll}
1)\ &
{{[[x_{2},\,x_{3}],\,[x_{3},\,x_{5}]]}=0};\quad
  {{[[x_{3},\,x_{4}],\,[x_{3},\,x_{5}]]}=0};\\
&  {{[[x_{3},\,[x_{1},\,x_{2}]],\,
    [[x_{2},\,x_{3}],\,[x_{3},\,x_{4}]]]}=0.}
\end{array}
$

$
\arraycolsep=0pt
\renewcommand{\arraystretch}{1.2}
\begin{array}{ll}
2)\ &
  {{[x_{3},\,[x_{2},\,x_{4}]]}-{[x_{4},\,[x_{2},\,x_{3}]]}=0};\quad
  {{[[x_{1},\,x_{2}],\,[x_{2},\,x_{3}]]}=0};\\
  &{{[[x_{1},\,x_{2}],\,[x_{2},\,x_{4}]]}=0};\quad
  {{[[x_{2},\,x_{3}],\,[x_{3},\,x_{5}]]}=0};\quad
  {{[[x_{3},\,x_{4}],\,[x_{3},\,x_{5}]]}=0.}
\end{array}
$

$ \arraycolsep=0pt
\renewcommand{\arraystretch}{1.2}
\begin{array}{ll}
3)\
  {{[[[x_{1},\,x_{2}],\,[x_{3},\,x_{5}]],\,
      [[x_{5},\,[x_{2},\,x_{3}]],\,[x_{5},\,
      [x_{3},\,x_{4}]]]]}=0.}
     \end{array}
     $

$ \arraycolsep=0pt
\renewcommand{\arraystretch}{1.2}
\begin{array}{ll}
4)\ &
  {{[[x_{1},\,x_{2}],\,[x_{2},\,x_{3}]]}=0};\quad
  {{[[x_{1},\,x_{2}],\,[x_{2},\,x_{4}]]}=0};\\
&  {{[[x_{5},\,[x_{2},\,x_{3}]],\,
  [[x_{2},\,x_{3}],\,[x_{2},\,x_{4}]]
      ]}=0.}
\end{array}
$

$ \arraycolsep=0pt
\renewcommand{\arraystretch}{1.2}
\begin{array}{ll}
5)\ &
  {{[[x_{4},\,[x_{1},\,x_{2}]],\,
  [[x_{2},\,x_{4}],\,[x_{3},\,x_{4}]]
      ]}=0};\\
      &{{[[x_{5},\,[x_{3},\,x_{4}]],\,
      [[x_{2},\,x_{4}],\,[x_{3},\,x_{4}]]]}=0.}
      \end{array}
      $

$ \arraycolsep=0pt
\renewcommand{\arraystretch}{1.2}
\begin{array}{ll}
6)\ &
  {{[[[x_{1},\,x_{2}],\,[x_{3},\,x_{5}]],\,
      [[x_{3},\,[x_{1},\,x_{2}]],\,[x_{4},\,
      [x_{1},\,x_{2}]]]]}=0.}
     \end{array}
     $


\item[]

$\underline{\fg(8, 3), \;\sdim=(55|50)}$


\item[]

$
\arraycolsep=0pt
\renewcommand{\arraystretch}{1.2}
\begin{array}{ll}
1)\ &
{{[x_{3},\,[x_{3},\,[x_{2},\,x_{3}]]]}=0};\quad
  {{[[x_{3},\,x_{4}],\,[x_{3},\,[x_{2},\,x_{3}]]]}=0};\\
& {{[[x_{3},\,[x_{1},\,x_{2}]],\,[x_{3},\,[x_{2},\,x_{3}]]]}=0};\\
&
{{[[[x_{1},\,x_{2}],\,[x_{3},\,x_{4}]],\,
    [[x_{3},\,x_{4}],\,[x_{4},\,[x_{2},\,x_{3}]]]]}=0};\\
&
{[[[x_{1},\,x_{2}],\,[x_{3},\,x_{4}]],\,
    [[x_{4},\,x_{5}],\,[x_{3},\,[x_{2},\,x_{3}]]]]}-\\
&
{[[[x_{2},\,x_{3}],\,[x_{4},\,x_{5}]],\,
    [[x_{3},\,x_{4}],\,[x_{3},\,[x_{1},\,x_{2}]]]]}=0.
\end{array}
$

$
\arraycolsep=0pt
\renewcommand{\arraystretch}{1.2}
\begin{array}{ll}
2)\ &
{{[[x_{3},\,x_{4}],\,[x_{4},\,x_{5}]]}=0};\quad
  {{[[x_{3},\,x_{4}],\,
    [x_{3},\,[x_{2},\,x_{3}]]]}=0};\\
&  {{[[x_{3},\,[x_{1},\,x_{2}]],\,
    [x_{3},\,[x_{2},\,x_{3}]]]}\,{=}\,0};
    \\
&
{{[[[x_{1},\,x_{2}],\,[x_{3},\,x_{4}]],\,
    [[x_{2},\,x_{3}],\,[x_{3},\,x_{4}]]]}\,{=}\,0};\\
&
{{[[[x_{4},\,x_{5}],\,[x_{3},\,[x_{1},\,x_{2}]]],\,
      [[x_{5},\,[x_{3},\,x_{4}]],\,
       [[x_{2},\,x_{3}],\,[x_{4},\,x_{5}]]]]}=0.}
\end{array}
$

$ \arraycolsep=0pt
\renewcommand{\arraystretch}{1.2}
\begin{array}{ll}
3)\ &
{{[[x_{3},\,x_{4}],\,[x_{4},\,x_{5}]]}=0};\quad
  {{[[x_{3},\,x_{4}],\,[[x_{2},\,x_{3}],\,[x_{3},\,x_{4}]]]}
     =0};\\
     &{{[[x_{3},\,[x_{1},\,x_{2}]],\,
      [[x_{2},\,x_{3}],\,[x_{3},\,x_{4}]]]}=0};\\
& {{[[[x_{1},\,x_{2}],\,[x_{3},\,x_{4}]],\,
      [[x_{3},\,x_{4}],\,[x_{4},\,
      [x_{2},\,x_{3}]]]]}=0};\\
 & {{[[[x_{4},\,x_{5}],\,[x_{3},\,[x_{1},\,x_{2}]]],\,
      [[x_{5},\,[x_{3},\,x_{4}]],\,
       [[x_{2},\,x_{3}],\,[x_{4},\,x_{5}]]]]}=0.}
\end{array}
$

$ \arraycolsep=0pt
\renewcommand{\arraystretch}{1.2}
\begin{array}{ll}
4)\ &
  {{[x_{3},\,[x_{2},\,x_{4}]]}-{[x_{4},\,[x_{2},\,x_{3}]]}=0};\quad
  {{[[x_{1},\,x_{2}],\,[x_{2},\,x_{3}]]}=0};\\
 & {{[[x_{1},\,x_{2}],\,[x_{2},\,x_{4}]]}=0};\quad
  {{[[x_{4},\,[x_{1},\,x_{2}]],\,
      [[x_{2},\,x_{4}],\,[x_{3},\,x_{4}]]]}=0};\\
      &
  {{[[[x_{1},\,x_{2}],\,[x_{4},\,x_{5}]],\,
      [[x_{5},\,[x_{2},\,x_{4}]],\,
      [x_{5},\,[x_{3},\,x_{4}]]]
      ]}=0.}
\end{array}
$

$ \arraycolsep=0pt
\renewcommand{\arraystretch}{1.2}
\begin{array}{ll}
5)\ &
  {{[[x_{1},\,x_{2}],\,[x_{2},\,x_{3}]]}=0};\quad
  {{[[x_{1},\,x_{2}],\,[x_{2},\,x_{4}]]}=0};\\
  &
  {{[[x_{2},\,x_{4}],\,[x_{4},\,x_{5}]]}=0};\quad
 \left (\mathrm{ad}_{[x_2,x_4]} \right )^2([x_2,x_3])=0;\\
     & {{[[x_{4},\,[x_{1},\,x_{2}]],\,
      [[x_{2},\,x_{4}],\,[x_{4},\,[x_{2},\,x_{3}]]]]}=0};\\
& {{[[[x_{1},\,x_{2}],\,[x_{4},\,x_{5}]],\,
      [[x_{5},\,[x_{2},\,x_{4}]],\,
       [[x_{2},\,x_{3}],\,[x_{4},\,x_{5}]]]]}=0.}
\end{array}
$

$ \arraycolsep=0pt
\renewcommand{\arraystretch}{1.2}
\begin{array}{ll}
6)\ &
  {{[[x_{2},\,x_{4}],\,[x_{4},\,x_{5}]]}=0};\quad
  {{[[x_{4},\,[x_{1},\,x_{2}]],\,
      [[x_{3},\,[x_{1},\,x_{2}]],\,[x_{4},\,[x_{1},\,x_{2}]]]
      ]}=0};\\
      & {{[[x_{4},\,[x_{1},\,x_{2}]],\,
      [[x_{4},\,[x_{1},\,x_{2}]],\,[x_{4},\,[x_{2},\,x_{3}]]]
      ]}=0};\\
      &{{[[[x_{1},\,x_{2}],\,[x_{4},\,x_{5}]],\,
      [[[x_{1},\,x_{2}],\,[x_{4},\,x_{5}]],\,
       [[x_{2},\,x_{3}],\,[x_{4},\,x_{5}]]]]}=0.}
       \end{array}
       $

       $
        \arraycolsep=0pt
\renewcommand{\arraystretch}{1.2}
\begin{array}{ll}
7)\ &
  {{[[x_{2},\,x_{4}],\,[x_{4},\,x_{5}]]}=0};\quad
  {{[[x_{1},\,x_{2}],\,[x_{2},\,[x_{2},\,x_{3}]]]}=0};\\
  &
  {{[[x_{2},\,x_{3}],\,[x_{4},\,[x_{1},\,x_{2}]]]}-
    {2\,
    [[x_{2},\,x_{4}],\,[x_{3},\,[x_{1},\,x_{2}]]]}=0;}\\
    &
  {{[[x_{2},\,x_{4}],\,[x_{2},\,[x_{2},\,x_{3}]]]}=0;}\\
 & {{[[[x_{1},\,x_{2}],\,[x_{4},\,x_{5}]],\,
      [[x_{5},\,[x_{2},\,x_{4}]],\,
       [[x_{2},\,x_{3}],\,[x_{4},\,x_{5}]]]]}=0.}
       \end{array}$

       $
       \arraycolsep=0pt
\renewcommand{\arraystretch}{1.2}
\begin{array}{ll}
8)\ &
  {{[[x_{1},\,x_{2}],\,[x_{2},\,x_{4}]]}=0};\quad
  {{[[x_{2},\,x_{3}],\,[x_{2},\,x_{4}]]}=0;}\\
&
{{[[x_{2},\,x_{4}],\,[x_{4},\,x_{5}]]}=0};\quad
  {{[[x_{1},\,x_{2}],\,[[x_{1},\,x_{2}],\,
  [x_{2},\,x_{3}]]]}=0;}\\
  & {{[[x_{4},\,[x_{1},\,x_{2}]],\,
      [[x_{1},\,x_{2}],\,[x_{3},\,[x_{1},\,x_{2}]]]]}=0;}\\
& {{[[[x_{1},\,x_{2}],\,[x_{4},\,x_{5}]],\,
      [[[x_{1},\,x_{2}],\,[x_{4},\,x_{5}]],\,
       [[x_{2},\,x_{3}],\,[x_{4},\,x_{5}]]]]}=0.}
      \end{array} $

      $
 \arraycolsep=0pt
\renewcommand{\arraystretch}{1.2}
\begin{array}{ll}
9)\ &
  {{[[x_{1},\,x_{2}],\,[x_{2},\,[x_{2},\,x_{3}]]]}=0};\\
 &{{[[x_{2},\,x_{3}],\,[x_{4},\,[x_{1},\,x_{2}]]]}-
    {2\,[[x_{2},\,x_{4}],\,[x_{3},\,[x_{1},\,x_{2}]]]}=0};\\
&
{{[[x_{2},\,x_{4}],\,[x_{2},\,[x_{2},\,x_{3}]]]}=0};\quad
  {{[[x_{4},\,[x_{1},\,x_{2}]],\,
      [[x_{2},\,x_{4}],\,[x_{4},\,[x_{2},\,x_{3}]]]]}=0.}
      \end{array}
      $

      $
\arraycolsep=0pt
\renewcommand{\arraystretch}{1.2}
\begin{array}{ll}
10)\ &
  {{[x_{2},\,[x_{1},\,x_{3}]]}-{[x_{3},\,[x_{1},\,x_{2}]]}=0};\quad
  {{[[x_{1},\,x_{2}],\,[x_{2},\,x_{4}]]}=0};\\
&  {{[[x_{4},\,[x_{1},\,x_{2}]],\,
      [[x_{1},\,x_{2}],\,[x_{1},\,x_{3}]]]}=0;}\quad
  {{[[x_{2},\,x_{3}],\,[x_{2},\,x_{4}]]}=0;}\\
& {{[[[x_{1},\,x_{2}],\,[x_{4},\,x_{5}]],\,
      [[[x_{1},\,x_{2}],\,[x_{4},\,x_{5}]],\,
       [[x_{2},\,x_{3}],\,[x_{4},\,x_{5}]]]]}=0.}
       \end{array}$

$  \arraycolsep=0pt
\renewcommand{\arraystretch}{1.2}
\begin{array}{ll}
11)\ &
{{[[x_{1},\,x_{2}],\,[x_{2},\,x_{4}]]}=0};\quad
  {{[[x_{2},\,x_{3}],\,[x_{2},\,x_{4}]]}=0};\\
 & \left ( \mathrm{ad}_{[x_1,x_2]}\right )^2([x_2,x_3])=0; \quad
  {{[[x_{4},\,[x_{1},\,x_{2}]],\,
      [[x_{1},\,x_{2}],\,[x_{3},\,[x_{1},\,x_{2}]]]]}=0};\\
&  {{[[x_{4},\,[x_{1},\,x_{2}]],\,
      [[x_{4},\,[x_{1},\,x_{2}]],\,
      [x_{4},\,[x_{2},\,x_{3}]]]]}=0.}
     \end{array}
     $

$ \arraycolsep=0pt
\renewcommand{\arraystretch}{1.2}
\begin{array}{ll}
12)\ &
  {{[[x_{1},\,x_{3}],\,[[x_{1},\,x_{3}],\,
  [x_{2},\,x_{3}]]]}=0};\;\;
  {{[[x_{4},\,[x_{2},\,x_{3}]],\,
      [[x_{1},\,x_{3}],\,[x_{2},\,x_{3}]]]}=0};\\
& {{[[[x_{1},\,x_{3}],\,[x_{2},\,x_{4}]],\,
      [[x_{1},\,x_{3}],\,[x_{2},\,[x_{1},\,x_{3}]]]]}=0};\\
&
{{[[[x_{4},\,x_{5}],\,[x_{2},\,[x_{1},\,x_{3}]]],\,
      [[[x_{2},\,x_{3}],\,[x_{4},\,x_{5}]],\,
       [[x_{4},\,x_{5}],\,[x_{2},\,[x_{1},\,x_{3}]]]]]}\,{=}\,0.}
   \end{array}
   $

$\arraycolsep=0pt
\renewcommand{\arraystretch}{1.2}
\begin{array}{ll}
13)\ &
  {{[x_{2},\,[x_{1},\,x_{3}]]}-{[x_{3},\,[x_{1},\,x_{2}]]}=0};\quad
  {{[[x_{1},\,x_{2}],\,[x_{2},\,x_{4}]]}=0};\\
 &{{[[x_{2},\,x_{3}],\,[x_{2},\,x_{4}]]}=0};\quad
 {{[[x_{2},\,x_{4}],\,[x_{4},\,x_{5}]]}=0};\\
 & {{[[x_{4},\,[x_{1},\,x_{2}]],\,
      [[x_{1},\,x_{2}],\,[x_{1},\,x_{3}]]]}=0};\\
&  {{[[x_{4},\,[x_{1},\,x_{2}]],\,
      [[x_{4},\,[x_{1},\,x_{2}]],\,
      [x_{4},\,[x_{2},\,x_{3}]]]]}=0.}
     \end{array}
     $

     $
\arraycolsep=0pt
\renewcommand{\arraystretch}{1.2}
\begin{array}{ll}
14)\ &
  {{[[x_{1},\,x_{3}],\,[x_{3},\,[x_{2},\,x_{3}]]]}=0};\quad
  {{[[x_{3},\,[x_{2},\,x_{3}]],\,
  [x_{4},\,[x_{2},\,x_{3}]]]}=0};\\
  &{{[[[x_{1},\,x_{3}],\,[x_{2},\,x_{3}]],\,
      [[x_{1},\,x_{3}],\,[x_{2},\,x_{4}]]]}=0};\\
&
{{[[[x_{4},\,x_{5}],\,[x_{2},\,[x_{1},\,x_{3}]]],\,
      [[x_{5},\,[x_{2},\,x_{4}]],\,
       [[x_{4},\,x_{5}],\,[x_{3},\,
       [x_{2},\,x_{3}]]]]]}=0.}
       \end{array}
       $

       $
       \arraycolsep=0pt
\renewcommand{\arraystretch}{1.2}
\begin{array}{ll}
15)\ &
  {{[[x_{2},\,x_{4}],\,[x_{4},\,x_{5}]]}=0};\quad
  {{[[x_{1},\,x_{3}],\,[[x_{1},\,x_{3}],\,
  [x_{2},\,x_{3}]]]}=0};\\
  &{{[[x_{4},\,[x_{2},\,x_{3}]],\,
      [[x_{1},\,x_{3}],\,[x_{2},\,x_{3}]]]}=0};\\
      &
  {{[[[x_{1},\,x_{3}],\,[x_{2},\,x_{4}]],\,
      [[x_{1},\,x_{3}],\,[x_{2},\,[x_{1},\,x_{3}]]]]}=0};\\
& {{[[[x_{1},\,x_{3}],\,[x_{2},\,x_{4}]],\,
      [[x_{4},\,[x_{2},\,x_{3}]],\,
       [[x_{1},\,x_{3}],\,[x_{2},\,x_{4}]]]]}=0.}
       \end{array}
       $

       $
        \arraycolsep=0pt
\renewcommand{\arraystretch}{1.2}
\begin{array}{ll}
16)\ &
     {{[x_{2},\,[x_{1},\,x_{3}]]}-
     {2\, [x_{3},\,[x_{1},\,x_{2}]]}=0}; \quad
     {{[[x_{2},\,x_{4}],\,[x_{4},\,x_{5}]]}=0};\\
 & \left (\mathrm{ad}_{[x_1,x_2]} \right )^2([x_1,x_3])=0;\quad
   \left (\mathrm{ad}_{[x_1,x_2]} \right )^2([x_2,x_3])=0.
\end{array}
     $

     $
 \arraycolsep=0pt
\renewcommand{\arraystretch}{1.2}
\begin{array}{ll}
17)\ &
  {{[[x_{2},\,x_{4}],\,[x_{4},\,x_{5}]]}=0};\quad
  {{[[x_{1},\,x_{3}],\,[x_{3},\,[x_{2},\,x_{3}]]]}=0};\\
 & {{[[x_{3},\,[x_{2},\,x_{3}]],\,
  [x_{4},\,[x_{2},\,x_{3}]]]}=0};\\
  & {{[[[x_{1},\,x_{3}],\,[x_{2},\,x_{3}]],\,
      [[x_{1},\,x_{3}],\,[x_{2},\,x_{4}]]]}=0};\\
& {{[[[x_{1},\,x_{3}],\,[x_{2},\,x_{4}]],\,
      [[x_{2},\,x_{4}],\,[x_{4},\,[x_{3},\,
      [x_{2},\,x_{3}]]]]]}=0.}
     \end{array}
     $

     $
     \arraycolsep=0pt
\renewcommand{\arraystretch}{1.2}
\begin{array}{ll}
18)\ &
  {{[[x_{2},\,x_{3}],\,[x_{2},\,x_{4}]]}=0};\quad
  {{[[x_{2},\,x_{4}],\,[x_{4},\,x_{5}]]}=0};\\
 & \left (\mathrm{ad}_{[x_1,x_3]} \right )^2([x_2,x_3])=0;\quad
   \left (\mathrm{ad}_{[x_2,x_3]}\right )^2([x_1,x_3])=0;\\
 & {{[[x_{2},\,[x_{1},\,x_{3}]],\,
     [[x_{1},\,x_{3}],\,[x_{2},\,[x_{1},\,x_{3}]]]]}=0};\\
&  {{[[x_{2},\,[x_{1},\,x_{3}]],\,
     [[x_{2},\,x_{3}],\,[x_{2},\,[x_{1},\,x_{3}]]]]}=0.}
      \end{array}
      $

      $
       \arraycolsep=0pt
\renewcommand{\arraystretch}{1.2}
\begin{array}{ll}
19)\ &
  {{[[x_{2},\,x_{3}],\,[x_{2},\,x_{4}]]}=0};\quad
  {{[[x_{2},\,x_{4}],\,[x_{4},\,x_{5}]]}=0};\\
& {{[[x_{1},\,x_{3}],\,[x_{3},\,[x_{2},\,x_{3}]]]}=0};\quad
  \left ( \mathrm{ad}_{[x_2,x_3]}\right )^2(x_3)=0; \\
& {{[[x_{2},\,[x_{1},\,x_{3}]],\,
      [[x_{1},\,x_{3}],\,[x_{2},\,x_{3}]]]}=0};\\
& {{[[x_{2},\,[x_{1},\,x_{3}]],\,
      [[x_{2},\,x_{3}],\,[x_{2},\,x_{3}]]]}=0.}
      \end{array}
      $

      $
       \arraycolsep=0pt
\renewcommand{\arraystretch}{1.2}
\begin{array}{ll}
20)\ &
  {{[[x_{2},\,x_{3}],\,[x_{2},\,x_{4}]]}=0};\quad
  \left (\mathrm{ad}_{[ x_1,x_3 ]} \right )^2(x_3)=0; \\
& {{[[x_{2},\,x_{3}],\,[x_{3},\,[x_{1},\,x_{3}]]]}=0};\quad
  {{[[x_{2},\,[x_{1},\,x_{3}]],\,[[x_{1},\,x_{3}],\,[x_{1},\,x_{3}]]]}=0};\\
& {{[[x_{2},\,[x_{1},\,x_{3}]],\,[[x_{1},\,x_{3}],\,[x_{2},\,x_{3}]]]}=0.}
\end{array}
      $

$
 \arraycolsep=0pt
\renewcommand{\arraystretch}{1.2}
\begin{array}{ll}
21)\ & {{[x_{3},\,[x_{3},\,x_{3}]]}=0};\quad
  {{[x_{3},\,[x_{2},\,[x_{1},\,x_{3}]]]}-
    {[[x_{1},\,x_{3}],\,[x_{2},\,x_{3}]]}=0};\\&
  {{[[x_{1},\,x_{3}],\,[x_{3},\,x_{3}]]}=0};\;\;
  {{[[x_{2},\,x_{3}],\,[x_{2},\,x_{4}]]}=0};\;\;
  {{[[x_{2},\,x_{3}],\,[x_{3},\,x_{3}]]}=0};\\&
  {{[[x_{2},\,[x_{1},\,x_{3}]],\,
  [x_{3},\,[x_{1},\,x_{3}]]]}=0};\quad
  {{[[x_{2},\,[x_{1},\,x_{3}]],\,[x_{3},\,[x_{2},\,x_{3}]]]}
=0.}
     \end{array}
     $


\item[]

$\underline{\fg(4, 6), \;\sdim =(66|32)}$

\item[]

$
\arraycolsep=0pt
\renewcommand{\arraystretch}{1.2}
\begin{array}{ll}
1)\ &
{{[[x_{3},\,x_{4}],\,[x_{4},\,x_{6}]]}=0};
\quad
  {{[[x_{4},\,x_{5}],\,[x_{4},\,x_{6}]]}=0};\\
&  {{[[x_{4},\,[x_{2},\,x_{3}]],\,
    [[x_{3},\,x_{4}],\,[x_{4},\,x_{5}]]]}=0.}
\end{array}
$

$
\arraycolsep=0pt
\renewcommand{\arraystretch}{1.2}
\begin{array}{ll}
2)\ &
{{[x_{4},\,[x_{3},\,x_{5}]]}-{[x_{5},\,[x_{3},\,x_{4}]]}=0};\quad
  {{[[x_{2},\,x_{3}],\,[x_{3},\,x_{4}]]}=0};\\
&
{{[[x_{2},\,x_{3}],\,[x_{3},\,x_{5}]]}=0};\quad
  {{[[x_{3},\,x_{4}],\,[x_{4},\,x_{6}]]}=0};\quad
{{[[x_{4},\,x_{5}],\,[x_{4},\,x_{6}]]}=0}.
\end{array}
$

$
\arraycolsep=0pt
\renewcommand{\arraystretch}{1.2}
\begin{array}{ll}
3)\ &
{{[[[x_{2},\,x_{3}],\,[x_{4},\,x_{6}]],\,
    [[x_{6},\,[x_{3},\,x_{4}]],\,
    [x_{6},\,[x_{4},\,x_{5}]]]]}=0}.
\end{array}
$

$\arraycolsep=0pt
\renewcommand{\arraystretch}{1.2}
\begin{array}{ll}
4)\ &
  {{[[x_{1},\,x_{2}],\,[x_{2},\,x_{3}]]}=0};\quad
  {{[[x_{2},\,x_{3}],\,[x_{3},\,x_{4}]]}=0};\quad
  {{[[x_{2},\,x_{3}],\,[x_{3},\,x_{5}]]}=0};\\&
  {{[[x_{6},\,[x_{3},\,x_{4}]],\,
  [[x_{3},\,x_{4}],\,[x_{3},\,x_{5}]]
      ]}=0.}
      \end{array}
      $

$\arraycolsep=0pt
\renewcommand{\arraystretch}{1.2}
\begin{array}{ll}
5)\ &
  {{[[x_{5},\,[x_{2},\,x_{3}]],\,
  [[x_{3},\,x_{5}],\,[x_{4},\,x_{5}]]]}=0};\\
 & {{[[x_{6},\,[x_{4},\,x_{5}]],\,[[x_{3},\,x_{5}],\,[x_{4},\,x_{5}]]]}=0.}
\end{array}
      $

$\arraycolsep=0pt
\renewcommand{\arraystretch}{1.2}
\begin{array}{ll}
6)\ &
  {{[[x_{1},\,x_{2}],\,[x_{2},\,x_{3}]]}=0};\\
&  {{[[[x_{2},\,x_{3}],\,[x_{4},\,x_{6}]],\,
     [[x_{4},\,[x_{2},\,x_{3}]],\,[x_{5},\,
     [x_{2},\,x_{3}]]]]}=0.}
\end{array}
     $

     $
\arraycolsep=0pt
\renewcommand{\arraystretch}{1.2}
\begin{array}{ll}
7)\ &
  {{[[[x_{4},\,x_{6}],\,[x_{3},\,[x_{1},\,x_{2}]]],\,
      [[[x_{1},\,x_{2}],\,[x_{3},\,x_{4}]],\,
       [[x_{1},\,x_{2}],\,[x_{3},\,x_{5}]]]]}=0.}
     \end{array}  $

\item[]

$\underline{\fg(6, 6),\;\sdim =(78|64)}$

\item[]

$
\arraycolsep=0pt
\renewcommand{\arraystretch}{1.2}
\begin{array}{ll}
1)\ &
{{[[x_{1},\,x_{2}],\,[x_{2},\,x_{3}]]}=0};\quad
  {{[[x_{3},\,x_{4}],\,[x_{4},\,x_{6}]]}=0};\quad
  {{[[x_{4},\,x_{5}],\,[x_{4},\,x_{6}]]}=0};\\
&  {{[[x_{4},\,[x_{2},\,x_{3}]],\,
    [[x_{3},\,x_{4}],\,[x_{4},\,x_{5}]]]}=0};\\
&
{{[[[x_{2},\,x_{3}],\,[x_{4},\,x_{6}]],\,
    [[x_{4},\,[x_{2},\,x_{3}]],\,[[x_{2},\,x_{3}],\,
    [x_{4},\,x_{5}]]]]}=0}.
\end{array}
$

$
\arraycolsep=0pt
\renewcommand{\arraystretch}{1.2}
\begin{array}{ll}
2)\ &
{{[[x_{3},\,x_{4}],\,[x_{4},\,x_{6}]]}=0};\quad
  {{[[x_{4},\,x_{5}],\,[x_{4},\,x_{6}]]}=0};\\
&  {{[[x_{4},\,[x_{2},\,x_{3}]],\,
    [[x_{3},\,x_{4}],\,[x_{4},\,x_{5}]]]}=0};\\
&
{{[[[x_{4},\,x_{6}],\,[x_{3},\,[x_{1},\,x_{2}]]],\,
    [[[x_{1},\,x_{2}],\,[x_{3},\,x_{4}]],\,
    [[x_{4},\,x_{5}],\,[x_{3},\,
    [x_{1},\,x_{2}]]]]]}=0}.
\end{array}
$

$
\arraycolsep=0pt
\renewcommand{\arraystretch}{1.2}
\begin{array}{ll}
7)\ &
{{[[[x_{2},\,x_{3}],\,[x_{4},\,x_{6}]],\,
    [[x_{6},\,[x_{3},\,x_{4}]],\,[x_{6},\,
    [x_{4},\,x_{5}]]]]}=0};\\
&
    [[[x_{4},\,x_{6}],\,[x_{3},\,[x_{1},\,x_{2}]]],\,
    [[[x_{4},\,x_{5}],\,[x_{3},\,[x_{1},\,x_{2}]]],\\
&
[[x_{4},\,x_{6}],\,[x_{3},\,[x_{1},\,x_{2}]]]]]=0.
\end{array}
$

$
\arraycolsep=0pt
\renewcommand{\arraystretch}{1.2}
\begin{array}{ll}
5)\ &
{{[[x_{1},\,x_{2}],\,[x_{2},\,x_{3}]]}\,{=}\,0};\\
&
  {{[[[x_{2},\,x_{3}],\,[x_{4},\,x_{6}]],\,
    [[x_{6},\,[x_{3},\,x_{4}]],\,
    [x_{6},\,[x_{4},\,x_{5}]]]]}\,{=}\,0};\\
&
{{[[[x_{2},\,x_{3}],\,[x_{4},\,x_{6}]],\,
    [[[x_{2},\,x_{3}],\,[x_{4},\,x_{5}]],\,
    [[x_{2},\,x_{3}],\,[x_{4},\,x_{6}]]]]}=0}.
\end{array}
$

$
\arraycolsep=0pt
\renewcommand{\arraystretch}{1.2}
\begin{array}{ll}
6)\ &
  {{[x_{4},\,[x_{3},\,x_{5}]]}-{[x_{5},\,[x_{3},\,x_{4}]]}=0};\quad
  {{[[x_{2},\,x_{3}],\,[x_{3},\,x_{4}]]}=0};\\
&
  {{[[x_{2},\,x_{3}],\,[x_{3},\,x_{5}]]}=0};\quad
  {{[[x_{3},\,x_{4}],\,[x_{4},\,x_{6}]]}=0};\quad
  {{[[x_{4},\,x_{5}],\,[x_{4},\,x_{6}]]}=0};\\
&
  {{[[[x_{4},\,x_{6}],\,[x_{3},\,[x_{1},\,x_{2}]]],\,
    [[[x_{1},\,x_{2}],\,[x_{3},\,x_{4}]],\,
    [[x_{1},\,x_{2}],\,[x_{3},\,x_{5}]]]]}=0}.
\end{array}
$

$ \arraycolsep=0pt
\renewcommand{\arraystretch}{1.2}
\begin{array}{ll}
3)\ &
   {{[[x_{1},\,x_{2}],\,[x_{2},\,x_{3}]]}=0};\quad
  {{[[x_{2},\,x_{3}],\,[x_{3},\,x_{4}]]}=0};\quad
  {{[[x_{3},\,x_{4}],\,[x_{4},\,x_{6}]]}=0};\\&
  {{[[x_{4},\,x_{5}],\,[x_{4},\,x_{6}]]}=0};\quad
  {{[[x_{3},\,x_{4}],\,[[x_{3},\,x_{4}],\,
  [x_{4},\,x_{5}]]]}=0};\\&
  {{[[x_{6},\,[x_{3},\,x_{4}]],\,
      [[x_{3},\,x_{4}],\,[x_{5},\,[x_{3},\,x_{4}]
      ]]]}=0.}
      \end{array}
      $

      $
       \arraycolsep=0pt
\renewcommand{\arraystretch}{1.2}
\begin{array}{ll}
4)\ &
  {{[x_{4},\,[x_{3},\,x_{5}]]}-{[x_{5},\,[x_{3},\,x_{4}]]}=0};\quad
  {{[[x_{1},\,x_{2}],\,[x_{2},\,x_{3}]]}=0};\\
& {{[[x_{2},\,x_{3}],\,[x_{3},\,x_{4}]]}=0};\quad
  {{[[x_{2},\,x_{3}],\,[x_{3},\,x_{5}]]}=0};\\
& {{[[x_{3},\,x_{4}],\,[x_{4},\,x_{6}]]}=0};\quad
  {{[[x_{4},\,x_{5}],\,[x_{4},\,x_{6}]]}=0};\\&
  {{[[[x_{2},\,x_{3}],\,[x_{4},\,x_{6}]],\,
      [[x_{4},\,[x_{2},\,x_{3}]],\,[x_{5},\,
      [x_{2},\,x_{3}]]]]}=0.}
     \end{array}
     $

      $
       \arraycolsep=0pt
\renewcommand{\arraystretch}{1.2}
\begin{array}{ll}
8)\ &
  {{[[x_{2},\,x_{3}],\,[x_{3},\,x_{4}]]}=0};\quad
  {{[[x_{3},\,x_{4}],\,[x_{4},\,[x_{4},\,x_{5}]]]}=0};\\
& {{[[x_{4},\,x_{6}],\,[x_{4},\,[x_{4},\,x_{5}]]]}=0;}\\&
  {{[[x_{4},\,x_{5}],\,[x_{6},\,[x_{3},\,x_{4}]]]}-
  {2\,[[x_{4},\,x_{6}],\,[x_{5},\,[x_{3},\,x_{4}]]]}=0}.
  \end{array}
  $

$ \arraycolsep=0pt
\renewcommand{\arraystretch}{1.2}
\begin{array}{ll}
9)\ &
  {{[x_{4},\,[x_{3},\,x_{5}]]}-{[x_{5},\,[x_{3},\,x_{4}]]}=0};\\
& {{[[x_{1},\,x_{2}],\,[x_{2},\,x_{3}]]}=0};\quad
  {{[[x_{3},\,x_{4}],\,[x_{4},\,x_{6}]]}=0};\\&
  {{[[x_{4},\,x_{5}],\,[x_{4},\,x_{6}]]}=0};\quad
  {{[[x_{6},\,[x_{3},\,x_{4}]],\,
  [[x_{3},\,x_{4}],\,[x_{3},\,x_{5}]]
      ]}=0.}
      \end{array}
      $

      $
      \arraycolsep=0pt
\renewcommand{\arraystretch}{1.2}
\begin{array}{ll}
10)\ &
  {{[[x_{1},\,x_{2}],\,[x_{2},\,x_{3}]]}=0};\quad
  {{[[x_{2},\,x_{3}],\,[x_{3},\,x_{4}]]}=0};\\&
  {{[[x_{6},\,[x_{3},\,x_{4}]],\,
      [[x_{5},\,[x_{3},\,x_{4}]],\,
      [x_{6},\,[x_{3},\,x_{4}]]]]}=0};\\&
      {{[[x_{6},\,[x_{3},\,x_{4}]],\,
      [[x_{6},\,[x_{3},\,x_{4}]],\,
      [x_{6},\,[x_{4},\,x_{5}]]]]}=0.}
     \end{array}
     $

     $
       \arraycolsep=0pt
\renewcommand{\arraystretch}{1.2}
\begin{array}{ll}
11)\ &
  {{[[x_{2},\,x_{3}],\,[x_{3},\,x_{4}]]}=0};\quad
  {{[[x_{2},\,x_{3}],\,[x_{3},\,x_{5}]]}=0};\\&
  {{[[x_{6},\,[x_{3},\,x_{4}]],\,
  [[x_{3},\,x_{4}],\,[x_{3},\,x_{5}]]
      ]}=0};\\&
      {{[[[x_{2},\,x_{3}],\,[x_{4},\,x_{6}]],\,
      [[x_{4},\,[x_{2},\,x_{3}]],\,
      [x_{5},\,[x_{2},\,x_{3}]]]]}
     =0.}
     \end{array}
     $

     $
      \arraycolsep=0pt
\renewcommand{\arraystretch}{1.2}
\begin{array}{ll}
12)\ &
  {{[[x_{1},\,x_{2}],\,[x_{2},\,x_{3}]]}=0};\quad
  {{[[x_{5},\,[x_{2},\,x_{3}]],\,
  [[x_{3},\,x_{5}],\,[x_{4},\,x_{5}]]
      ]}=0};\\&
      {{[[x_{6},\,[x_{4},\,x_{5}]],\,
      [[x_{3},\,x_{5}],\,[x_{4},\,x_{5}]]]}=0};\\&
  {{[[[x_{4},\,x_{6}],\,[x_{5},\,[x_{2},\,x_{3}]]],\,
      [[x_{5},\,[x_{2},\,x_{3}]],\,
      [[x_{2},\,x_{3}],\,[x_{4},\,x_{5}]]
       ]]}=0.}
       \end{array}
       $

       $
        \arraycolsep=0pt
\renewcommand{\arraystretch}{1.2}
\begin{array}{ll}
13)\ &
  {{[[x_{1},\,x_{2}],\,[x_{2},\,x_{3}]]}=0};\quad
  {{[[x_{2},\,x_{3}],\,[x_{3},\,x_{4}]]}=0};\\
& {{[[x_{2},\,x_{3}],\,[x_{3},\,x_{5}]]}=0;}\quad
  {{[[x_{6},\,[x_{3},\,x_{4}]],\,
    [[x_{3},\,x_{4}],\,[x_{3},\,x_{5}]]]}=0;}\\
& {{[[[x_{4},\,x_{6}],\,[x_{3},\,[x_{1},\,x_{2}]]],\,
    [[[x_{1},\,x_{2}],\,[x_{3},\,x_{4}]],\,
    [[x_{1},\,x_{2}],\,[x_{3},\,x_{5}]]]]}=0.}
\end{array}
       $

       $
       \arraycolsep=0pt
\renewcommand{\arraystretch}{1.2}
\begin{array}{ll}
14)\ &
  {{[[x_{5},\,[x_{2},\,x_{3}]],\,[[x_{3},\,x_{5}],\,[x_{4},\,x_{5}]]]}=0};\\&
  {{[[x_{6},\,[x_{4},\,x_{5}]],\,[[x_{3},\,x_{5}],\,[x_{4},\,x_{5}]]]}=0};\\&
    [[[x_{3},\,[x_{1},\,x_{2}]],\,[x_{6},\,[x_{4},\,x_{5}]]],\,[[[x_{1},\,x_{2}],\,[x_{3},\,x_{5}]],\,\\
& 
[[x_{4},\,x_{5}],\,[x_{3},\,[x_{1},\,x_{2}]]]]]= 0.
\end{array}
       $

       $
        \arraycolsep=0pt
\renewcommand{\arraystretch}{1.2}
\begin{array}{ll}
15)\ &
  {{[[x_{2},\,x_{3}],\,[x_{3},\,x_{4}]]}=0};\quad
  {{[[x_{3},\,x_{4}],\,[x_{4},\,x_{5}]]}=0};\\
& {{[[x_{3},\,x_{4}],\,[x_{4},\,x_{6}]]}=0};\quad
  {{[[x_{4},\,x_{6}],\,[[x_{4},\,x_{5}],\,[x_{4},\,x_{6}]]]}=0};\\
& {{[[x_{6},\,[x_{3},\,x_{4}]],\,[[x_{4},\,x_{6}],\,
    [x_{6},\,[x_{4},\,x_{5}]]]]}=0.}
\end{array}
      $

      $
        \arraycolsep=0pt
\renewcommand{\arraystretch}{1.2}
\begin{array}{ll}
16)\ &
  {{[[x_{1},\,x_{2}],\,[x_{2},\,x_{3}]]}=0};\quad
  {{[[x_{2},\,x_{3}],\,[x_{3},\,x_{5}]]}=0};\\
&
{{[[x_{3},\,x_{5}],\,[[x_{3},\,x_{5}],\,[x_{4},\,x_{5}]]]}=0};\;\;
  {{[[x_{6},\,[x_{4},\,x_{5}]],\,[[x_{3},\,x_{5}],\,[x_{4},\,x_{5}]]]}=0;}\\&
  {{[[[x_{3},\,x_{5}],\,[x_{4},\,x_{6}]],\,[[x_{3},\,x_{5}],\,[x_{4},\,[x_{3},\,x_{5}]]]]}=0.}
\end{array}
      $

      $
      \arraycolsep=0pt
\renewcommand{\arraystretch}{1.2}
\begin{array}{ll}
17)\ &
  {{[[x_{1},\,x_{2}],\,[x_{2},\,x_{3}]]}=0};\\
& {{[[[x_{2},\,x_{3}],\,[x_{4},\,x_{6}]],\,[[x_{4},\,[x_{2},\,x_{3}]],\,
    [x_{5},\,[x_{2},\,x_{3}]]]]}=0};\\&
  {{[[[x_{4},\,x_{6}],\,[x_{3},\,[x_{1},\,x_{2}]]],\,
    [[[x_{1},\,x_{2}],\,[x_{3},\,x_{4}]],\,[[x_{1},\,x_{2}],\,[x_{3},\,x_{5}]]]]}=0.}
\end{array}
$

 $
      \arraycolsep=0pt
\renewcommand{\arraystretch}{1.2}
\begin{array}{ll}
18)\ &
  {{[x_{5},\,[x_{4},\,x_{6}]]}-{[x_{6},\,[x_{4},\,x_{5}]]}=0};\quad
  {{[[x_{3},\,x_{4}],\,[x_{4},\,x_{5}]]}=0};\\&
  {{[[x_{3},\,x_{4}],\,[x_{4},\,x_{6}]]}=0};\quad
  {{[[x_{6},\,[x_{3},\,x_{4}]],\,
  [[x_{4},\,x_{6}],\,[x_{5},\,x_{6}]]
      ]}=0.}
\end{array}%
\vspace{-3mm}
     $

\end{arab}

}

\clearpage

\normalsize

\section{Proofs for $p=2$: Lie algebras}

\ssec{$\fg=\fsl(n+1)$, $p=2$} The elements $E_{ij}$, where $1\leq
i<j\leq n+1$, form a basis of the algebra $\fn$. In particular,
$x_i=E_{i,i+1}$. Clearly, we have
\[
{}[E_{ij},E_{kl}]=\delta_{jk}E_{il}+\delta_{il}E_{kj}.
\]

Let $\fhov$ be the algebra of diagonal matrices. The elements
$E_{ii}$, where $1\leq i\leq n+1$, form a basis of $\fhov$. Let the
$\omega_i$ be the dual basis elements.

We consider the weights of $\fn$ with respect to $\fhov$. The weight
of $E_{ij}$ is equal to $\omega_i+\omega_j$.

Recall several facts about homology.

\sssbegin{Lemma}\label{l1} Set
\[M_c=\{E_{i_1j_1},\dots,E_{i_mj_m}\}\text{~ for a basic chain
$c=E_{i_1j_1}\ww\dots\ww E_{i_mj_m}$.}\] If for any $E_{ij}\in M_c$
and any $k$ such that $i<k<j$, at least one of the elements $E_{ik}$
and $E_{kj}$ lies in $M_c$, then $c$ can not appear with non-zero
coefficient in decomposition of a boundary with respect to basic
chains.\end{Lemma}

\begin{proof} Clearly, it suffices to show that $c$ can not
appear with non-zero coefficient in the decomposition of the
differential of a basic chain with respect to  basic chains. It
follows from the formula for the differential $d$ that any basic
chains that appears with non-zero coefficient in decomposition of
the differential of a basic chain $F$ with respect to  basic chains,
can be obtained from $F$ by replacing $E_{ik}$ and $E_{kj}$ by
$E_{ij}$ for some $i,j,k$. If $c$ satisfies the hypothesis of the
Lemma, then $c$ can not be obtained in such a way from any $F$.
\end{proof}

The elements of $C_2(\fn;\Kee)$ have weights of two types:
$\omega_i+\omega_j$ and $\omega_i+\omega_j+\omega_k+\omega_l$.
Consider them:

\underline{I. A weight $\alpha=\omega_i+\omega_j$, where $1\leq
i<j\leq n+1$.} The following chains form a basis of $C\nc_\alpha$:
\[
\begin{array}{lll}
E_{ik}\ww E_{kj},&i<k<j;&d(E_{ik}\ww E_{kj})=E_{ij};\\
E_{ki}\ww E_{kj},&1\leq k<i;&d(E_{ki}\ww E_{kj})=0;\\
E_{ik}\ww E_{jk},&j<k\leq n+1;&d(E_{ik}\ww E_{jk})=0. \end{array}
\]

Thus, the following cycles form a basis of $C\nc_\alpha$:
\[
\begin{array}{ll}
E_{ik}\ww E_{kj}+E_{i,k+1}\ww E_{k+1,j},&i<k<j-1;\\
E_{ik}\ww E_{jk},&j<k\leq n+1 . \end{array}
\]
We consider them:

1) $E_{ik}\ww E_{kj}+E_{i,k+1}\ww E_{k+1,j}=d(E_{ik}\ww E_{k,k+1}\ww
E_{k+1,j})$, so this is a boundary.

2) $E_{ki}\ww E_{kj}$, where $1\leq k<i$; in this case, we consider
three subcases:

a) $j-i>1$: In this case, $E_{ki}\ww E_{kj}=d(E_{ki}\ww E_{k,j-1}\ww
E_{j-1,j})$.

b) $i-k>1$: In this case, $E_{ki}\ww E_{kj}=d(E_{k,i-1}\ww
E_{i-1,i}\ww E_{kj})$.

c) $i-k=j-i=1$, i.e., $i=k+1;~ j=k+2$. In this case, according to
Lemma \ref{l1}, the basic chain $E_{ki}\ww E_{kj}$ can not appear
with non-zero coefficient in decomposition of a boundary with
respect to basic chains; so this is a non-trivial cycle. It gives us
the relation
\[
{}[E_{k,k+1},E_{k,k+2}]=0,\qquad\text{i.e.,~} [x_k,[x_k,
x_{k+1}]]=0.
\]
Here $k\in\overline{1, n-1}$.

3) This case is completely analogous to the previous one; it gives
us the relation
\[
{}[x_k,[x_{k-1},x_k]]=0,
\]
where $k\in\overline{2, n}$.

\underline{II. A weight
$\alpha=\omega_i+\omega_j+\omega_i+\omega_j$, where $1\leq
i<j<k<l\leq n+1$.} Clearly, the space $C\nc_\alpha$ has the
following basis:
\[
c_{\alpha,1}=E_{ij}\ww E_{kl},\qquad c_{\alpha,2}=E_{ik}\ww
E_{jl},\qquad c_{\alpha,3}=E_{il}\ww E_{jk}.
\]
All this three chains are cycles, i.e., $Z\nc_\alpha=C\nc_\alpha$.
Here we have three subcases:

1) $j-i>1$. Then
\[
\begin{array}{l}
c_{\alpha,1}=d(E_{i,i+1}\ww E_{i+1,j}\ww E_{kl});\\
c_{\alpha,2}=d(E_{i,i+1}\ww E_{i+1,k}\ww E_{jl});\\
c_{\alpha,3}=d(E_{i,i+1}\ww E_{i+1,l}\ww E_{jk}).\end{array}
\]

2) $l-k>1$. Then, similarly to the previous case,
\[
\begin{array}{l}
c_{\alpha,1}=d(E_{ij}\ww E_{k,l-1}\ww E_{l-1,l});\\
c_{\alpha,2}=d(E_{ik}\ww E_{j,l-1}\ww E_{l-1,l});\\
c_{\alpha,3}=d(E_{jk}\ww E_{i,l-1}\ww E_{l-1,l}).\end{array}
\]

3) $j-i=l-k=1$, i.e., $j=i+1; l=k+1$. Then, from Lemma \ref{l1},
$c_{\alpha,1}$ is a non-trivial cycle. It gives the relation
\[
{}[E_{i,i+1},E_{k,k+1}]=0,\qquad\text{i.e.,~}[x_i,x_k]=0.
\]
Here $i, k\in\overline{1, n}$, and $k-i\geq 2$.

For the other cycles, we need to consider the two subcases:

a) $k-j>1$. Then \[ c_{\alpha,2}=d(E_{i,k-1}\ww E_{k-1,k}\ww
E_{jl});\quad c_{\alpha,3}=d(E_{il}\ww E_{j,k-1}\ww E_{k-1,k}).
\]

b) $k-j=1$, i.e., $i=j-1$; $k=j+1$; $l=i+2$. It is easy to see (like
in the proof of Lemma \ref{l1}) that the only two chains such that
$c_{\alpha,2}$ or $c_{\alpha,3}$ appear with non-zero coefficients
in the decomposition of their differentials with respect to  basic
chains are
\[
E_{j-1,j}\ww E_{j,j+1}\ww E_{j,j+2}\quad \text{and}\quad
E_{j-1,j+1}\ww E_{j,j+1}\ww E_{j+1,j+2}.
\]
The differentials of both these chains are equal to
$c_{\alpha,2}+c_{\alpha,3}$. So we can consider one of the chains
$c_{\alpha,2}$ or $c_{\alpha,3}$ as a non-trivial cycle. The cycle
$c_{\alpha,2}$ gives the relation
\[
{}[E_{j-1,j+1},E_{j,j+2}]=0,\qquad\text{i.e.,~}
[[x_{j-1},x_j],[x_j,x_{j+1}]]=0,
\]
and $c_{\alpha,2}$ gives an equivalent (taking other relations into
account) relation
\[
{}[E_{j-1,j+2},E_{j,j+1}]=0 \qquad\text{i.e.,~}
[[x_{j-1},[x_j,x_{j+1}]],x_j]=0.
\]
Here $j\in\overline{2, n-1}$.

\ssec{Proofs: Lie superalgebras} In the exceptional cases, the
relations are ob\-tained by means of \texttt{SuperLie}. For the
$\fsl$ series, the arguments of the non-super case are applicable.
For the other series, the answers are conjectural but we tested them
by means of \texttt{SuperLie} for small values of superdimensions,
and hence are sure.

\end{document}